\documentclass[11pt]{article}

\usepackage[T1]{fontenc}
\usepackage[latin1]{inputenc}
\usepackage{mathtools,amsthm,amssymb}
\usepackage{mathabx}
\usepackage{mathrsfs}

\usepackage{xcolor}
\usepackage{authblk}

\usepackage{hyperref}

\usepackage[margin=1in]{geometry}

\newtheorem{prop}{Proposition}[section]
\newtheorem{cor}[prop]{Corollary}
\newtheorem{defi}[prop]{Definition}

\newtheorem{lem}[prop]{Lemma}

\newtheorem{theo}[prop]{Theorem}

\newcommand{\vertiii}[1]{{\left\vert\kern-0.25ex\left\vert\kern-0.25ex\left\vert #1
    \right\vert\kern-0.25ex\right\vert\kern-0.25ex\right\vert}}

\newcommand{\DD}{\mathbb{D}}
\newcommand{\VV}{\mathbb{V}}

\newcommand{\La}{ {\cal L }}

\newcommand{\Ea}{ {\cal E }}

\newcommand{\Va}{ {\cal V }}
\newcommand{\Ua}{ {\cal U }}

\newcommand{\Qa}{ {\cal Q }}

\newcommand{\Xa}{ {\cal X }}

\newcommand{\Ta}{ {\cal T}}
\newcommand{\Ha}{ {\cal H }}

\newcommand{\Pa}{ {\cal P }}
\newcommand{\Za}{ {\cal Z }}
\newcommand{\Ya}{ {\cal Y }}
\newcommand{\Wa}{ {\cal W }}

\def \RR{\mathbb{R}}

\def \EE{\mathbb{E}}

\def \CC{\mathbb{C}}
\def \LL{\mathbb{L}}

\def \WW{\mathbb{W}}
\def \BB{\mathbb{B}}

\numberwithin{equation}{section}

\begin{document}

\title{On One-Dimensional Riccati Diffusions}

\author[$1$]{A.N. Bishop\thanks{A.N. Bishop is also an adjunct fellow at the Australian National University (ANU). He was supported by the Australian Research Council (ARC) via a Discovery Early Career Researcher Award (DE-120102873).}}
\author[$2$]{P. Del Moral\thanks{P. Del Moral was supported in part by funding from the Simons Foundation and the Centre de Recherches Mathématiques, through the Simons-CRM scholar-in-residence program.}}
\author[$3$]{K. Kamatani}
\author[$4$]{B. R\'emillard\thanks{B. R\'emillard was supported in part by funding from the Natural Sciences and Engineering Research Council of Canada.}}
\affil[$1$]{{\small University of Technology Sydney (UTS) and CSIRO, Australia}}
\affil[$2$]{{\small INRIA, Bordeaux Research Center, France and UNSW, Sydney, Australia}}
\affil[$3$]{{\small Graduate School of Engineering Science, Osaka University, Osaka, Japan}}
\affil[$4$]{{\small HEC Montr\'eal, Montr\'eal, Canada}}
\date{}

\maketitle

\begin{abstract}
This article is concerned with the fluctuation analysis and the stability properties of a class of one-dimensional Riccati diffusions. These one-dimensional stochastic differential equations exhibit a quadratic drift function and a non-Lipschitz continuous diffusion function. We present a novel approach, combining tangent process techniques, Feynman-Kac path integration, and exponential change of measures, to derive sharp exponential decays to equilibrium. We also provide uniform estimates with respect to the time horizon, quantifying with some precision the fluctuations of these diffusions around a limiting deterministic Riccati differential equation. These results provide a stronger and almost sure version of the conventional central limit theorem. We illustrate these results in the context of ensemble Kalman-Bucy filtering. To the best of our knowledge, the exponential stability and the fluctuation analysis developed in this work are the first results of this kind for this class of nonlinear diffusions.
\end{abstract}

\section{Introduction}

Let $W_t$ be a Wiener process, $A\in\RR$, $R\wedge S>0$, and $U\wedge V\geq 0$ some  given parameters.
We consider the diffusion process on the non-negative half-line $\RR_+=[0,\infty[$  defined for any $X_0\in\RR_+$ by the stochastic differential equation
\begin{eqnarray}
dX_t=\Lambda(X_t)~dt+\sigma_{\epsilon}(X_t)~dW_t, \label{mp-ref}
\end{eqnarray}
with some Riccati-type drift function
\begin{eqnarray}
\Lambda(x)=2Ax+R-Sx^2\quad \mbox{\rm and diffusion term}\quad \sigma_{\epsilon}(x):=\epsilon~\sqrt{x(U+Vx^2)}. \label{def-sigma-Lambda}
\end{eqnarray}
and some parameter $\epsilon\in\RR_+$. When $\epsilon=0$, $x_0\in\RR_+$, the diffusion (\ref{mp-ref}) reduces to the deterministic Riccati dynamical system
defined on $\RR_+$ by the equation
\begin{equation}
\partial_t \,x_t=\Lambda(x_t)=-S~(x_t-\varpi_+)~(x_t-\varpi_-),
\label{nonlinear-KB-Riccati}
\end{equation}
with the equilibrium states $(\varpi_-,\varpi_+)$ defined by
   \begin{equation}\label{def-varpi}
S\,\varpi_-:=A-\lambda/2~<~0~<~S\,\varpi_+:=A+\lambda/2 \quad \mbox{\rm with}\quad
\lambda :=2\sqrt{A^2+RS}.
 \end{equation}
We let $\Phi_{s,t}^{\epsilon}(X_s)=X_t$, with $s\leq t$, denote the stochastic flow associated with the diffusion defined by \eqref{mp-ref}. We also let $\phi_{s,t}(x_s)=x_t $ be the semigroup associated with the Riccati equation (\ref{nonlinear-KB-Riccati}). We may write $\Phi_{0,t}^{\epsilon}=\Phi_{t}^{\epsilon}$ and $\phi_{0,t}=\phi_{t}$, etc.

Under mild conditions discussed later, the Riccati diffusion $X_t$ is well defined on the half line $\RR_+$. In this case, the origin is a regular and repellent state in the sense that the process can start at $0$, but will never return to the origin.

Note the case $V = 0$ implying $\sigma_{\epsilon}(x):=\epsilon\,\sqrt{Ux}$ may act as a basic canonical prototype for a quadratic (Riccati-type) diffusion equation on the non-negative half-line.

The analysis of one-dimensional diffusions of the form (\ref{mp-ref}) acts as a basic prototype for the study of various quadratic Riccati-type diffusions arising in multivariate statistics, signal processing, and econometrics and financial mathematics. Other quadratic diffusion models, different from (\ref{mp-ref}) also appear in the literature. For example, backward-type matrix Riccati diffusions arise in linear-quadratic optimal control problems with random coefficients; see, e.g., \cite{bismut1976,Kohlmann2003,hu2003indefinite}. A different class of random Riccati equations arises in network control and filtering with random observation losses; see \cite{sinopoli2004kalman} and references therein. Note also that the Cox-Ingersoll-Ross process (i.e. the Wishart process in one-dimension \cite{Cox/Ingersoll/Ross:1985a}) can be viewed as special, linear, simplification of this model with $A<0<R\wedge U$ and no quadratic term $S=0=V$. This latter case illustrates that taking $A<0$ stable (i.e. Hurwitz stable) may significantly simplify the derivation of certain fluctuation and stability estimates and reduces the broader applicability of related results.

\subsection{A Related Diffusion Equation and Ensemble Kalman-Bucy Filtering}

The stochastic Riccati equations defined by \eqref{mp-ref} are motivated by applications in signal processing and data assimilation problems, and more particularly in the stochastic analysis of ensemble Kalman-Bucy-type \cite{evensen2003ensemble} filters (abbreviated {\tt EnKF}). In this context, up to a change of probability space, the matrix-valued version of the stochastic Riccati equation (\ref{mp-ref}) represents the evolution of the sample covariance associated with these filters. For more details, we refer to ~\cite{DelMoral/Tugaut:2016}, as well Section \ref{sec-statementEnKF} and Section~\ref{EnKf-sec} in the present article. The one-dimensional case (\ref{mp-ref}) represents the flow of the sample variance.

We are also interested in a certain stochastic Ornstein-Uhlenbeck process of the form,
\begin{equation}\label{def-intro-Z}
dZ_t=\displaystyle\frac{1}{2}~\partial \Lambda(X_t)~Z_t~dt+
\varsigma_{\overline{\epsilon}}(X_t)~dW^{\prime}_t\quad\mbox{\rm with}\quad
\varsigma^2_{\overline{\epsilon}}:=\varsigma^2+\sigma_{\overline{\epsilon}}^2\quad\mbox{\rm and}\quad
\varsigma^2(x):=R+Sx^2,
\end{equation}
with a Wiener process $W^{\prime}_t$ independent of $W_t$ and some parameter $\overline{\epsilon}\in[0,1]$. Note that (\ref{def-intro-Z}) is coupled to (\ref{mp-ref}) through both the drift and diffusion terms. We call this process a stochastic Ornstein-Uhlenbeck process because of this coupling; i.e. the coefficients of this Ornstein-Uhlenbeck process are themselves stochastic. The stochastic flow of the $(\RR_+\times\RR)$-valued diffusion  $(X_t,Z_t)$ defined by the stochastic differential equations (\ref{mp-ref}) and (\ref{def-intro-Z}) is denoted by
\begin{equation}
\Pi^{(\epsilon,\overline{\epsilon})}_t(x,z)=\left(\Phi^{\epsilon}_t(x),\Psi^{(\epsilon,\overline{\epsilon})}_t(x,z)\right). \label{flows-X-Z}
\end{equation}

In the context of the {\tt EnKF}, the parameters $\epsilon=2/\sqrt{N}$ and $\overline{\epsilon}={1}/{\sqrt{N+1}}$ are related to the population size of a system of $(N+1)$ interacting particles. Up to a change of probability space, the difference (error) between the {\tt EnKF} sample mean and the true signal state is described by a stochastic Ornstein-Uhlenbeck process of the form (\ref{def-intro-Z}).

When $\overline{\epsilon}=0=\epsilon$, equivalently when $(U,V)=0$ or $N\rightarrow\infty$,  the diffusion $Z_t:=\mathscr{Z}_t$ reduces to the difference (error) between the conventional Kalman-Bucy filter (see \cite{Bishop/DelMoral:2016}) and the true signal state. It is  given by the non-homogeneous Ornstein-Uhlenbeck process
  \begin{equation}\label{def-Z0}
d\mathscr{Z}_t=\displaystyle\frac{1}{2}~\partial \Lambda(x_t)~\mathscr{Z}_t~dt+
\varsigma(x_t)~dW^{\prime}_t
    \end{equation}
where $x_t$ denotes the solution of the deterministic Riccati equation (\ref{nonlinear-KB-Riccati}). For more details on the derivation and origin of these coupled diffusion processes, and their particle interpretations, we refer to Section \ref{sec-statementEnKF} and later Section~\ref{EnKf-sec}. Again, we reiterate that (\ref{mp-ref}) is of separate interest on its own, as a model for quite general quadratic (Riccati-type) diffusion equations.

\subsection{Objectives}

The main results and paper organisation are given in Section \ref{sec-mainresults}.

To the best of our knowledge, the time-uniform fluctuation and stability analysis of Riccati-type diffusion equations has not been addressed in the literature. This article addresses this problem in the one-dimensional setting. This diffusion may be associated with the {\tt EnKF} sample variance, and we accommodate the general case, allowing the underlying signal (defined by $A$) to be unstable. 

\subsubsection{Uniform Bias and Fluctuation Analysis}

The first objective of this article is to analyze the fluctuations of the pair of diffusion processes $(X_t,Z_t)$ when the parameters $\epsilon$ and $\overline{\epsilon}$ tend to $0$. We provide uniform fluctuation estimates w.r.t. the time parameter, as well as sharp exponential decay rates for the fluctuations, in the sense that they converge to the exponential decay of the Riccati equation (\ref{nonlinear-KB-Riccati}) towards its fixed point $\varpi_+$. See a statement of the main results in Section \ref{sec-statementfluctuation}.

A matrix-valued version of (\ref{mp-ref}) is studied in \cite{Bishop/DelMoral/Niclas:2017} where complete Taylor-type expansions with remainder are developed at any order, and non-asymptotic fluctuation, bias and central-limit theorems are also given. However, most of the uniform results in \cite{Bishop/DelMoral/Niclas:2017} only hold when $A$ is stable. In this article, first and second-order non-asymptotic Taylor-type expansions are derived. In contrast to \cite{Bishop/DelMoral/Niclas:2017}, the time-uniform estimates in this article do not require stability of $A$, i.e. $A<0$ in this case. The estimates here are valid for any value of $A\in\RR$, as soon as $R\wedge S>0$. However, unlike \cite{Bishop/DelMoral/Niclas:2017} we only consider a second-order expansion of the diffusion here. The analysis in \cite{Bishop/DelMoral/Niclas:2017} also extends to full path-wise fluctuation properties. See the later work also in \cite{2018arXiv180800235} which deals with matrix-valued Riccati diffusions and also accommodates unstable $A$ matrices.

Some general perturbation results for vector-valued stochastic differential equations are studied in \cite{hutzenthaler2014}. However, the analysis of \cite{hutzenthaler2014} is not focused on time-uniform estimates, but rather on a fluctuation analysis of perturbed stochastic processes under some generalised coefficient (drift/diffusion) conditions. Here we restrict our fluctuations analysis to quadratic stochastic differential equations of the Riccati-type (\ref{mp-ref}), and to the specified stochastic Ornstein-Uhlenbeck process (\ref{def-intro-Z}).

Previewing later results, we note that for any $n\geq 1$, any integrable initial state $\EE(X_0)<\infty$, and any $A\in\RR$ we have the uniform moment estimates
$$
(n-1)\,\frac{V}{S}\,\epsilon^2<2 ~~~~\Longrightarrow~~~~\sup_{t\geq 0}{\EE\left(X_t^n\right)}<\infty.
$$
The existence of these moments ensures there is no finite-time explosion, as well as the positive-recurrence of the diffusion (\ref{mp-ref}); e.g. see~\cite{Ikeda/Watanabe:1983}. Conversely, starting in the stationary regime associated with the invariant measures introduced in Section~\ref{ricc-diff-stab-sec}, we have
\begin{eqnarray*}
 (n-2)\,\frac{V}{S}\,\epsilon^{2}\geq2~~~~\Longrightarrow~~~~ \forall t\geq 0,~~
 \EE\left(X_t^n\right)=\infty.
 \end{eqnarray*}
Following~\cite{bakry1986critere}, we may also arrive at these moment properties using the exponential stability estimates given in Section~\ref{ricc-diff-stab-sec}.

\subsubsection{Contraction and Stability Analysis}

The second objective of this article is to analyse the stability properties of Riccati diffusions of the form (\ref{mp-ref}). See a statement of the main results in Section \ref{sec-statementstability}.

Since the pioneering articles of~\cite{feller1954diffusion,feller1954general}, the theory of one-dimensional diffusions has been developed in numerous directions, including using operator theory, spectral analysis, and classical semigroups and stochastic differential techniques. It is beyond the scope of this article to review these developments. We refer to the seminal book of~\cite{chen2006eigenvalues} and the references therein.

Several general conditions for the exponential decay to equilibrium of one-dimensional diffusions have been proposed in the literature. Some of them are based on the existence of a Poincar\'e or Hardy-type inequalities w.r.t the invariant measure, or related variational-type integral criteria w.r.t the scale and the speed measures; see for instance~\cite{chen2006eigenvalues}, as well as the articles~\cite{barthe2003sobolev,muckenhoupt1972hardy} and~\cite{cattiaux2010functional} for some particular classes of heavy tailed limiting distributions. Another more technical route is to estimate the transition densities of the process, and to find judicious Lyapunov functions as in~\cite{mattingly2002ergodicity}, or to use coupling and transport techniques as in~\cite{cheng2017exponential,gentil2005modified,gozlan2012transport}. Another strategy is to relate the spectrum of the diffusion with the one of the Schr\"odinger-type operators as in~\cite{shigekawa2013spectra}, or in~\cite{wu2001large}.

The sophisticated approaches discussed above are often not adapted for deriving precise estimates of the spectral gap of nonlinear diffusions with nonlinear and non-uniformly elliptic diffusions functions. In our case, the Riccati diffusions (\ref{mp-ref}) have a quadratic drift $\Lambda$ and the diffusion function $\sigma_1$ is not globally Lipschitz. In addition, the function $\sigma_1$ is not uniformly positive and the origin of the Riccati diffusion (\ref{mp-ref}) corresponds to a Neumann-type boundary condition. We note that the drift function $\Lambda=\partial F$ is the derivative of the double-well drift function
$$
	F(x)=-\frac{S}{3}~x~(x-\chi_{-})~(x- \chi_+),
$$
with roots $(\chi_{-},\chi_{+})$ given by
$$
\chi_{-}:=\frac{3A}{2S}-\left[
\left(\frac{3A}{2S}\right)^2+\frac{3R}{S}\right]^{1/2}< 0< \chi_+:=\frac{3A}{2S}+\left[
\left(\frac{3A}{2S}\right)^2+
\frac{3R}{S}\right]^{1/2}.
$$
This latter point implies that $F$ is not convex and the diffusion function is not constant, and hence the conventional tools of the Bochner-Bakry-Emery theory~\cite{bakry1983diffusions} cannot be applied since these typically require that both $-\partial ^2F$ and $\sigma_1$ are uniformly lower bounded.

Also notice that the diffusion has a polynomial growth of order $3/2$ as soon as $V>0$. It follows by \cite{hutzenthaler2011strong} that a basic Euler time-discretization may blow up, regardless of the boundedness properties of the diffusion.  Although the Euler schemes may blow up, it is quite possible that their tamed versions converge. In this one-dimensional case, one may also look at exact simulation methods of the kind discussed in \cite{Beskos2005}.

In the present article we develop a new approach for studying the stability of the equation (\ref{mp-ref}), combining tangent processes with Feynman-Kac path integrals and using a judicious exponential change of probability measure. We obtain sharp estimates of exponential decays of Riccati diffusions for small values of the parameter $\epsilon$.

\subsubsection{Applications to Ensemble Kalman Filters}

The {\tt EnKF} can be interpreted as the mean-field particle approximation of a nonlinear McKean-Vlasov-type diffusion. These probabilistic models were introduced in~\cite{mckean1966class}. For a detailed discussion on these models and their  domains of application, we refer the reader to the lecture notes of~\cite{Sznitman:1991,meleard1996asymptotic}, and the monograph~\cite{del2013mean}. The {\tt EnKF} is a key numerical methods for solving high-dimensional forecasting and data assimilation problems; see, e.g., \cite{evensen2003ensemble}. We refer to (some of) the seminal methodology papers in \cite{evensen1994orig,houtekamer1998data,burgers1998analysis,andersonanderson1999,houtekamer2001seq,hamill01,anderson2001ensemble,anderson2003local,tippetsqrt2003,sakov2008deterministic,Reich2013,Yang2016}. This list is by no means exhaustive; see also \cite{evensen2009book,kalnay2003atmospheric,law2015data,Reich2015book} for more background, and the detailed chronological list of references in Evensen's text \cite{evensen2009book}. Our analysis here captures a one-dimensional, rather toy (i.e. linear-Gaussian), {\tt EnKF} setup.

Convergence and large-sample asymptotics of the discrete-time {\tt EnKF} has been studied in~\cite{le2011large}, and in~\cite{mandel2011convergence}, in the sense of taking the number of particles $(N+1)\rightarrow\infty$. Non-linear state-space models are accommodated in this sense in \cite{le2011large}. In~\cite{law2016deterministic}, the authors extend this idea to continuous-time non-Gaussian state-space models (e.g. certain non-linear diffusions). In \cite{tong2016stability}, the authors analyse the long-time behaviour of the {\tt EnKF}, with finite ensemble size, using Foster-Lyapunov techniques. Applying these results to basic linear-Gaussian filtering problems, the analysis in \cite{tong2016stability} would require stability of the signal model, i.e. stability of $A$. In \cite{Kelly2014}, again the long-time behaviour of the {\tt EnKF} is analysed (with and without so-called variance inflation) under a class of quadratic dissipative system models; and which again if linearised equates to a form of stability on the signal model. In \cite{Kelly2014} time-uniform results follow under a sufficiently large inflation regime. See also~\cite{majda-2,tong-2} for related stability analysis in the presence of adaptive covariance inflation and projection techniques.

Ensemble Kalman methods for inverse problems have also been considered in the literature \cite{Iglesias2013} with some related analysis \cite{2017a,2017b}. See these references for further details on this topic.

In the linear-Gaussian filtering domain specifically, uniform error estimates w.r.t. the time horizon have been developed in~\cite{DelMoral/Tugaut:2016}; for any ensemble size, without inflation, but once again under the similarly strong assumption that the signal is stable; i.e. $A<0$ in the one-dimensional case. And as noted previously, non-asymptotic fluctuation, bias and central-limit theorems on the {\tt EnKF} sample covariance are given in \cite{Bishop/DelMoral/Niclas:2017}; with time-uniform results holding under a stability condition on $A$. So-called variance inflation and localization are considered from a purely mathematical vantage in the linear-Gaussian setting in \cite{Bishop/DelMoral/Pathiraja:2017}. We refer further to \cite{DelMoral/Tugaut:2016,Bishop/DelMoral/Pathiraja:2017,Bishop/DelMoral/Niclas:2017} for some related and further background and references on the mathematical properties of the {\tt EnKF} in the linear-Gaussian setting. See also \cite{2018arXiv180800235} for a matrix-valued extension of the work considered herein, which seeks to remove the strong stability results on the signal required in prior work.

Given (\ref{mp-ref}), (\ref{def-intro-Z}) and the preceding discussion, we note the convergence properties of the {\tt EnKF} rely heavily on the fluctuations of a Riccati diffusion. In this article, we study the fluctuation and stability of the one-dimensional model (\ref{mp-ref}) as the time horizon $t$ tends to $\infty$, and as $\epsilon$ tends to $0$. {\em We make no assumption on the stability of the underlying signal; i.e. on the stability of $A\in\RR$}. We emphasize that the Riccati diffusion is reversible w.r.t. a probability measure for any value of $A$ as soon as $R\wedge S>0$. However, the nature of these invariant measures depends strongly on the choice of the parameters $(U,V)$. We refer to Section~\ref{inv-meas-sec} for a detailed exposition of these measures, and for an explicit description of the different classes of stationary measure captured by different choices of $(U,V)$. Different choices of $(U,V)$ correspond to different variants of the {\tt EnKF}.

Of course, we note that the one-dimensional models discussed in this article do not capture faithfully the higher dimensional problems typically considered in the filtering and data assimilation literature. Nevertheless, this analysis provides some insight on the fluctuation and the stability properties of the {\tt EnKF} when the signal itself is unstable (i.e. $A>0$ in this case). The stability properties of matrix-valued Riccati diffusions with unstable signal models has been considered in the more recent article~\cite{2018arXiv180800235}. The analysis developed in \cite{2018arXiv180800235} extends the one-dimensional analysis here to multidimensional data assimilation problems and ensemble Kalman-Bucy filters with unstable signals. However, in the multidimensional setting, the decay rates to equilibrium are not sharp, and the stationary measures are not given in closed form.

\subsubsection{Remarks}

To the best of our knowledge, the uniform fluctuation/moment estimates and the exponential rates to equilibrium discussed in this article are the first results of this type for this class of models.

Moreover, while applications in data assimilation and the {\tt EnKF} are partial motivators for this work, {\em the one-dimensional diffusion (\ref{mp-ref}) is also of interest in its own mathematical right, as a prototype for quadratic (Riccati-type) stochastic differential equations}.

\subsection{Some Preliminary Notation}

Here we introduce some notation necessary for the statement of our main results which are stated  in the next section. Firstly, given a probability measure $\pi$ on $\RR$ and some function $f\in \LL_1(\pi)$, we write
$\pi(f)$ the Lebesgue integral
$$
\pi(f):=\int~f(x)~\pi(dx)
\quad\mbox{\rm and we let $\theta$ be the identity function $\theta(x):=x$.}
$$
Throughout this article $\partial^nf$ stands for the $n$-th derivative of some smooth function $x\mapsto f(x)$ w.r.t. the parameter $x$. When $n=1$ we write  $\partial f$ instead of $\partial^1f$.

Given some real valued stochastic process $X_t$,  whenever they exist for any $n\geq 1$ and any time horizon $t\geq 0$,
we set
$$
\vert X\vert:=\sup_{t\geq 0}\vert X_t\vert\qquad
\vertiii{X_t}_{n}:=\EE\left[\Vert X_t\Vert^n\right]^{1/n}\quad\mbox{\rm and}\quad
\vertiii{X}_{n}:=\sup_{t\geq 0}\vertiii{X_t}_{n}.
$$

The $\sigma$-distance on $\RR_+$ and the corresponding Wasserstein distance (a.k.a. Kantorovich-Monge-Rubinstein metric) is
defined by the formula
$$
d_{\sigma}(x,y):=\left\vert\int_x^y~\sigma^{-1}(z)~dz\right\vert\quad\mbox{\rm and}\quad
\DD_{\sigma}\left(\mu_1,\mu_2\right):=\inf{\left\{\int~d_{\sigma}(x_1,x_2)~\mu(d(x_1,x_2))\right\}}.
$$
In the above display, $\sigma$ stands for some positive function on $\RR_+$, and the infimum is taken over all coupling probability measures $\mu$ on $\RR_+^2$ with the first and second marginals equal to $\mu_1$ and $\mu_2$.
The variational form of the distance $d_{\sigma}$ is given by the formula
$$
d_{\sigma}(x,y)=\sup{\left\{\vert f(x)-f(y)\vert~f\in C^{\infty}(\RR_+)~~\mbox{\rm s.t.}~~
\Vert \sigma \partial f\Vert\leq 1\right\}}.
$$
In the above display $\Vert f\Vert=\sup_{x\geq 0}\vert f(x)\vert$ stands for the uniform norm.
In the same vein, we have the Kantorovich-Rubinstein duality relation:
\begin{eqnarray}\label{k-r-duality}
\DD_{\sigma}\left(\mu_1,\mu_2\right)&=&\sup{\left\{\vert \mu_1(f)-\mu_2(f)\vert~:~f\in \mbox{\rm Lip}_{\sigma}(\RR_+)\right\}}.
\end{eqnarray}
In the above display $ \mbox{\rm Lip}_{\sigma}(\RR_+)$ stands for the space of Lipschitz functions $f$ on the metric space $(\RR_+,d_{\sigma})$ with unit Lipschitz constant; that is such that
$$
\vert f(x)-f(y)\vert\leq d_{\sigma}(x,y).
$$
The Kantorovich-Rubinstein theorem (\ref{k-r-duality}) on arbitrary compact metric spaces has been presented in~\cite{kant-1,kant-2}. The extension to separable metric spaces can be found in~\cite{acosta,dudley,dudley-2}, see also~\cite{kellerer-2} for an extension to general metric spaces and Radon measures. Further details on optimal transport and Wasserstein distance can be found in the books~\cite{Villani-1,Villani-2}.

The article discusses several $\LL_n$-mean error fluctuation
estimates as well as a series of exponential asymptotic stability inequalities.  Special attention is paid to the quantitative nature
of the fluctuation and stability results. We track closely the dependency of the
estimation constants and the convergence decays in terms of the parameters of the model.

We shall distinguish two classes of parameters. Firstly, we introduce parameters that depend solely on the kinetic parameters of the drift function $(A,R,S)$. Secondly, we introduce parameters that also depend on the choice of the diffusion parameters $(U,V)$.

$\bullet$ The first set of parameters depending on the drift parameters $(A,R,S)$ is the collection of parameters $( \imath,\jmath,\imath_{\kappa},\jmath_{\kappa})$  indexed by $\kappa\geq 0$  and defined by
\begin{equation}\label{def-iota-iota-kappa}
\begin{array}{rclcccl}
  \imath &:=&\displaystyle\frac{\jmath}{\sqrt{1+\jmath^2}}&&
\jmath&:=&\displaystyle\frac{A}{\sqrt{RS}}\qquad
\imath_{\kappa}:=\kappa(1+\imath )\quad\mbox{\rm and}\quad\displaystyle
\frac{ \jmath_{\kappa}+1}{\imath_{\kappa}+1}:= (1+\imath)^2~\left(1+\jmath^2\right) \end{array}.
\end{equation}
The above parameters can alternatively be defined in terms of the  equilibrium states $(\varpi_-,\varpi_+)$
introduced in  (\ref{def-varpi}) using the formulae
\begin{equation}\label{definition-varpi}
\varpi:=1-\varpi_+/\varpi_-=2(1-\imath)^{-1} ~~\Longleftrightarrow~~  \imath=1-{2}/{\varpi}.
\end{equation}
Observe the uniform estimates w.r.t. the model parameters
$$
 -1\leq \imath\leq 1\qquad  0\leq \frac{\jmath_{\kappa}+1}{\jmath^2+1}\leq 4(1+2\kappa)\quad\mbox{\rm as well as}\quad
0\leq \imath_{\kappa}\leq 2\kappa.
$$

$\bullet$ Define $\overline{U}:=U/R$ and $\overline{V}:=V/S$. We also consider the collection of positive parameters $(\zeta,\zeta_{\kappa})$ dependent on the diffusion parameter $(U,V)$ (as well as on $(A,R,S)$) and defined by the formulae
\begin{equation}\label{def-zetas}
\zeta:=\frac{\imath+1}{\jmath^2+1}~\overline{U}
\quad\mbox{\rm
and}
\quad  \zeta_{\kappa}:=~
\frac{\imath_{\kappa}+1}{\jmath^2+1}~\overline{U}~+~
\frac{ \jmath_{\kappa}+1}{\jmath^2+1}~\overline{V}.
\end{equation}
Observe the uniform estimates w.r.t. the drift parameter $A$ given by
$$
0\leq \zeta \leq~2~\overline{U}\quad\mbox{\rm with}\quad 0\leq \zeta_{\kappa}\leq (2\kappa+1)\left(\overline{U}+4~ \overline{V}\right).
$$

As mentioned above, these parameters allow us to quantify with some precision the constants arising in the fluctuation and the stability inequalities in terms of the model parameters. For instance, when $A=0$ the above parameters resume to
$$
 \imath= \jmath=0\quad\mbox{\rm  and}\quad  \imath_{\kappa}=\jmath_{\kappa}=\kappa~~\Longrightarrow~~ \varpi=2\quad\mbox{\rm  and}\quad
\zeta=\overline{U}\quad\mbox{\rm  and}\quad  \zeta_{\kappa}=(\kappa+1)~(\overline{U}+\overline{V}).
$$
Also observe that
\begin{equation}\label{def-Di}
\widehat{\sigma}(x):=x^{ \imath+1} ~~\Longrightarrow~~
d_{\widehat{\sigma}}(x,y):=\vert  \imath\vert^{-1}~\vert x^{-{ \imath}}-y^{-{ \imath}}\vert,
\end{equation}
with the convention $ d_{\widehat{\sigma}}(x_1,x_2)=\vert \log{x_1}-\log{x_2}\vert$ when $\imath=0$.

Finally we write $c,c_n,c_{i,n},c_{i,n}(x),\ldots$
  some positive universal constants and parameters whose values may vary from line to line, but they only depend on some parameters $i,n,x$, etc, as well as on the parameters of the Riccati processes $(A,R,S,U,V)$, but importantly not on the time horizon.

\section{Statement of the Main Results and Article Organisation}\label{sec-mainresults}

\subsection{A Uniform Bias and Fluctuation Theorem} \label{sec-statementfluctuation}

The first main objective of this article is to quantify the bias and the fluctuations of the diffusion process $X_t$ in (\ref{mp-ref}) around the limiting Riccati equation $x_t$ in (\ref{nonlinear-KB-Riccati}) as $\epsilon\rightarrow 0$. These regularity properties are also used to analyze the fluctuations of the stochastic Ornstein-Uhlenbeck process $Z_t$ in (\ref{def-intro-Z}) around the limiting diffusion $\mathscr{Z}_t$ introduced in (\ref{def-Z0}). In both situations, we provide a series of refined uniform $\LL_n$-type estimates w.r.t. the time horizon.

Our first main result takes basically the following form.

\begin{theo}\label{t1-intro}
The $n$-th moments of the stochastic flow $\Phi^{\epsilon}(x)$ are uniformly bounded as soon as $n\geq 1$ and $\epsilon\geq 0$ are chosen so that $(n-1)\,\epsilon^2\,\overline{V}<2$.
In this situation, for any $x\in\RR_+$ we have the uniform estimates
$$
\vertiii{ \Phi^{\epsilon}(x)-\phi(x)}_n\vee \vertiii{ Z-\mathscr{Z}}_n~\leq~ c_{1,n}(x)~\epsilon.
$$
In addition we have the uniform fluctuation and bias estimates
$$
\vertiii{ \Phi^{\epsilon}(x)-\phi(x)-\epsilon~\VV(x)}_n\leq c_{2,n}(x)~\epsilon^{2}\quad
\mbox{and}\quad \vert \EE\left(\Phi^{\epsilon}(x)\right)-\phi(x)-\epsilon^2~\WW(x)\vert\leq c_{3,n}(x)~\epsilon^3,
$$
for some processes $(\VV_t(x),\WW_t(x))$ and some finite functions $c_{i,n}(x)$ which are explicitly defined in terms of the model parameters.
\end{theo}

For a more precise statement and a detailed description of   $(\VV_t(x),\WW_t(x))$ and  $c_{i,n}(x)$ we refer the reader to Section~\ref{fluctuation-sec}
and Section~\ref{2-d-section}; see for instance Theorem~\ref{theo--1intro} and Corollary~\ref{cor-1-EnKF}.

In contrast with the multivariate stochastic analysis developed in~\cite{Bishop/DelMoral/Niclas:2017}, we emphasize that the uniform estimates w.r.t. the time horizon in the one-dimensional models considered in this article do not require stability of $A$.

  The stochastic analysis developed in the present work relies on specific properties of one-dimensional diffusions. Section~\ref{properties-Riccati-section} provides a brief review on the regularity properties of Riccati semigroups. Their smoothness and the continuity properties are discussed in Section~\ref{smoothness-sec} and in Section~\ref{robustness-sec}. The stochastic analysis developed in this article combines Riccati semigroup techniques with fluctuating random fields methods. This methodology, as well as a precise description of the bias and fluctuation estimates of stochastic Riccati diffusions are described in Section~\ref{fluctuation-sec}. The analysis of stochastic Ornstein-Uhlenbeck processes is discussed in Section~\ref{2-d-section}

\subsection{Contraction and Stability Theorems} \label{sec-statementstability}

 The second objective in this article is to analyze the stability properties of the Riccati diffusion $X_t$ in (\ref{mp-ref}) as $t\rightarrow\infty$. We denote by $P_t^{\epsilon}$ be the semigroup associated with the stochastic flow $\Phi_{t}^{\epsilon}$; that is
 for any bounded measurable function $f$ and any probability measure $\mu$ on $\RR_+$  we have
\begin{equation}\label{def-sg}
P_t^{\epsilon}(f)(x):=\EE\left(f(\Phi_{t}^{\epsilon}(x))\right)\quad\mbox{\rm and}\quad \mu=\mbox{\rm Law}(X_0)~\Longrightarrow ~\mu P^{\epsilon}_t:=\mbox{\rm Law}(X_t).
\end{equation}

The stability of the Riccati diffusion (\ref{mp-ref}) is quite generally considered in Section \ref{ricc-diff-stab-sec}. In Section~\ref{2-d-section}, these stability properties are used to analyze the long-time behaviour of the stochastic Ornstein-Uhlenbeck process (\ref{def-intro-Z}).

   The description of the reversible measures $\pi_{\epsilon}$ of the semigroup $P_t^{\epsilon}$ are discussed in some details in Section~\ref{inv-meas-sec}. We preview that discussion and note that $\pi_{\epsilon}$ is a heavy-tailed distribution whenever $V>0$. When $V=0$, the stationary measure $\pi_{\epsilon}$ is a weighted Gaussian distribution restricted to the half line. This means that we must carefully account for the values of the parameters $(U,V)$ in the diffusion function (\ref{def-sigma-Lambda}) in every estimate.

To describe precisely our main results we need to introduce some terminology. The stability of the
linear process (\ref{def-intro-Z}) is dictated by the stability of the exponential semigroups defined by
\begin{equation}\label{def-Ea-st}
\Ea_{s,t}^{\epsilon}(x):=\exp{\left[\int_s^t\frac{1}{2}~\partial \Lambda(\Phi^{\epsilon}_u(x))du\right]}\quad\mbox{\rm and}\quad
\Ea_{s,t}(x):=\exp{\left[\int_s^t\frac{1}{2}~\partial \Lambda(\phi_u(x))du\right]}.
\end{equation}
When $s=0$, we simplify notation and we write $(\Ea_{t}^{\epsilon}(x),\Ea_{t}(x))$ instead of $(\Ea_{0,t}^{\epsilon}(x),\Ea_{0,t}(x))$.

We also recall that the decay rate to equilibrium of $x_t$ is given for any time horizon $t\geq \upsilon>0$  by
$$
c_{1,\upsilon}~\exp{\left[-\lambda \, t\right]}~\leq~
\sup_{x\geq 0}{\vert \phi_t(x)-\varpi_{+}\vert}~\vee~ \sup_{x\geq 0}{\Ea_t(x)^2}~\leq~c_{2,\upsilon}\,\exp{\left[-\lambda \, t\right]}.
$$
The proof of this assertion is provided in Section~\ref{smoothness-sec}; see also the explicit formula (\ref{closed-form}).

Thus we can expect that
the Riccati diffusion (\ref{mp-ref}) tends to equilibrium with an exponential rate that converges towards $\lambda $ as the parameter $\epsilon\rightarrow 0$.
We consider the parameters
\begin{equation}\label{def-lambda-epsilon-iota}
 \widehat{\lambda}_{\epsilon} :=
\lambda \left(1-\frac{\epsilon^2}{2}~\zeta\right)
\quad\mbox{\rm and}\quad
\widehat{\lambda}_{\epsilon,\kappa} := \lambda ~\left(1-\frac{\epsilon^2}{2}~\zeta_{\kappa}\right).
\end{equation}
where it is implicitly assumed that $\epsilon\in\RR_+$ is chosen s.t. $\widehat{\lambda}_{\epsilon,\kappa} \geq 0$. In this case, we have
$$
\kappa\geq 1\quad\Longrightarrow\quad
\zeta\leq \zeta_{\kappa}\quad\Longrightarrow\quad
 \widehat{\lambda}_{\epsilon} \,\geq \widehat{\lambda}_{\epsilon,\kappa}.
 $$
With this notation, our second main result takes basically the following form.
\begin{theo}\label{t2-intro}
There exists some parameter ${\epsilon_{\star}}\in\RR_+$ such that
for any time horizon $t> 0$, any $\epsilon\in [0,{\epsilon_{\star}}\,]$, and any probability measures $\mu_1,\mu_2$ on $\RR_+$ we have
the contraction inequality
\begin{equation}\label{Wasserstein-estimate-overline}
\DD_{\widehat{\sigma}}\left(\mu_1P^{\epsilon}_t,\mu_2P^{\epsilon}_t\right) ~\leq ~\exp{\left[-\widehat{\lambda}_{\epsilon} \,t\right]}~ \DD_{\widehat{\sigma}}\left(\mu_1,\mu_2\right).
\end{equation}
In addition, for any $x,x_1,x_2\in\RR_+$ we have
$$
\displaystyle{
\vertiii{\Ea_{t}^{\epsilon}(x)^{2}}_n}
\leq c_{1,n}(x)\,
\exp{\left[-\widehat{\lambda}_{\epsilon,n} \, t\right]}
$$
and
$$
\vertiii{\Phi^{\epsilon}_t(x_1)-\Phi^{\epsilon}_t(x_2)}_n\leq ~c_{2,n}~d_{\widehat{\sigma}}(x_1,x_2)\,\exp{\left[-\widehat{\lambda}_{\epsilon} \,t\right]}
$$
and for any $z=(z_1,z_2)\in (\RR_+\times\RR)^2$ we have
$$
\vertiii{\Psi^{(\epsilon,\overline{\epsilon})}_t(z_1)-\Psi^{(\epsilon,\overline{\epsilon})}_t(z_2)}_{n/2}^2\leq  c_{3,n}(z)~\exp{\left[-\widehat{\lambda}_{\epsilon,n} ~t\right]}.
$$
\end{theo}

For a more precise statement of these results, and a description of the parameters ${\epsilon}_{\star}$,
$c_{1,n}(x)$, $c_{2,n}$ and $c_{3,n}(z)$, we refer the reader to Section~\ref{contract-ineq-sec}, Section \ref{exp-semigroups-sec} and Section~\ref{2-d-section}; see for instance Theorem~\ref{theo-3-intro}, Theorem~\ref{theo-expo-sg-estimate},  Theorem~\ref{theo-4-intro}, and Theorem~\ref{theo-OU-proc}.

Poincar\'e inequalities and contraction estimates w.r.t. the Wasserstein distance associated with the diffusion function $\sigma_1$
can also be derived under more restrictive conditions. For example, we have the following result.
\begin{theo}\label{theo-2-intro}
Assume that $V=0$ and $0\leq\epsilon^2\,\overline{U}\leq 2$ and $A\in\mathbb{R}$ is chosen such that
$$
\lambda_{\epsilon}\,=-A+\sqrt{3RS\left(1-\frac{\epsilon^2}{2}~\overline{U}\right)} ~>~0.
$$
In this case, for any probability measures $\mu_1,\mu_2$ on $\RR_+$ we have the contraction inequality
\begin{equation}\label{Wasserstein-estimate}
\DD_{\sigma_1}\left(\mu_1P^{\epsilon}_t,\mu_2P^{\epsilon}_t\right)\leq \exp{\left[-\lambda_{\epsilon}\,t\right]}~
\DD_{\sigma_1}\left(\mu_1,\mu_2\right)
\end{equation}
In addition, the reversible measure $\pi_{\epsilon}$ of  the semigroup $P^{\epsilon}_t$  satisfies the Poincar\'e inequality
\begin{equation}\label{poincare-inequality}
\begin{array}{l}
\displaystyle 2\lambda_{\epsilon}\,\left[\pi_{\epsilon}(f^2)-\pi_{\epsilon}(f)^2\right]\leq~\pi_{\epsilon}\left(\sigma_{\epsilon}^2~(\partial f)^2\right)\\
\\
\displaystyle\Longrightarrow~~
\pi_{\epsilon}\left(\left[P_t^{\epsilon}(f)-\pi_{\epsilon}(f)\right]^2\right)^{1/2}\leq  \exp{\left[-\lambda_{\epsilon}\,t\right]}~\pi_{\epsilon}\left(\left[f-\pi_{\epsilon}(f)\right]^2\right)^{1/2}.
\end{array}
\end{equation}
\end{theo}
We refer to Lemma~\ref{lem-3} and Theorem~\ref{theo-2-inside} in Section~\ref{contract-ineq-sec} for a more detailed exposition of this result with various combinations of $(U,V)$. For example, other estimates of $\lambda_{\epsilon}$ are given in Lemma~\ref{lem-3} in Section~\ref{contract-ineq-sec} for other combinations of the parameters $(U,V)$.

Surprisingly, when $V=0$, and for negative values of $A$, the parameter $\lambda_{\epsilon}$ is greater than the decay rate $\lambda$ of the deterministic Riccati semigroup. That is, the Riccati diffusion converges faster to the invariant measure than the deterministic Riccati does to its fixed point. More specifically, under the assumptions of the above theorem, when $U=R$ we have
$$
 d_{\sigma_1}(x_1,x_2)=\frac{2}{\sqrt{R}}~\vert \sqrt{x_1}-\sqrt{x_2}\vert\quad\mbox{\rm and}\quad
\frac{\lambda_{\epsilon}}{\lambda}=\frac{\vert\jmath\vert+\sqrt{3(1-\epsilon^2/2)}}{\sqrt{1+\jmath^2}}>1.
$$
In this case, a closer inspection shows that the invariant measure $\pi_{\epsilon}$ with these parameters has lighter-weighted Gaussian-type tails; see (\ref{iota-pi-2}) in Section~\ref{inv-meas-sec}.

In a different direction, whenever $U=0$, Theorem~\ref{theo-2-inside} shows that the parameter $\lambda_{\epsilon}$ is larger than
the decay rate  $\lambda$ of the Riccati semigroup for any positive values of $A$. For instance, using the estimates stated in Lemma~\ref{lem-3} when $A\geq 0=U$ and $V=S$  we have
$$
 d_{\sigma_1}(x_1,x_2)=\frac{2}{\sqrt{S}}~\left\vert \frac{1}{\sqrt{x_1}}-\frac{1}{\sqrt{x_2}}\right\vert\quad\mbox{\rm and}\quad
\frac{\lambda_{\epsilon}}{\lambda}=\frac{\jmath+\sqrt{3(1-3\epsilon^2/2)}}{\sqrt{1+\jmath^2}}>1.
$$

Finally, we point to the later work also in \cite{2018arXiv180800235}, for a matrix-valued (partial) extension of the stability analysis considered herein, and which also accommodates unstable $A$ matrices.

\subsection{Applications to Ensemble Kalman Filters} \label{sec-statementEnKF}

Lastly, Section~\ref{EnKf-sec} is dedicated to the illustration of these results in the context of ensemble Kalman-Bucy filters ({\tt EnKF}) which are of interest in data assimilation problems. To underline the impact of our results, we preview several corollaries  which can be derived as direct consequences of the main theorems stated above.

  Consider a time-invariant linear-Gaussian filtering model of the following form
\begin{equation}\label{lin-Gaussian-diffusion-filtering}
d\mathscr{X}_t=A\,\mathscr{X}_t~dt+R^{1/2}\,d\mathscr{W}_t\quad\mbox{\rm and}\quad
d\mathscr{Y}_t=B\,\mathscr{X}_t~dt+\Sigma^{1/2}\,d\mathscr{V}_{t},
\end{equation}
where $(\mathscr{W}_t,\mathscr{V}_t)$ is an $2$-dimensional Brownian motion, $\mathscr{X}_0$ is a Gaussian random variable with mean and variance $(\EE(\mathscr{X}_0),P_0)$
(independent of $(\mathscr{W}_t,\mathscr{V}_t)$), $\Sigma,R> 0$, $A,B\in \RR$, $\mathscr{Y}_0=0$ and we set $S=B^2/\Sigma$.
We let $\Ya_t=\sigma\left(\mathscr{Y}_s,~s\leq t\right)$ be the $\sigma$-algebra filtration generated by the observations.
The conditional distribution $\eta_t=\mbox{\rm Law}\left(\mathscr{X}_t~|~\Ya_t\right)$ of the signal internal states $\mathscr{X}_t$ given $\Ya_t$ is a Gaussian distribution with a conditional mean
and a conditional variance  given by
$$
\mathscr{M}_t:=\EE\left(\mathscr{X}_t~|~\Ya_t\right)\quad \mbox{\rm and}\quad
\mathscr{P}_t:=\EE\left(\left[\mathscr{X}_t-\EE\left(\mathscr{X}_t~|~\Ya_t\right)\right]^2\right).
$$

Ensemble Kalman-Bucy filters can be interpreted as a (non-unique) mean field particle approximation of the Kalman-Bucy filtering equation. We refer to Section~\ref{EnKf-sec} for a more description of the sample mean $\widehat{\mathscr{M}}_t$ and the sample variance $\widehat{\mathscr{P}}_t$ associated with the three different versions of the {\tt EnKF} discussed in the present article. We consider also the ``re-centered'' process $\widehat{\mathscr{Z}}_t:=(\widehat{\mathscr{M}}_t-\mathscr{X}_t)$ that measures the difference between the sample mean  and the true signal state.

As shown in Sections~\ref{McKean-Vlasov-sec} and \ref{mean-field-EnKF-sec},  the evolution equations of the three different versions of $(\widehat{\mathscr{P}}_t,\widehat{\mathscr{Z}}_t)$ coincide with the ones of $(X_t,Z_t)$ defined in (\ref{mp-ref}) and in (\ref{def-intro-Z}) with the three combinations of parameters $(U,V)\in\left\{(R,S),(R,0),(0,0)\right\}$. Each of these three combinations corresponds to a particular instance of the {\tt EnKF} algorithm; see Section \ref{McKean-Vlasov-sec}.

We begin with uniform fluctuation estimates w.r.t. the time horizon.
\begin{cor}
 Suppose that $\widehat{\mathscr{M}}_0=\mathscr{M}_0$ and $\widehat{\mathscr{P}}_0=\mathscr{P}_0$. For any $n\geq 1$ and for sufficiently large $N\geq1$ we have the uniform estimates,
 $$
\sqrt{N}~{\vertiii{ \widehat{\mathscr{P}}-\mathscr{P}}_n}\leq c_{1,n}\quad \mbox{and}\quad
\sqrt{N}~{\vertiii{\widehat{\mathscr{M}}-\mathscr{M}}_n}\leq c_{2,n}.
 $$
 \end{cor}
 The l.h.s. assertion in the above corollary is a direct consequence of Theorem~\ref{t1-intro}. The r.h.s. estimate is given in Corollary~\ref{cor-ref-int}.

To give a flavour of our stability results in this context, the next corollary of Theorem~\ref{t2-intro} concerns the stability and the fluctuation of the {\tt EnKF} associated with the parameters $(U,V)=(R,0)$. In this case, the parameters introduced in (\ref{def-lambda-epsilon-iota}) resume to
$$
\zeta:=\frac{\imath+1}{\jmath^2+1}\leq2
\quad\mbox{\rm
and}
\quad  \zeta_{\kappa}:=
\frac{\imath_{\kappa}+1}{\jmath^2+1}~
\leq 2\kappa+1,
$$
with the parameters $( \imath,\jmath,\imath_{\kappa}) $ defined in (\ref{def-iota-iota-kappa}).
We denote by $(\widehat{\mathscr{P}}_t,\widehat{\mathscr{M}}_t)$ and $(\widehat{\mathscr{P}}_t^{\,\prime},\widehat{\mathscr{M}}^{\,\prime}_t)$ the processes starting from two possibly different initial conditions.

\begin{cor}
 For any $n\geq 1$ and any $t> 0$ we have
   \begin{equation}\label{enkf-intro-2}
N>
2\zeta_{2n}~~\Longrightarrow~~ \vertiii{\widehat{\mathscr{M}}_t-\widehat{\mathscr{M}}_t^{\,\prime}}_n~\leq ~c_{1,n}\,
\exp{\left[-\lambda ~\left(\frac{1}{2}-\frac{\zeta_{2n}}{N}\right)~t\right]},
 \end{equation}
and
 \begin{equation}\label{enkf-intro-1}
 N> 2\zeta~~\Longrightarrow~~ \vertiii{\widehat{\mathscr{P}}_t-\widehat{\mathscr{P}}_t^{\,\prime}}_n~\leq ~c_{2,n}\,
\exp{\left[
-\lambda ~\left(1-\frac{2\zeta}{N}\right)~t\right]}.
 \end{equation}
 \end{cor}
The first assertion is a consequence of Corollary~\ref{enkf-intro-2-proof}. The second assertion is a direct consequence of Corollary~\ref{enkf-intro-1-proof}. The other cases of $(U,V)\in\left\{(R,S),(R,0),(0,0)\right\}$ are considered in Section \ref{stability-EnKF-sec}.

The exponential decay of the exponential semigroup (\ref{def-Ea-st}) discussed in Theorem~\ref{t2-intro}
play a central role in the stability of the pair process $(\widehat{\mathscr{P}}_t,\widehat{\mathscr{Z}}_t)$. For large time horizons the Lyapunov exponent of the stochastic
Ornstein-Uhlenbeck process (\ref{def-intro-Z}) can be estimated by the formula
\begin{equation}\label{Lyap-expo-def}
\frac{1}{t}\log{\Ea_{t}^{\epsilon}(\widehat{\mathscr{P}}_0)}=
\frac{1}{t}~\int_0^t (A-\widehat{\mathscr{P}}_sS)~ds~\simeq_{t\rightarrow\infty}~ A-\pi_{\epsilon}(\theta)S,
\end{equation}
where $ \pi_{\epsilon}$ stands for the reversible measure (\ref{iota-pi-1}) when $V>0$ and in (\ref{iota-pi-2}) when $V=0$.
The next corollary is a restatement of Corollary~\ref{cor-ref-adrian} and provides estimates of the Lyapunov exponent (\ref{Lyap-expo-def}).
 \begin{cor}
Assume that $V=0$ and let  $\mbox{\rm Law}(\widehat{\mathscr{P}}_0)=\pi_{\epsilon}$ be the reversible probability measure  defined in (\ref{iota-pi-2}) with $\epsilon=2/\sqrt{N}$. In this situation, for any $t> 0$ we have
 $$
N> 4~~\Longrightarrow~~ -\sqrt{A^2+RS}\leq~ A-\EE(\widehat{\mathscr{P}}_t)~S~\leq~ -\sqrt{A^2+RS~\left(1-\frac{4}{N}\right)}<0.
 $$
 Now assume that $V=S$ and let  $\mbox{\rm Law}(\widehat{\mathscr{P}}_0)=\pi_{\epsilon}$ be the reversible probability measure  defined in (\ref{iota-pi-1}) with $\epsilon=2/\sqrt{N}$. In this situation, for any $t> 0$ we have
 $$
N> 4~~\Longrightarrow~~ -\sqrt{A^2+RS}\leq~ A-\EE(\widehat{\mathscr{P}}_t)~S~\leq~ -\frac{\sqrt{A^2+RS~\left(1-\left({4}/{N}\right)^2\right)}-{4A}/{N}}{1+4/N}<0.
 $$
\end{cor}

\section{Properties of Riccati Semigroups}\label{properties-Riccati-section}
\subsection{Smoothness Properties}\label{smoothness-sec}

We recall that the Riccati semigroup $\phi_{t}$ associated with the parameters  $(A,R,S)$  is given in closed form by the formula
\begin{equation}\label{closed-form}
\phi_t(x)=\varpi_++(x-\varpi_+)~\frac{(\varpi_+-\varpi_-)~e^{-\lambda \, t}}{(\varpi_+-\varpi_-)~e^{-\lambda \, t}+(x-\varpi_-)~(1-e^{-\lambda \, t})}.
\end{equation}
The regularity properties Riccati semigroups are rather well understood in any dimension. We refer to~\cite{Bishop/DelMoral:2016,bd-CARE} for a review on the stability properties
of Riccati semigroups and related Kalman-Bucy diffusion processes. In the one-dimensional case, these properties can be easily checked using the
explicit form given above. For instance, for any $\upsilon>0$ we have the uniform estimates
\begin{equation}\label{unif-phi-star}
0<\phi_{\upsilon}(0)~\leq~ \inf_{x\geq 0}\inf_{t\geq \upsilon}\phi_t(x) ~\leq~ \sup_{x\geq 0}\sup_{t\geq \upsilon}\phi_t(x)~\leq~
\varpi_++
(e^{\lambda \, \upsilon}-1)^{-1}~(\varpi_+-\varpi_-)\leq c_{\upsilon}~(2\varpi_+-\varpi_-).
\end{equation}
When $x\leq \varpi_+$ we have $\phi_t(x)\leq \varpi_+$. When $x> \varpi_+$ for any $t\geq \upsilon$ we have
$$
\frac{\phi_t(x)-\varpi_+}{\varpi_+-\varpi_-}~=~\frac{x-\varpi_+}{(\varpi_+-\varpi_-)e^{\lambda \, t}+(x-\varpi_+)~(e^{\lambda \, t}-1)}~\leq~
\frac{e^{-\lambda \, t}}{1-e^{-\lambda \, \upsilon}}
$$
We also have
 $$
 \phi_t(x)=\varpi_++\frac{\varpi_+-\varpi_-}{\lambda}~\partial_t\log{\left[(\varpi_+-\varpi_-)~e^{-\lambda t}+(x-\varpi_-)~(1-e^{-\lambda t})\right]}
 $$
 This yields for any $s\leq t$ the formula
 \begin{eqnarray*}
 \Ea_{t}(x)&=&\exp{\left[-\lambda t/2\right]}~\frac{1}{(\varpi_+-\varpi_-)~e^{-\lambda t}+(x-\varpi_-)~
 (1-e^{-\lambda t})}\\
&\leq &\exp{\left[-\lambda t/2\right]}~\frac{1}{\varpi_+~e^{-\lambda t}-\varpi_-} \quad\Longrightarrow~~
 \Ea_{t}^{\epsilon}(x)^2\leq (-\varpi_-)^{-1}\exp{\left[-\lambda t\right]}~ \end{eqnarray*}
We also have the rather crude estimates
\begin{equation}\label{unif-phi-star-2}
x\wedge \varpi_+~\leq~ \inf_{t\geq 0}{\phi_t(x)}~\leq~ \phi_{\star}(x):=\sup_{t\geq 0}{\phi_t(x)}=\varpi_+\vee x.
\end{equation}

Observe that the inverse flow $\phi^{-1}_t(x)=1/\phi_t(x)$ satisfies the same equation as in (\ref{nonlinear-KB-Riccati}) by replacing
$(A,R,S)$ by $(A_-,R_-,S_-):=(-A,S,R)$. The extension of this result to the inverse of the stochastic flow is given below.
The proof is a direct application of Ito's formula, thus it is skipped.

\begin{lem}\label{inverse-ricc-lem}
Assume that $R\geq \epsilon^2U$. For any $x>0$, the inverse semigroup $\Phi_t^{-\epsilon}(x):=1/\Phi_{t}^{\epsilon}(x)$ satisfies the stochastic Riccati equation
\begin{equation}\label{inverse-phi-eps}
d\Phi^{-\epsilon}_t(x)=\Lambda_{-\epsilon}\left(\Phi^{-\epsilon}_t(x)\right)~dt+\sigma_{-\epsilon}\left(\Phi^{-\epsilon}_t(x)\right)~dW_t,
\end{equation}
with the drift and the diffusion functions
$$
\Lambda_{-\epsilon}(x)=2A_-~x+R_{-\epsilon}-S_{-\epsilon}x^2, \quad \quad\sigma_{-\epsilon}(x)=\epsilon~\sqrt{ x~\left[U_-~+V_-~x^2\right]},
$$
defined respectively in terms of the parameters
\begin{equation}\label{inverse-phi-eps-ARS}
A_-=-A\qquad R_{-\epsilon}:=(S+\epsilon^2~V)\qquad S_{-\epsilon}:=(R-\epsilon^2~U)\quad\mbox{and}\quad (U_-,V_-)=(V,U).
\end{equation}

\end{lem}

Using elementary differentiations we check the following lemma.
\begin{lem}
For any $x\in\RR_+$ and $n\geq 1$ we have
 \begin{eqnarray}
\partial ^n\phi_t(x)&=&n!~(-1)^{n+1}~\frac{(\varpi_+-\varpi_-)^2~\left[1-e^{-\lambda \, t}\right]^{n-1}~e^{-\lambda \, t}}{
  \left[e^{-\lambda \, t}~(\varpi_+-\varpi_-)+(1-e^{-\lambda \, t})~(x-\varpi_-)\right]^{n+1}}~.
  \label{estimate-second-derivative}
 \end{eqnarray}
 In addition we have the exponential semigroup formula
 \begin{equation}\label{Ea-first-derivative}
 \partial  \phi_t(x)=\Ea_{t}(x)^2\quad\mbox{and the estimates}\quad
  \vert \partial ^n\phi_t(x)\vert\leq n!~\varpi^2\vert\varpi_-\vert^{-(n-1)}~\exp{\left[-\lambda \, t\right]}.
 \end{equation}
 \end{lem}
 The l.h.s. formula in (\ref{Ea-first-derivative}) comes from the evolution equation
 $$
 \partial_t\left(\partial \phi_t(x)\right)=2(A-S\phi_t(x))~\partial \phi_t(x).
 $$

 For instance we have
  \begin{equation}
\varpi^2~\geq  e^{\lambda \, t}~\partial  \phi_t(x)\geq \varpi(x)^2\quad\mbox{\rm with the function}\quad\varpi(x):=1-\frac{(\varpi_+\vee x)-\varpi_+}{(\varpi_+\vee x)-\varpi_-}.
\label{estimate-first-derivative}
   \end{equation}
In the same vein, we have $\varpi^2~e^{-\lambda  t}/\varpi_-\leq 2^{-1}\partial ^2\phi_t(x)<0$.
 This shows that $x\mapsto \phi_t(x)$ is a concave increasing function.

 \subsection{Robustness Properties}\label{robustness-sec}

   Let $\overline{\phi}_t(x)$ be the Riccati semigroup associated with some parameters $(\overline{R},\overline{S})$ such that
 $$
  R\geq \overline{R}\quad \mbox{\rm and}\quad \overline{S}\geq S.
 $$

  We denote by $(\overline{\lambda},\overline{\varpi}_+,\overline{\varpi}_-)$ the parameters defined as $(\lambda ,{\varpi}_+,{\varpi}_-)$ by replacing
  $(R,S) $ by $(\overline{R},\overline{S})$. To simplify the presentation we write $(\overline{\lambda},\lambda)$
  instead of  $\left(\overline{\lambda}(\overline{R},\overline{S}), \lambda \right)$. In this notation it is easily checked that
     \begin{equation}\label{lower-upper-fix}
 \varpi_-\leq \overline{\varpi}_-<0<\overline{\varpi}_+\leq  \varpi_+.
     \end{equation}

  \begin{prop}\label{prop-comparison}
 For any $x\in\RR_+$, we have the estimate
     \begin{equation}\label{lower-estimate}
   0~\leq~ \phi_t(x)-\overline{\phi}_t(x)~\leq ~2~(\lambda+\overline{\lambda})^{-1}~\overline{\varpi}~\varpi~\left(
[~R-\overline{R}~]+[~\overline{S}-S~]~\phi_{\star}(x)^2\right).
   \end{equation}
   \end{prop}

  \proof
  We have
   \begin{eqnarray*}
  \partial_t\left[\phi_t(x)-\overline{\phi}_t(x)\right]
  &=&\left(2A-\overline{S}\left[\phi_t(x)+\overline{\phi}_t(x)\right]\right)\left[\phi_t(x)-\overline{\phi}_t(x)\right]+[~R-\overline{R}~]+[~\overline{S}-S~]~\phi_t(x)^2.
   \end{eqnarray*}
   On the other hand, we have
   $$
      \begin{array}{l}
\displaystyle
   ~\exp{\left[\int_s^t~\left(A-\overline{S}\phi_u(x)~du\right)\right]}
      ~\exp{\left[\int_s^t~\left(A-\overline{S}~\overline{\phi}_u(x)~du\right)\right]}\qquad\qquad\\
      \\
 \displaystyle  \qquad\qquad   \leq  \overline{\varpi}~\varpi~e^{-[\lambda+\overline{\lambda}] (t-s)/2}~
      \exp{\left[\int_s^t~(S-\overline{S})\phi_u(x)~du\right]}\leq  \overline{\varpi}~\varpi~e^{-[\lambda+\overline{\lambda}] (t-s)/2}~~\Longleftarrow~~
      S\leq \overline{S}.
        \end{array}
   $$
   This yields the formula
   $$
   \begin{array}{l}
\displaystyle 0\leq  \phi_t(x)-\overline{\phi}_t(x)\\
\\
~~~=\displaystyle~\int_0^t~\exp{\left[\int_s^t~\left(2A-\overline{S}\left[\phi_u(x)+\overline{\phi}_u(x)\right]~du\right)\right]}~\left(
[~R-\overline{R}~]+[~\overline{S}-S~]~\phi_s(x)^2\right)~ds,
   \end{array}
   $$
   from which the proof of the proposition is easily completed.\qed

\section{Fluctuation Analysis of Riccati Diffusions}\label{fluctuation-sec}

\subsection{Fluctuation Random Fields}
Consider the first and second order fluctuation random fields defined by the formulae
\begin{eqnarray*}
\VV_t^{\epsilon}(x)&:=&\epsilon^{-1}\left[\Phi_{t}^{\epsilon}(x)-\phi_t(x)\right],
\\
\WW_t^{\epsilon}(x)&:=&\epsilon^{-1}\left[\VV^{\epsilon}_t(x)-\VV_t(x)\right],
\quad \mbox{\rm and set}\quad
\overline{\WW}_t^{\epsilon}(x):=
\epsilon^{-1}\left[\EE\left(\WW_t^{\epsilon}(x)\right)-\WW_t(x)\right].
\end{eqnarray*}
In the above display, $\VV_t(x)$ and $\WW_t(x)$ stands for the processes defined by
\begin{eqnarray*}
\VV_t(x)&:=&\int_0^t~\left(\partial \phi_{t-s}\right)(\phi_{s}(x))~\sigma_1\left(\phi_{s}(x)\right)~dW_s,\\
\WW_t(x)&:=&\frac{1}{2}~\int_0^t~\left(\partial ^2\phi_{t-s}\right)(\phi_{s}(x))~\sigma_1^2\left(\phi_{s}(x)\right)~ds~<0.
\end{eqnarray*}
The bias and the fluctuation of $\Phi_{t}^{\epsilon}$ around $\phi_t$ as $\epsilon\rightarrow 0$ are encapsulated respectively
in the deterministic and the stochastic processes $\WW_t(x)$ and $\VV_t(x)$. More precisely we have
$$
\begin{array}{l}
\Phi_{t}^{\epsilon}(x)=\phi_t(x)+\epsilon~\VV_t(x)+\epsilon^2~\WW_t^{\epsilon}(x)
\quad\mbox{\rm and}\quad
~\EE\left[\Phi_{t}^{\epsilon}(x)\right]=\phi_t(x)+\epsilon^2~\WW_t(x)+\epsilon^3~\overline{\WW}_t^{\epsilon}(x).
\end{array}
$$

The first part of the article is concerned with quantitative and uniform estimates of the fluctuation random fields introduced above.
Most of the estimates developed in the article are expressed in terms of judiciously chosen collections of Riccati semigroups.

\subsection{Uniform Fluctuation Estimates}\label{section-intro-sg-fluctuation}
The objective of this section is to analyze the fluctuations of the random fields $\VV_{t}^{\epsilon}$ and $\WW_t^{\epsilon}$ in terms
of the collection of Riccati semigroups  defined below.

\begin{defi}
We let $\phi^{(\epsilon,n)}$ be the collection of Riccati semigroups indexed by the parameters
$\epsilon\in\RR_+$ and $n\in\RR$, and defined as $\phi_t$ by replacing $(R,S)$ by the parameters
\begin{eqnarray}
(R^{(\epsilon,n)},S^{(\epsilon,n)})&:=&(R,S) +(n-1)\,\frac{\epsilon^2}{2}\,\left(U,-V\right).\label{def-RS-epsilon-n}
\end{eqnarray}
We also denote by $(\varpi^{(\epsilon,n)}_-,\varpi^{(\epsilon,n)}_+,\phi^{(\epsilon,n)}_{\star},\ldots)$, the objects defined similarly to $(\varpi_-,\varpi_+,\phi_{\star},\ldots)$
but with $(R,S)$ replaced by the parameters $(R^{(\epsilon,n)},S^{(\epsilon,n)})$ in the corresponding definition.
\end{defi}
The semigroups $\phi_t^{(\epsilon,n)}$, resp. $\phi_t^{(\epsilon,-n)}$ indexed by $n\geq 0$
 are well founded as soon as
 \begin{equation}\label{def-RS-epsilon-n-well-founded}
 (n-1)~\epsilon^2~\overline{V}< 2\quad\mbox{\rm and respectively}\quad
 (n+1)~\epsilon^2\,\overline{U}< 2.
  \end{equation}
 Observe that when $V=0$ the flow $\phi^{(\epsilon,n)}$ is well defined for any $\epsilon \in \RR_+$ and any $n\geq 0$.
In addition we have $$\epsilon=0~~\Longrightarrow~~\phi^{(0,n)}_t=\phi_t={\phi}^{(0,-n)}_t\quad\mbox{\rm and}\quad
\phi_t=\phi^{(\epsilon,1)}_t.
$$

The $\LL_n$-norm of the fluctuation random fields $\VV^{\epsilon}_t(x)$ will be estimated in term of the collection of functions
$v^{\epsilon}_{n}(x)$ defined by
\begin{eqnarray*}
v^{\epsilon}_{n}(x)&:=&\varpi_{\lambda}~\left[
\frac{ \epsilon}{\sqrt{2\lambda }}~\sigma_1^2\left(\phi^{(\epsilon,3n)}_{\star}(x)\right)-\frac{n}{2}~\varpi_-~\sigma_1\left(\phi^{(\epsilon,{3n}/{2})}_{\star}(x)\right)
  \right]\\
  &&\longrightarrow_{\epsilon\rightarrow 0}~\frac{n~\varpi^2}{\sqrt{2\lambda }}~\sigma_1\left(\phi_{\star}(x)\right)\quad\mbox{\rm with the parameter
  $
  \displaystyle\varpi_{\lambda}=-\frac{\varpi^2}{\varpi_-}~\sqrt{\frac{2}{\lambda }}~
  $.}
  \end{eqnarray*}

 The $\LL_n$-norm of the random field $\WW^{\epsilon}_t(x)$ will be estimated in term of the collection of functions:
 $$
 w^{\epsilon}_{n}:= w^{\epsilon}_{1,n}(x)+{w}^{\epsilon}_{2,n}(x)\quad\mbox{\rm and}\quad   v_{\epsilon,n}(x):=1\vee {\vertiii{\VV^{\epsilon}(x)}}_n,
 $$
\begin{eqnarray*}
\mbox{\rm with}\quad w^{\epsilon}_{1,n}(x)&:=&\varpi_{\lambda}~\left[~n~v_{\epsilon,n}(x)~\sigma_1\left(\phi_{\star}(x)\right)+~\frac{1}{\sqrt{2\lambda }}~\sigma_1^2\left(\phi^{(\epsilon,3n)}_{\star}(x)\right)\right], \\
{w}^{\epsilon}_{2,n}(x)&:=&\frac{3n}{2}~\varpi_{\lambda}~~\left[
U~\phi^{-}_{\star}(x)+\sqrt{UV}+V~\phi^{(\epsilon,n)}_{\star}(x)/2
\right]^{1/2}~\left[2\epsilon~v_{\epsilon,4n}^2-v_{\epsilon,2n}(x)~\varpi_-\right].
\end{eqnarray*}
In the above display, $\phi^{-}_{\star}(x)$  stands for the supremum of the inverse $\phi^{-1}_{t}(x)$ Riccati semigroup w.r.t. the time horizon.
Finally, the estimate of the bias $\overline{\WW}_t^{\epsilon}(x)$ developed below is expressed in terms of the functions:
$$
 \begin{array}{l}
\displaystyle\overline{w}^{\epsilon}(x)\\
\displaystyle:= \overline\varpi_{\lambda}~v_{\epsilon,4}(x)^2~
\left[\sigma_1^2\left(\phi_{\star}(x)\right)+\left({U}/{3}+4V~\phi^{(\epsilon,4)}_{\star}(x)^2\right)~\left(
3\epsilon-\varpi_-\right)\right]~\quad\mbox{\rm with}\quad\overline\varpi_{\lambda}:=\frac{3}{\lambda }~\left(\frac{\varpi}{\varpi_-}\right)^2.
\end{array}
$$
We are now in a position to state the main result of this section.
\begin{theo}\label{theo--1intro}
For any $x\in\RR_+$, $\epsilon\in \RR_+$, and any $n\geq 1$ such that (\ref{def-RS-epsilon-n-well-founded}) holds, we have the norm estimates
\begin{equation}\label{bias-n}
{\phi}^{(\epsilon,-1)}_t(x)\leq \vertiii{\Phi_{t}^{\epsilon}(x)}_n\leq \phi^{(\epsilon,n)}_t(x)
\quad\mbox{and}\quad
{\phi}^{(\epsilon,-n)}_t(x)
\leq
\vertiii{\Phi_{t}^{\epsilon}(x)^{-1}}_n^{-1}\leq \phi_t(x).
\end{equation}
In addition, we have the uniform fluctuation estimates
\begin{equation}
{\vertiii{\VV^{\epsilon}(x)}}_n\leq ~v^{\epsilon}_{n}(x)
\qquad
\displaystyle {\vertiii{\WW^{\epsilon}(x)}}_n\leq ~w^{\epsilon}_{n}(x)\quad\mbox{and}\quad
\sup_{t\geq 0}{\left\vert\,\overline{\WW}_t^{\epsilon}(x)\,\right\vert}\leq ~\overline{w}^{\epsilon}(x).
\label{unif-est-intro}
\end{equation}
\end{theo}
The proof of the $n$-th moment estimates (\ref{bias-n}) and the uniform estimates (\ref{unif-est-intro}) are lengthy and technical. They are provided respectively in Section~\ref{defi-phi-n-eps-section}--\ref{proof-theo-intro}.

The estimates stated in the above theorem are sharp when $\epsilon\rightarrow 0$, in the sense that
all the Riccati semigroups discussed above converge to  $\phi_t$ as $\epsilon\downarrow 0$. We end this section with some comments
about these properties.

Firstly, observe that  for any parameters $0\leq \epsilon_2\leq \epsilon_1$ and $1\leq n_1\leq n_2$
satisfying (\ref{def-RS-epsilon-n-well-founded})  we have the following monotonicity properties:
\begin{equation}\label{monotone-epsilon-n}
{\phi}^{(\epsilon_1,-n_2)}_t\leq{\phi}^{(\epsilon_1,-n_1)}_t\leq{\phi}^{(\epsilon_2,-n_1)}_t\leq \phi_t
\leq \phi^{(\epsilon_2,n_1)}_t\leq \phi^{(\epsilon_1,n_1)}_t\leq  \phi^{(\epsilon_1,n_2)}_t.
\end{equation}
The above inequalities are
direct consequences of the estimate (\ref{lower-estimate}) stated in Proposition~\ref{prop-comparison}. The difference between $\phi^{(\epsilon,n)}_t$
and $\phi_t$  can be quantified in terms of the parameter $\epsilon$ using the estimates stated in Proposition~\ref{prop-comparison}.

When $V>0$ we emphasize that the Riccati semigroups $\phi^{(\epsilon,n)}_t(x)$ are only defined when the l.h.s. condition in (\ref{def-RS-epsilon-n-well-founded}) is satisfied. Therefore, we cannot expect to have an uniform estimate of $\phi^{(\epsilon,n)}_{\star}(x)$ w.r.t. the fluctuation parameter $\epsilon$ for any values of $n$.

\section{Stability Analysis of Riccati Diffusions}\label{ricc-diff-stab-sec}

\subsection{Reversible Probability Measures}\label{inv-meas-sec}
The infinitesimal generator $L$ of the diffusion $X_t$ defined in (\ref{mp-ref}) is given  by the differential operator
\begin{eqnarray}\label{liouville-ref-0}
L=\frac{1}{2}~\sigma_{\epsilon}^2~\partial ^2+\Lambda~\partial
~\Longrightarrow ~
df(X_t)=L(f)(X_t)~dt+dM_t(f).
\end{eqnarray}
Recall that the martingale $dM_t(f):=\sigma_{\epsilon}(X_t)~\partial f(X_t)$ has an angle bracket given for any
smooth functions $f$ and $g$ on $\RR$ by
$$
\partial_t\langle M(f),M(g)\rangle_t=\Gamma_{L}(f,g)(X_t),
$$
with the ``carr\'e du champ'' operator
$$
\Gamma_{L}(f,g):=L(fg)-gL(f)-fL(g)=\sigma_{\epsilon}^2~\partial f~\partial g.
$$
The second order differential operator $L$ discussed above can be rewritten as
\begin{eqnarray}\label{liouville-ref}
L(f)
&=&\frac{1}{2}~\sigma_{\epsilon}^2~e^{-\,\Ua_{\epsilon}}~\partial \left( e^{\,\Ua_{\epsilon}}~\partial f\right)\quad\mbox{\rm with}\quad
\Ua_{\epsilon}(x):=2~\int_{p}^{x}~\frac{\Lambda(y)}{\sigma^2_{\epsilon}(y)}~dy,
\end{eqnarray}
where $p$ is an arbitrary point in $]0,\infty[$. We set
$$
m_{\epsilon}(x)=:\int_0^x~q_{\epsilon}(y)~dy\quad\mbox{\rm with}\quad q_{\epsilon}(x):=2\sigma^{-2}_{\epsilon}(x)~\exp{\left(\Ua_{\epsilon}(x)\right)}\quad
\mbox{\rm and}\quad s_{\epsilon}(x):=\int_p^x \exp{\left(-\Ua_{\epsilon}(y)\right)}~dy.
$$
The unnormalized measure with density $m_{\epsilon}(x)$ is often called a speed measure, and $s_{\epsilon}$ a scale function. We further assume that $0\leq\epsilon^2 \overline{U}<2$.

The Sturm-Liouville formulation of the generator $L$ given in (\ref{liouville-ref}) shows that a reversible measure of the Riccati diffusion (\ref{mp-ref}) is given by the formula
$$
\pi_{\epsilon}(dx)=\frac{1}{\Za_{\epsilon}}~1_{\RR_+}(x)~\frac{1}{\sigma^2_{\epsilon}(x)}~\exp{\left(\Ua_{\epsilon}(x)\right)}~dx,
$$
where $\Za_{\epsilon}$ stands for some normalizing constant.
Different type of reversible probability distributions can obtained depending on the choice of the parameters.
\begin{itemize}
\item When $U\wedge V>0$, the measure $
\pi_{\epsilon}$ is the heavy tailed probability measure
\begin{equation}\label{iota-pi-1}
\pi_{\epsilon}(dx)~\propto~1_{\RR_+}(x)~\frac{x^{\frac{2}{\epsilon^2}\frac{R}{U}-1}}{\left[U+Vx^2\right]^{1+
\frac{1}{\epsilon^2}\left(\frac{R}{U}+\frac{S}{V}\right)}}~\exp{\left[\frac{4}{\epsilon^2}~\frac{A}{\sqrt{UV}}~\tan^{-1}\left(x~\sqrt{\frac{V}{U}}\right)\right]}~dx.
\end{equation}
Observe that
\begin{equation}\label{iota-1-pi}
\int_0^{\infty}~x^{- \imath}~\pi_{\epsilon}(dx)<\infty ~~\Longleftrightarrow~~ - \imath~\epsilon^2~ \overline{V}<2\quad \mbox{\rm and}\quad
 \imath~\epsilon^2 ~\overline{U}<2,
\end{equation}
with the parameter $ \imath$ defined in (\ref{def-iota-iota-kappa}). More specifically, when $A<0$ the condition $ \imath~\epsilon^2~ \overline{U}<2$ is automatically satisfied. Conversely, when  $A>0$, the condition $-\imath\,\epsilon^2 \,\overline{V}<2$ is satisfied.
\item When $U>0$ and $V=0$, the probability measure $
\pi_{\epsilon}$ reduces to
\begin{equation}\label{iota-pi-2}
\pi_{\epsilon}(dx)~\propto~1_{\RR_+}(x)~x^{\frac{2}{\epsilon^2}\frac{R}{U}-1}~\exp{\left[-\frac{S}{U\epsilon^2}~\left(x-2~\frac{A}{S}\right)^2\right]}~dx.
\end{equation}
Observe that
\begin{equation}\label{iota-2-pi}
\int_0^{\infty}~x^{- \imath}~\pi_{\epsilon}(dx)<\infty ~~\Longleftrightarrow~~
 \imath~\epsilon^2 ~\overline{U}<2.
\end{equation}
In contrast with (\ref{iota-1-pi}), the r.h.s. condition of \eqref{iota-2-pi} is automatically satisfied for any $A<0$.
\item  When  $V>0$ and $U=0$, the probability measure $
\pi_{\epsilon}$ takes the form
$$
\pi_{\epsilon}(dx)~\propto~1_{\RR_+}(x)~x^{-\left[\frac{2}{\epsilon^2}\frac{S}{V}+3\right]}~\exp{\left[-\frac{R}{V\epsilon^2}~\left(\frac{1}{x}+2~\frac{A}{R}\right)^2\right]}~dx.
$$
Observe that
$$
\int_0^{\infty}~x^{- \imath}~\pi_{\epsilon}(dx)<\infty ~~\Longleftrightarrow~~ \imath~\epsilon^2~ \overline{V}<2.
$$
\item When $U>0$ and $V=0=S>A$,  the probability measure $
\pi_{\epsilon}$ reduces to the Gamma distribution
$$
\pi_{\epsilon}(dx)~\propto~1_{\RR_+}(x)~x^{\frac{2}{\epsilon^2}\frac{R}{U}-1}~\exp{\left(\frac{4}{\epsilon^2}~\frac{A}{U}~x\right)}~dx.
$$
\end{itemize}
When $\epsilon^2\,\overline{U}<2$, the function
 $\partial_x s_{\epsilon}$ is not integrable around the origin so that $0$
 is repelling. On the other hand, the function
$q_{\epsilon}$ is integrable around the origin so that $0$ is also a regular boundary state.
We also have
$\lim_{x\rightarrow\infty}m_{\epsilon}(x)<\infty$ and $\lim_{x\rightarrow\infty}s_{\epsilon}(x)=\infty$. This shows that the boundary states $0$ and $\infty$
are both regular and repellent. For a more thorough discussion on the classification of boundary states we refer to~\cite{abundo1997some,Ito-McKean-65}, and the more recent review articles~\cite{fukushima2015feller} and \cite{peskir2014boundary}.

Applying (\ref{bias-n}) with $n=1$ and letting $t\rightarrow\infty$ we have the bias estimate
$$
\frac{A+\sqrt{A^2+RS~\left(1-\epsilon^2\overline{U}\right)\left(1+\epsilon^2\overline{V}\right)}}{S\left(1+\epsilon^2\overline{V}\right)}~\leq ~\pi_{\epsilon}(\theta)~\leq~ \varpi_+.
$$
with $\theta(x):=x$. This yields the following corollary.
\begin{cor}\label{cor-ref-adrian}
For any $\epsilon^2\,\overline{U}<1$, we have
\begin{equation}
A-\varpi_+S~\leq~ A-\pi_{\epsilon}(\theta)S~\leq~ -\frac{\sqrt{A^2+RS~\left(1-\epsilon^2\overline{U}\right)\left(1+\epsilon^2\overline{V}\right)}-\epsilon^2~A\overline{V}
}{1+\epsilon^2\overline{V}}.
\end{equation}
\end{cor}

\subsection{Tangent Processes}\label{tangent-proc-sec}
The second part of the article is dedicated to the stability properties of the Riccati diffusion process.
Firstly, observe that the long time behavior of one-dimensional Riccati semigroups is encapsulated into the exponential decays of the tangent process $\tau_t(x)$ defined below
\begin{equation}\label{deterministic-expo-A-PS}
\begin{array}{rcl}
(\ref{Ea-first-derivative})\quad\mbox{\rm and}\quad (\ref{estimate-first-derivative})&\Longrightarrow&\displaystyle
\varpi(x)^2~e^{-\lambda t}~\leq~ \tau_t(x):=\partial  \phi_t(x)=
\Ea_{t}(x)^2  ~\leq~ \varpi^2~\exp{\left[-\lambda \, t\right]}\\
&&\\
&\Longrightarrow&\displaystyle~\frac{1}{t}~\log{\tau_t(x)}\longrightarrow_{t\rightarrow\infty}~\partial \Lambda(\varpi_+)=-\lambda.
\end{array}
\end{equation}
 In the same manner, $x\rightarrow\Phi_{t}^{\epsilon}(x)$ is almost surely differentiable, and the tangent process
 $$
 \Ta_t^{\epsilon}(x):=\partial \Phi_{t}^{\epsilon}(x)
~~\Longrightarrow~~
d\,\Ta_t^{\epsilon}(x)=\left[\partial \Lambda(\Phi_{t}^{\epsilon}(x))~dt+\partial \sigma_{\epsilon}(\Phi_{t}^{\epsilon}(x))~dW_t\right]~\Ta_t^{\epsilon}(x)
$$ is given by the exponential formula
\begin{eqnarray}
\sigma_{1}(x)~\Ta_t^{\epsilon}(x)&=&\sigma_{1}(\Phi_{t}^{\epsilon}(x))~ \exp{\left[-\int_0^t~\Ha^{\epsilon}(\Phi^{\epsilon}_s(x))~ds\right]}\leq \sigma_{1}(\Phi_{t}^{\epsilon}(x))~\exp{\left[-\lambda_{\epsilon}\,t\right]},\label{tangent-epsilon}
\end{eqnarray}
with the potential function $\Ha^{\epsilon}$ on $\RR_+$ defined by
\begin{equation}\label{def-Va-epsilon}
\Ha^{\epsilon}:=\sigma_{\epsilon}^{-1}L(\sigma_{\epsilon})-\partial \Lambda=-\partial \Lambda+\frac{\partial \sigma_1}{\sigma_1}~\Lambda+
\frac{\epsilon^2}{2}~\sigma_1~\partial ^2\sigma_1\quad\mbox{\rm and}\quad\lambda_{\epsilon}\,:=\inf_{x\geq 0}{\Ha^{\epsilon}(x)}
\end{equation}
The r.h.s. estimate in (\ref{tangent-epsilon}) provides an almost sure exponential decay of the tangent process.
In addition,  using (\ref{def-sg}) for any bounded differentiable function $f$ on $\RR_+$ and any $x>0$ we have the formula
\begin{equation}\label{derivative-sg}
 \partial P_t^{\epsilon}(f)(x)=\EE\left[\,\partial f(\Phi_{t}^{\epsilon}(x))~\Ta_t^{\epsilon}(x)\,\right].
 \end{equation}
To get one step further in our discussion, let $\overline{\Phi}_{t}^{\epsilon}$ be the stochastic Riccati flow associated with the parameters
\begin{equation}\label{def-Repsilon-delta}
\overline{R}_{\epsilon}:=R\left(1+\frac{\epsilon^2}{2}~\overline{U}\right)\quad\mbox{\rm and}\quad \overline{S}_{\epsilon}=S\left(1-\frac{3\epsilon^2}{2}~\overline{V}\right).
\end{equation}
 In this notation, the stochastic version of (\ref{deterministic-expo-A-PS}) is given by the following proposition.
\begin{prop}\label{FK-intro-prop}
For any bounded measurable function $f$ on $\RR_+$, any $x\in\RR_+$, and any time horizon $t\geq 0$
we have the Feynman-Kac formula
\begin{equation}\label{FK-intro}
\EE\left[f(\Phi_{t}^{\epsilon}(x))\,\Ta_t^{\epsilon}(x)\right]=\EE\left[f(\overline{\Phi}_{t}^{\epsilon}(x))\,\exp{\left[\int_0^t~\partial \Lambda\left(\,\overline{\Phi}^{\epsilon}_s(x)\right)~ds\right]}\right].
\end{equation}
\end{prop}
The proof of the above proposition is provided in Section~\ref{FK-intro-prop-proof}.

\subsection{Contraction Inequalities}\label{contract-ineq-sec}

The next lemma that can be used to quantify the almost sure decay rate  of the tangent process (\ref{tangent-epsilon}) in terms
 of the drift parameters $(A,R,S)$ and the diffusion parameters $(U,V)$.
\begin{lem}\label{lem-3}
Let $\lambda_{\epsilon}$ be the parameter defined in (\ref{def-Va-epsilon}).
The following assertions hold:
\begin{eqnarray*}
V=0\quad \mbox{and}\quad \epsilon^2\,\overline{U}\leq 2&\Longrightarrow&
\lambda_{\epsilon}\,=-A+\sqrt{3RS\left(1-\frac{\epsilon^2}{2}~\overline{U}\right)};\\
U=0\quad \mbox{and}\quad \epsilon^2~\overline{V}\leq 2/3&\Longrightarrow&
\lambda_{\epsilon}\,=A+\sqrt{3RS\left(
1-\frac{3\epsilon^2}{2}~\overline{V}~\right)};\\
 \epsilon^2\,\overline{U}\leq 2~~~\mbox{and}~~~ \frac{3\epsilon^2\overline{V}}{1+3\epsilon^2V}\leq 2/3
 &\Longrightarrow&
\lambda_{\epsilon}\,\geq -\vert A\vert+\sqrt{RS\left[1-\frac{\epsilon^2}{2}~\overline{U}\right]\left[1+3~\epsilon^2~\overline{V}\left(S-\frac{3}{2}\right)\right]}.
\end{eqnarray*}
\end{lem}
The proof of this lemma is technical, and it is given in the Appendix~\ref{proof-lem-3}.

\begin{theo}\label{theo-2-inside}
Assume that $\lambda_{\epsilon}$ defined in (\ref{def-Va-epsilon}) is positive. The Wasserstein contraction inequality (\ref{Wasserstein-estimate}), and the Poincar\'e inequality (\ref{poincare-inequality}), stated in Theorem~\ref{theo-2-intro}
are satisfied.
\end{theo}
The proof of the preceding theorem is given in Section~\ref{theo-2-intro-proof}. The preceding theorem is more general than the example case given in Theorem~\ref{theo-2-intro} with $V=0$.

The estimates stated in Lemma~\ref{lem-3} and  Theorem~\ref{theo-2-inside} are clearly not satisfactory when $U\not=0$ and $A>0$ as they require that
$A<\sqrt{RS}$. However, when $U=0$, the estimates are not really useful when $A<0$, as they require that $\vert A\vert <\sqrt{3RS}$.
Roughly speaking, to improve the
estimates discussed above, we need to interpolate between these two cases. More precisely, observe that the semigroup derivative formula
(\ref{derivative-sg}) is expressed  in terms of the exponential type process (\ref{tangent-epsilon}).  Formula (\ref{derivative-sg}) can be seen as the expectation w.r.t. to law of the process $\Phi_{t}^{\epsilon}(x)$ weighted by some exponential potential function.
In terms of importance sampling techniques, these Feynman-Kac formulae can be seen as the integral w.r.t. a twisted process which is more likely to visit regions
with low energy type $\Ha^{\epsilon}$-values. The state regions with negative $\Ha^{\epsilon}$-values are more likely to be visited but the
lower bound $\lambda_{\epsilon}\,<0$ doesn't give any information on the killing rate of the process.

To obtain some more useful estimates, we seek a judicious change of measure under which the potential function is an absorbing potential.
To describe precisely these twisted models we need to introduce some additional notation.
We set
\begin{equation}\label{Riota-ref}
\widehat{R}_{\epsilon}=R\left(1-\epsilon^2~\left( \frac{1}{2}+\imath\right)~\overline{U}\right)\quad \mbox{\rm and}\quad
\widehat{S}_{\epsilon}=S\left(1-\epsilon^2~\left( \frac{1}{2}-\imath\right)~\overline{V}\right).
\end{equation}
Let $\widehat{\epsilon}$ be the smallest parameter $\epsilon\in \RR_+$ such that $\widehat{R}_{\epsilon}\wedge \widehat{S}_{\epsilon}>0$.
\begin{defi}
For any $\epsilon\in [0,\widehat{\epsilon}]$, let
 $\widehat{\Phi}^{\epsilon}_t(x)$ be the stochastic Riccati flow  associated with the parameters $(\widehat{R}_{\epsilon},\widehat{S}_{\epsilon})$, and denote by $\widehat{P}^{\epsilon}_t$ the corresponding transition semigroup
 $$
\widehat{P}^{\epsilon}_t(f)(x)=\EE\left(f\left(\widehat{\Phi}^{\epsilon}_t(x)\right)\right).
 $$
 \end{defi}

We also let $\widehat{\Ha}_{\epsilon}$ be the collection of potential functions on $\RR_+$  defined by
\begin{eqnarray*}
\widehat{\Ha}_{\epsilon}(x)
=2\imath \left[A-\left[1+\frac{\epsilon^2}{2}~ \imath_1~\overline{V}\right]~Sx\right]+~ \imath_1~
\left(\left[1-\frac{\epsilon^2}{2}~ \imath_1~\overline{U}\right]~\frac{R}{x}+\left[1+\frac{\epsilon^2}{2}~ \imath_1~\overline{V}\right]~Sx\right).
\end{eqnarray*}
and we consider the the tangent-type process $\widehat{\Ta}^{\epsilon}_t(x)$ defined by
\begin{equation}\label{tangent-overline-estimate}
\widehat{\sigma}(x)~\widehat{\Ta}^{\epsilon}_t(x):=\widehat{\sigma}\left[\widehat{\Phi}^{\epsilon}_t(x)\right]\,\exp{\left[-\int_0^t\widehat{\Ha}_{\epsilon}(\widehat{\Phi}^{\epsilon}_s(x))~ds\right]}\quad\mbox{\rm with}\quad \widehat{\sigma}(x):=x^{ \imath_1 }.
\end{equation}

We are now in a position to state the main result of this section.
\begin{theo}\label{theo-3-intro}
For any $x\in\RR_+$, any bounded differentiable function $f$ on $\RR_+$, any $\epsilon\in [0,\widehat{\epsilon}\,]$, and any time horizon $t\geq 0$,
we have the Feynman-Kac  formulae
\begin{equation}\label{tangent-overline}
\displaystyle\EE\left[f\left(\Phi^{\epsilon}_t(x)\right)~\Ta^{\epsilon}_t(x)\,\right]=
\EE\left[f\left(\widehat{\Phi}^{\epsilon}_t(x)\right)~\widehat{\Ta}^{\epsilon}_t(x)\,\right]\quad\mbox{and}\quad \inf_{x\geq 0}\widehat{\Ha}_{\epsilon}(x)\geq \widehat{\lambda}_{\epsilon}.
\end{equation}
In this situation, for any smooth function $f$, we have the commutation formula
$$
\widehat{\sigma}\,\left|\partial P_t^{\epsilon}(f)\right|
\,\leq\, \exp{\left[-\widehat{\lambda}_{\epsilon} \,t\right]}~\widehat{P}^{\epsilon}_t\left(\widehat{\sigma}\,\left|\partial f\right|\,\right) ~~~\Longrightarrow~~~(\ref{Wasserstein-estimate-overline}).
$$
\end{theo}
The proof of the above theorem is given in Section~\ref{theo-3-intro-proof}. The proof of the implication of (\ref{Wasserstein-estimate-overline})  in Theorem~\ref{t2-intro} follows the same line of arguments as the proof of (\ref{Wasserstein-estimate}), so it is omitted.

Then next corollary is a direct consequence of \eqref{Wasserstein-estimate-overline} and the estimates \eqref{iota-1-pi} and \eqref{iota-2-pi}.
 \begin{cor}
For any $x>0$ we have
 $$
  \epsilon^2~\left[\zeta\vee\left( \imath\,\overline{U}\right)\vee \left(- \imath\,\overline{V}\right)\right]<2
~~\Longrightarrow~~ \lim_{t\rightarrow\infty}\DD_{\widehat{\sigma}}\left(\delta_xP^{\epsilon}_t,\pi_{\epsilon}\right)=0
 $$
where $\pi_{\epsilon}$ is the measure  defined in (\ref{iota-pi-1}) when $V>0$, and in (\ref{iota-pi-2})
 when $V=0$.
 \end{cor}

\subsection{Exponential Semigroups}\label{exp-semigroups-sec}

Our next result is expressed in terms of the collection of uniformly bounded functions
\begin{equation}\label{rho-widehat-def}
\displaystyle
\rho_{\epsilon,\kappa}(x):=\displaystyle\left[
1
+\frac{\varpi_+}{{\phi}^{\,(\epsilon,-( \imath_{\kappa}\vee 1))}_{1}(0)\vee x}\right]
^{ \imath_1}
\exp{\left[3\lambda \right]}~\leq~\displaystyle \widehat{\rho}_{\epsilon,\kappa}:=\sup_{x\geq 0}\rho_{\epsilon,\kappa}(x)=\rho_{\epsilon,\kappa}(0).
\end{equation}
By (\ref{monotone-epsilon-n}) one can check that $\epsilon\rightarrow \rho_{\epsilon,\kappa}(0)$ is a decreasing function. For any $\kappa> 0$, we define $\epsilon_{\kappa}$ be the smallest parameter $\epsilon \in \RR_+$ s.t. $\widehat{\lambda}_{\epsilon,\kappa} \geq 0$, i.e.
\begin{equation}\label{condi-unif-laplace-estimates-bis}
\epsilon_{\kappa} ~:=~ \inf\left\{ \epsilon\in\RR_+~~\mathrm{s.t.}~~\frac{\epsilon^2}{2}~\zeta_{\kappa}<1\right\}
\end{equation}

We have the following moment stability result on the stochastic exponential semigroup. The stability of the deterministic version of this semigroup is a fundamental result in Kalman-Bucy filtering due originally to Bucy; see \cite{Bishop/DelMoral:2016,bd-CARE}.

\begin{theo}\label{theo-expo-sg-estimate}
For any $\kappa> 0$, let $\epsilon_{\kappa}$ be defined as in (\ref{condi-unif-laplace-estimates-bis}). For any $x\in\RR_+$, any time horizon $t\geq 0$, and any $\epsilon \in [0,\epsilon_{\kappa}]$, we have the Laplace estimates
\begin{equation}\label{unif-laplace-estimates-bis}
\begin{array}{l}
 \varpi(x)^2\,\exp{\left[-\lambda \,t\right]}~\leq~
\displaystyle{
\EE\left[
\Ea_{t}^{\epsilon}(x)^{2\kappa}\right]^{1/\kappa}}
~\leq~ \rho_{\epsilon,\kappa}(x)\,
\exp{\left[-\widehat{\lambda}_{\epsilon,\kappa} ~t\right]}.
\end{array}
\end{equation}
\end{theo}
The proof of the preceding theorem is provided in Appendix~\ref{theo-expo-sg-estimate-proof}. Observe that
\begin{equation}\label{comparison-lambda-hat}
\kappa_1\leq \kappa_2\quad\mbox{\rm and}\quad
\epsilon_1\leq \epsilon_2
~\Longrightarrow~
\lambda_{\kappa_2,\epsilon_2} \leq \lambda_{\kappa_1,\epsilon_1} \quad\mbox{\rm and}\quad
\rho_{\epsilon_1,\kappa_1}(x)\leq \rho_{\epsilon_2,\kappa_2}(x).
\end{equation}

 To get one step further, note that
\begin{equation}\label{linearization-Phi}
\Lambda(x_1)-\Lambda(x_2)=\left[(A-Sx_1)+(A-Sx_2)\right]~(x_1-x_2).
\end{equation}
 For any $x_1,x_2\in\RR_+$ we set
$\Ea_{s,t}^{\epsilon}(x_1,x_2)=\Ea_{s,t}^{\epsilon}(x_1)\Ea_{s,t}^{\epsilon}(x_2)$
with the exponential semigroup $\Ea_{s,t}^{\epsilon}(x)$ defined in (\ref{def-Ea-st}).
 Using (\ref{linearization-Phi}) we have
$$
\Phi^{\epsilon}_t(x_1)-\Phi^{\epsilon}_t(x_2)=\Ea_{0,t}^{\epsilon}(x_1,x_2) (x_1-x_2)+\epsilon~
\int_0^t~
\Ea_{s,t}^{\epsilon}(x_1,x_2)~\left[\sigma_{1}\left(\Phi^{\epsilon}_t(x_1)\right)-\sigma_{1}\left(\Phi^{\epsilon}_t(x_2)\right)\right]~dW_s.
$$
This implies that
$$
\EE\left[\Phi^{\epsilon}_t(x_1)-\Phi^{\epsilon}_t(x_2)\right]=\EE\left[\Ea_{0,t}^{\epsilon}(x_1,x_2)\right]~(x_1-x_2).
$$
Combining Theorem~\ref{theo-expo-sg-estimate} with Cauchy-Schwartz inequality  we readily get  the following proposition.
\begin{prop}
Let $\epsilon_1$ be a parameter defined as in (\ref{condi-unif-laplace-estimates-bis}). For any $x_1\geq x_2\in\RR_+$, any $\epsilon\in [0,\epsilon_1]$, and any time horizon $t\geq 0$
we have
$$
\displaystyle 0~\leq~ \EE\left[\Phi^{\epsilon}_t(x_1)\right]-\EE\left[\Phi^{\epsilon}_t(x_2)\right]~\leq~
\widehat{\rho}_{\epsilon,1}
\exp{\left[-\widehat{\lambda}_{\epsilon,1} ~t\right]}~(x_1-x_2),
$$
with the parameters $\widehat{\lambda}_{\epsilon,1} $ and $\widehat{\rho}_{\epsilon,1}$ introduced in (\ref{def-lambda-epsilon-iota}) and (\ref{rho-widehat-def}).
\end{prop}

Recall that $\widehat{\Phi}^{\epsilon}_t(x)$ is the stochastic Riccati flow  associated with the parameters $(\widehat{R}_{\epsilon},\widehat{S}_{\epsilon})$. Following Theorem~\ref{theo--1intro}, the $n$-moments of
this flow can be estimated in terms of the  Riccati semigroups $\widehat{\phi}^{(\epsilon, n)}_t(x)$ associated with the parameters
\begin{eqnarray}
\left(\widehat{R}^{(\epsilon, n)},\widehat{S}^{(\epsilon, n)}\right)&=&(\widehat{R}_{\epsilon},\widehat{S}_{\epsilon})+(n-1)~\frac{\epsilon^2}{2}~\left(U,-V\right).
\label{def-phi-check}
\end{eqnarray}
Let $(\widehat{\varpi}^{(\epsilon,n)}_-,\widehat{\varpi}^{(\epsilon,n)}_+)$ be the parameters defined as $(\varpi_-,\varpi_+)$
by replacing the pair $(R,S)$ by (\ref{def-phi-check}). In this notation, we have the following theorem.

\begin{theo}\label{theo-4-intro}
For any $x\in\RR_+$, any $t>0$, any $n\geq 1$, and any $\epsilon\in\RR_+$ s.t. $$
\epsilon^2\,\left[~\zeta
\vee\left(\left(\imath_{n-2}+2\right)\,\overline{V}\right)\right]<2,
$$
we have the contraction inequality
\begin{equation}\label{Wasserstein-estimate-overline-bis}
\vertiii{\Phi^{\epsilon}_t(x_1)-\Phi^{\epsilon}_t(x_2)}_n\leq ~\left[2 \widehat{\varpi}^{(\epsilon, \imath_n)}_+-\widehat{\varpi}^{(\epsilon, \imath_n)}_-\right]^{ \imath_1}~d_{\widehat{\sigma}}(x_1,x_2)~\exp{\left[-\widehat{\lambda}_{\epsilon} \,t\right]}~.
\end{equation}

\end{theo}
The proof of the preceding theorem is given in Section~\ref{theo-4-intro-proof}.

\section{Stochastic Ornstein-Uhlenbeck Processes}\label{2-d-section}

\subsection{Stability Properties}
Recalling that $X_t$ is $\pi_{\epsilon}$-reversible as soon as $\mbox{\rm Law}(X_0)=\pi_{\epsilon}$, we can easily check that
\begin{eqnarray*}
\Psi_t^{(\epsilon,\overline{\epsilon})}(X_0,Z_0)&\stackrel{law}{=}&\overline{\Psi}^{(\epsilon,\overline{\epsilon})}_t(X_0,Z_0)\\&:=&\displaystyle \Ea_{t}^{\epsilon}(X_0)~Z_0+
\int_0^t~\displaystyle \Ea_{s}^{\epsilon}(X_0)~\varsigma_{\overline{\epsilon}}(X_s)~d\Wa^{\prime}_s~~\longrightarrow_{t\rightarrow\infty}~~\int_0^{\infty}~\displaystyle \Ea_{s}^{\epsilon}(X_0)~\varsigma_{\overline{\epsilon}}(X_s)~d\Wa^{\prime}_s.
\end{eqnarray*}

Given the Riccati diffusion $X_t$, the process $Z_t$ is the non-homogeneous
Ornstein-Uhlenbeck process; that is, we have that
\begin{equation}\label{OU-formula}
\begin{array}{l}
\displaystyle
\Psi^{(\epsilon,\overline{\epsilon})}_t(x,z)=\Ea_{0,t}^{\epsilon}(x)~z+\int_0^t~
\Ea_{s,t}^{\epsilon}(x)~\varsigma_{\overline{\epsilon}}(\Phi^{\epsilon}_s(x))~d\Wa^{\prime}_s\\
\\
\displaystyle\Longrightarrow~~
\Psi^{(\epsilon,0)}_t(x,z)=\Ea_{0,t}^{\epsilon}(x)~z+\int_0^t~
\Ea_{s,t}^{\epsilon}(x)~\varsigma(\Phi^{\epsilon}_s(x))~d\Wa^{\prime}_s.
\end{array}
\end{equation}
We have the following uniform moment bound and stability result.
\begin{theo}\label{theo-OU-proc}
For any $(x,z)\in (\RR_+\times\RR)$, any $n\geq 1$, any $t\geq0$, any $\epsilon\in [0,\epsilon_{2n}]$, and any $\overline{\epsilon}\in [0,1]$, we have,
\begin{eqnarray}
{\vertiii{\Psi^{(\epsilon,\overline{\epsilon})}(x,z)}_n}&\leq&c_{1,n}~(1+\vert z\vert+x^{3/2}).\label{first-unif-Psi-formula}
\end{eqnarray}
with $\epsilon_{2n}$ defined in (\ref{condi-unif-laplace-estimates-bis}). In addition, for any $(x_1,z_1),(x_2,z_2)\in (\RR_+\times\RR)$ we have
\begin{eqnarray}
{\vertiii{\Psi^{(\epsilon,\overline{\epsilon})}_t(x_1,z_1)-\Psi^{(\epsilon,\overline{\epsilon})}_t(x_2,z_2)}_n}&\leq& c_{2,n}~g(x_1,z_1,x_2,z_2)
\exp{\left[-\widehat{\lambda}_{\epsilon,2n} ~t/2\right]},~\label{first-contraction-formula}
\end{eqnarray}
with the function
$$
g(x_1,z_1,x_2,z_2):=(1+(\vert z_1\vert\wedge\vert z_2\vert))(1+x_1+x_2)^3~\frac{\vert x_1-x_2\vert}{x_1\wedge x_2}+\vert z_1-z_2\vert.
$$
\end{theo}
The proof of this theorem follows from standard stochastic calculus tools, and is in Appendix~\ref{theo-OU-proc-proof}.

\subsection{Fluctuation-Type Properties}

Several interesting results can be derived from the Ornstein-Uhlenbeck  formula (\ref{OU-formula}):

$\bullet$ The robustness properties w.r.t. the parameter $\overline{\epsilon}$ are encapsulated into the centered
Gaussian process defined by the formula
$$
\Psi^{(\epsilon,\overline{\epsilon})}_t(x,z)-\Psi^{(\epsilon,0)}_t(x,z)=\overline{\epsilon}^2~\int_0^t~
\Ea_{s,t}^{\epsilon}(x)~\left[
\frac{(\sigma_{1}/\varsigma)^2}{1+\sqrt{1+\overline{\epsilon}^2(\sigma_{1}/\varsigma)^2}}\right](\Phi^{\epsilon}_s(x))
~
\varsigma\left(\Phi^{\epsilon}_s(x)\right)~d\Wa^{\prime}_s.
$$
\begin{cor}\label{cor-1-EnKF}
Keeping the hypotheses and notation of Theorem \ref{theo-OU-proc}, we have
\begin{equation}
\begin{array}{l}
\displaystyle
{\vertiii{\Psi^{(\epsilon,\overline{\epsilon})}(x,z)-\Psi^{(\epsilon,0)}(x,z)}_{n}}
\displaystyle\leq \overline{\epsilon}^2~c_n~\left(1+x
\right).
\end{array}
\end{equation}
\end{cor}

\proof
Combining Burkholder-Davis-Gundy inequality and the generalized Minkowski
inequality for any $n\geq 1$ we have
\begin{eqnarray*}
\EE\left[\left\vert\Psi^{(\epsilon,\overline{\epsilon})}_t(x,z)-\Psi^{(\epsilon,0)}_t(x,z)\right\vert^n\right]^{2/n}
&\leq& \overline{\epsilon}^4\,n^{2}\,R^{-1}\,\int_0^t~
\EE\left[\Ea_{s,t}^{\epsilon}(x)^{n}~
\sigma_{1}^{2n}(\Phi^{\epsilon}_s(x))\right]^{2/n}~ds\\
&\leq& \overline{\epsilon}^4\,n^{2}\,R^{-1}\,\int_0^t~
\EE\left[\Ea_{s,t}^{\epsilon}(x)^{2n}\right]^{1/n}~\EE\left[
\sigma_{1}^{4n}(\Phi^{\epsilon}_s(x))\right]^{1/n}~ds.
\end{eqnarray*}
Using the uniform moment estimates (\ref{bias-n}) and the Laplace estimates (\ref{unif-laplace-estimates-bis}) we conclude that
$$
\begin{array}{l}
\displaystyle
\vertiii{\Psi^{(\epsilon,\overline{\epsilon})}(x,z)-\Psi^{(\epsilon,0)}(x,z)}_{n}~
\displaystyle\leq ~\overline{\epsilon}^2~n~\widehat{\rho}_{\epsilon,2n}\,(2R\widehat{\lambda}_{\epsilon,2n} )^{-1/2}\,\left[U\,\phi^{(\epsilon,2n)}_{\star}(x)+V\,\phi^{(\epsilon,6n)}_{\star}(x)
\right].
\end{array}
$$
This ends the proof of the corollary.\qed

$\bullet$ When $\epsilon=0=\overline{\epsilon}$ the flow $\Psi^{(0,0)}_t(x,z)$ reduces to the non-homogeneous
Ornstein-Uhlenbeck driven by the Riccati flow $\phi_t(x)$; that is, we have that
$$
\Psi^{(0,0)}_t(x,z)=\Ea_{0,t}(x)~z+\int_0^t~
\Ea_{s,t}(x)~\varsigma(\phi_s(x))~d\Wa^{\prime}_s.
$$

\begin{cor}\label{cor-2-EnKF}
For any $(x,z)\in (\RR_+\times\RR)$, any $n\geq 1$, any $t\geq0$, any $\epsilon\in [0,\epsilon_{3n}]$, and any $\overline{\epsilon}\in [0,1]$, we have,
\begin{equation}
\displaystyle
{\vertiii{\Psi^{(\epsilon,\overline{\epsilon})}(x,z)-\Psi^{(0,0)}(x,z)}_{n}}
\displaystyle~\leq~c_{n}~\left[1+
x\right]~\left[\overline{\epsilon}^2~+\epsilon~
~v^{\epsilon}_{4n}(x)~~(1+z)\right].
\end{equation}
with $\epsilon_{3n}$ defined in (\ref{condi-unif-laplace-estimates-bis}).
\end{cor}
The proof of the preceding corollary is in the Appendix~\ref{cor-2-EnKF-proof}.

\section{Ensemble Kalman-Bucy Filters}\label{EnKf-sec}

Because of their practical importance, this section is dedicated to the illustration of our main results within the {\tt EnKF} framework. We consider the filtering problem introduced in
(\ref{lin-Gaussian-diffusion-filtering}).

\subsection{A Class of McKean-Vlasov Diffusions} \label{McKean-Vlasov-sec}

For any probability measure $\eta$ on $\RR$ we let $\Pa_{\eta}$ be the $\eta$-variance
$$
\Pa_{\eta}:=\eta\left([\theta-\eta(\theta)]^2\right)\quad\mbox{\rm and}\quad
2\,\Qa_{\eta}:=R\,\Pa_{\eta}^{-1} -\Pa_{\eta}\,S~\quad\mbox{\rm with}\quad S:=B~\Sigma^{-1}B,
$$
as soon as $\Pa_{\eta}>0$, with the identity function $\theta(x):=x$.
We now consider three different classes of conditional nonlinear McKean-Vlasov type diffusion processes
\begin{eqnarray}
(1)\qquad d\overline{X}_t&=&A~\overline{X}_t~dt~+~R^{1/2}~d\overline{W}_t+\Pa_{\overline{\eta}_t}~B~\Sigma^{-1}~\left[dY_t-\left(B\overline{X}_tdt+\Sigma^{1/2}~
d\overline{V}_{t}\right)\right];\nonumber\\
(2)\qquad  d\overline{X}_t&=&A~\overline{X}_t~dt~+~R^{1/2}~d\overline{W}_t+\Pa_{\overline{\eta}_t}~B ~\Sigma^{-1}~\left[dY_t-B\left(\frac{\overline{X}_t+
\overline{\eta}_t(\theta)}{2}\right)~dt\right];\nonumber\\
(3)\qquad  d\overline{X}_t&=&\left[A\overline{X}_t+\Qa_{\overline{\eta}_t}~\left(\overline{X}_t-\overline{\eta}_t(\theta)\right)\right]~dt+\Pa_{\overline{\eta}_t}~B ~\Sigma^{-1}~\left[dY_t-B~
\overline{\eta}_t(\theta)~dt\right].
\label{Kalman-Bucy-filter-nonlinear-ref}
\end{eqnarray}
In all cases  $(\overline{W}_t, \overline{V}_t,\overline{X}_0)$ are independent copies of $(\mathscr{W}_t, \mathscr{V}_t,\mathscr{X}_0)$ (thus independent of
 the signal and the observation path) and
 \begin{equation}\label{def-nl-cov}
\overline{\eta}_t=
\mbox{\rm Law}(\overline{X}_t~|~\Ya_t).
\end{equation}
These diffusions are time-varying Ornstein-Uhlenbeck processes \cite{DelMoral/Tugaut:2016} and consequently $ \overline{\eta}_t$ is Gaussian; see also \cite{Bishop/DelMoral:2016}.
These Gaussian distributions have the same conditional mean $\mathscr{M}_t=\overline{\eta}_t(\theta)$ and the same conditional variance $\mathscr{P}_t=\Pa_{\eta_t}=\Pa_{\overline{\eta}_t}$. They satisfy
the Kalman-Bucy filter
$$
d\mathscr{M}_t=\frac{1}{2}~\partial \Lambda(\mathscr{P}_t)~\mathscr{M}_t~dt+\mathscr{P}_t~B \Sigma^{-1}~d\mathscr{Y}_t\quad\mbox{\rm with the Riccati equation}\quad
\partial_t\mathscr{P}_t=\Lambda\left(\mathscr{P}_t\right).
$$

For a more detailed discussion on the origins, and the non-uniqueness of these nonlinear McKean interpretations, see~\cite{evensen2009book,sakov2008deterministic,Reich2013,Taghvaei2016ACC}. In particular, the case (1) corresponds to the limiting object in the continuous-time version of the `vanilla' {\tt EnKF} \cite{evensen2009book}; while (2) is the limiting continuous-time object of the (perhaps confusingly named) `deterministic' {\tt EnKF} of \cite{sakov2008deterministic}, see also \cite{Reich2013,Yang2016}; and (3) is a fully deterministic optimal-transport-inspired version of an ensemble Kalman-Bucy filter \cite{Reich2013,Taghvaei2016ACC}.

\subsection{Mean Field Particle Interpretations}\label{mean-field-EnKF-sec}

Ensemble Kalman-Bucy filters ({\tt EnKF}) coincide with the mean-field particle interpretation of the nonlinear diffusion processes defined in \eqref{Kalman-Bucy-filter-nonlinear-ref}. To be more precise, let $(\overline{W}^i_t,\overline{V}^i_t,\overline{X}_0^i)_{1\leq i\leq N+1}$ be $(N+1)$ independent copies of $(\overline{W}_t,\overline{V}_t,\overline{X}_0)$.

The {\tt EnKF} associated with the first class $(1)$ of nonlinear process $\overline{X}_t$ defined in (\ref{Kalman-Bucy-filter-nonlinear-ref})  is given by the Mckean-Vlasov type interacting diffusion process
\begin{eqnarray}\label{fv1-3}
d\overline{X}_t^i=A~\overline{X}_t^i~dt+R^{1/2}~d\overline{W}_t^i+\widehat{\mathscr{P}}_t~B \Sigma^{-1}\left[d\mathscr{Y}_t-\left(B ~\overline{X}_t^i~ dt+\Sigma^{1/2}~d\overline{V}^i_{t}\right)\right],
\end{eqnarray}
with $1\leq i\leq N+1$, $N\geq1$, and the rescaled particle variance
\begin{equation}\label{fv1-3-2}
\begin{array}{l}
\displaystyle \widehat{\mathscr{P}}_t:=\left(1+N^{-1}\right)\,\Pa_{\overline{\eta}^{N}_t}
~~\quad\mbox{\rm with}\quad~~
\displaystyle\overline{\eta}^{N}_t:=(N+1)^{-1}\sum_{1\leq i\leq N+1}\delta_{\overline{X}_t^i}.
\end{array}\end{equation}
Let $\widehat{\mathscr{M}}_t=\overline{\eta}^{N}_t(\theta)$ be the particle estimate of the conditional mean $\mathscr{M}_t$.
From~\cite{DelMoral/Tugaut:2016}, and via the representation theorem (Theorem 4.2~\cite{Karatzas/Shreve:1991}; see also~\cite{Doob:1990}), there exists a filtered probability space enlargement under which we have
\begin{eqnarray}
d\widehat{\mathscr{P}}_t&=&\Lambda(\widehat{\mathscr{P}}_t)~dt+\sigma_{\epsilon}(\widehat{\mathscr{P}}_t)~d\Wa_t,
\nonumber\\
d\widehat{\mathscr{M}}_t&=&\frac{1}{2}~\partial \Lambda(\widehat{\mathscr{P}}_t)~\widehat{\mathscr{M}}_t~dt+\widehat{\mathscr{P}}_t~B \Sigma^{-1}~d\mathscr{Y}_t+ \overline{\epsilon}~\sigma_{1}(\widehat{\mathscr{P}}_t)~d\Wa^{\prime}_t,
 \quad\label{EnKF-1}
\end{eqnarray}
with the parameters
$$
 \overline{\epsilon}:=\frac{\epsilon}{\sqrt{\epsilon^2+4}} ~~\quad\mbox{\rm and}\quad~~ (U,V)=(R,S).
$$
In the above display $(\Wa_t,\Wa^{\prime}_t)$ stands for a $2$-dimensional Wiener process and
$$
\epsilon:=\frac{2}{\sqrt{N}}~~\Longrightarrow ~~ \overline{\epsilon}=\frac{1}{\sqrt{N+1}}.
$$
In the same vein,  the {\tt EnKF} associated with the second class $(2)$ of  nonlinear process $\overline{X}_t$
discussed in (\ref{Kalman-Bucy-filter-nonlinear-ref}) is given by the Mckean-Vlasov type interacting diffusion process
\begin{eqnarray}\label{fv1-3-bis}
d\overline{X}_t^i=A~\overline{X}_t^i~dt+R^{1/2}~d\overline{W}_t^i+\widehat{\mathscr{P}}_t~B \Sigma^{-1}\left[d\mathscr{Y}_t-2^{-1}B\left(\overline{X}_t^i+\overline{\eta}^N_t(\theta)\right) dt\right],
\end{eqnarray}
with $1\leq i\leq N+1$ and the rescaled particle variance
\begin{equation}\label{fv1-3-2-bis}
\begin{array}{l}
\displaystyle \widehat{\mathscr{P}}_t:=\left(1+N^{-1}\right)\,\Pa_{\overline{\eta}^{N}_t}
~~\quad\mbox{\rm with}\quad~~
\displaystyle\overline{\eta}^{N}_t:=(N+1)^{-1}\sum_{1\leq i\leq N+1}\delta_{\overline{X}_t^i}.
\end{array}
\end{equation}
Let $\widehat{\mathscr{M}}_t=\overline{\eta}^{N}_t(\theta)$ be the particle estimate of the conditional mean $m_t$.
Arguing as above, there exists a filtered probability space enlargement under which we find
\begin{eqnarray}
d\widehat{\mathscr{P}}_t&=&\Lambda(\widehat{\mathscr{P}}_t)~dt+\sigma_{\epsilon}(\widehat{\mathscr{P}}_t)~d\Wa_t,\nonumber\\
 d\widehat{\mathscr{M}}_t&=&\frac{1}{2}~\partial \Lambda(\widehat{\mathscr{P}}_t)~\widehat{\mathscr{M}}_t~dt+\widehat{\mathscr{P}}_t~B \Sigma^{-1}~d\mathscr{Y}_t+\sigma_{\overline{\epsilon}}(\widehat{\mathscr{P}}_t)~d\Wa^{\prime}_t\quad
\mbox{\rm with}
 \quad (U,V)=(R,0)\label{EnKF-2}.
 \end{eqnarray}
In the last case (3), the particle mean  $\widehat{\mathscr{M}}_t=\overline{\eta}^{N}_t(\theta)$ and the particle variance $\widehat{\mathscr{P}}_t$ associated with the McKean interpretation (3) discussed in (\ref{Kalman-Bucy-filter-nonlinear-ref}) satisfy exactly the same equations as the Kalman-Bucy filter with the associated deterministic Riccati equation \cite{Bishop/DelMoral:2016}. The ``randomness'' only comes from the initial conditions. Thus, the stability analysis of this last class of models (3) reduces to the one of the Kalman-Bucy filter and the associated Riccati equation. This reduction is true in any multidimensional filtering setting also. Several exponential estimates, in any dimension, can be found in the article of~\cite{Bishop/DelMoral:2016,bd-CARE}; see also \cite{deWiljes2018} for analysis and applications of this approach in the nonlinear filtering setting. The fluctuation analysis of this third class of {\tt EnKF} model can also be developed easily by combining the Lipschitz-type estimates w.r.t. the initial state presented in~\cite{Bishop/DelMoral:2016,bd-CARE}, with conventional sample mean error estimates based on independent copies of the initial values, see for instance~\cite{del2017wishart} for $\chi^2$-type estimates associated with sample covariance estimates.

\begin{figure}[!ht]
	\centering\resizebox*{0.75\textwidth}{0.25\textheight}{\includegraphics{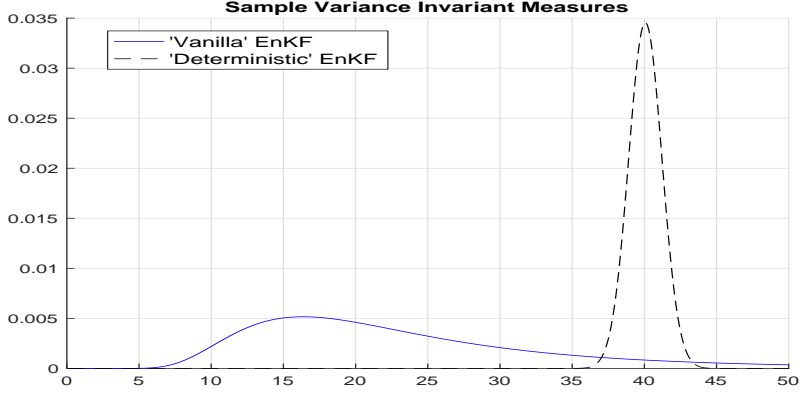}}
	\caption{The invariant measure of the sample variance of the `vanilla' {\tt EnKF} in case (1), versus that of the `deterministic' {\tt EnKF} in case (2).}
	\label{fig:invariantmeasureEnKFvsDEnKF}
\end{figure}

The invariant measure of $\widehat{\mathscr{P}}_t$ in (\ref{fv1-3-2}), (\ref{EnKF-1}) associated with the `vanilla' {\tt EnKF} in case (1) is given by (\ref{iota-pi-1}) with $(U,V)=(R,S)$. Similarly, the invariant measure of $\widehat{\mathscr{P}}_t$ in (\ref{fv1-3-2-bis}), (\ref{EnKF-2}) associated with the `deterministic' {\tt EnKF} in case (2) is given by (\ref{iota-pi-2}) with $(U,V)=(R,0)$. As an illustrative example, take $A=20$ (i.e. the underlying signal model is highly unstable), $R=S=1$ and $N=6\Rightarrow\epsilon=2/\sqrt{6}$. In Figure \ref{fig:invariantmeasureEnKFvsDEnKF} we compare the invariant measure for the flow of the sample variance in each case.

Of course, we notice in Figure \ref{fig:invariantmeasureEnKFvsDEnKF} the heavy tails of the invariant measure (\ref{iota-pi-1}) for the `vanilla' {\tt EnKF} sample variance, and conversely the Gaussian-type tails in the case (\ref{iota-pi-2}) of the `deterministic' {\tt EnKF}. Note also the positioning of the mode/mean in each case. In case (1) of the `vanilla' {\tt EnKF}, $n$-th order moments exist only when $(2n-4)/N$ is strictly less than one (in this case for $n<5$); while all moments exist in case (2) for the `deterministic' {\tt EnKF}.

\subsection{Fluctuation and Stability Properties}\label{stability-EnKF-sec}

This section is dedicated to the stochastic analysis of the diffusion processes discussed in (\ref{EnKF-1}) and (\ref{EnKF-2}).
We set
$$
\mathscr{Z}_t:=\mathscr{M}_t-\mathscr{X}_t\quad\mbox{\rm and}\quad
\widehat{\mathscr{Z}}_t:=
\widehat{\mathscr{M}}_t-\mathscr{X}_t.
$$
By the representation theorem, the evolution equations associated with the three different classes of McKean interpretation
discussed in (\ref{Kalman-Bucy-filter-nonlinear-ref}) are given by

\begin{eqnarray}
d\widehat{\mathscr{P}}_t&=&\Lambda(\widehat{\mathscr{P}}_t)~dt+\sigma_{\epsilon}(\widehat{\mathscr{P}}_t)~d\Wa_t,\nonumber\\
d\widehat{\mathscr{Z}}_t&=&\frac{1}{2}~\partial \Lambda(\widehat{\mathscr{P}}_t)~\widehat{\mathscr{Z}}_t~dt+
\varsigma^2_{\overline{\epsilon}}(\widehat{\mathscr{P}}_t)~d\Wa^{\prime}_t,
\label{class-EnKF}
 \end{eqnarray}
 with the diffusion functions $\sigma_{\epsilon}$ and $\varsigma^2_{\overline{\epsilon}}$ defined in (\ref{def-sigma-Lambda}) and
 (\ref{def-intro-Z}) and the respective parameter in each case given by
 $$
 (1)\qquad (U,V)= (R,S)\qquad (2)\qquad (U,V)= (R,0)\quad \mbox{\rm and}\quad (3)\qquad (U,V)= (0,0).
 $$
This shows that the bivariate process
 $(\widehat{\mathscr{P}}_t,\widehat{\mathscr{Z}}_t)$ is a
 two dimensional diffusion
driven by  independent Wiener processes with a first component that doesn't depend on the second.
These models are clearly encapsulated in the class of 2-dimensional
stochastic processes $(X_t,Z_t)$ on $(\RR_+\times\RR)$ discussed in (\ref{mp-ref}) and (\ref{def-intro-Z}).
Also observe that
$$
\widehat{\mathscr{Z}}_t-\mathscr{Z}_t=\widehat{\mathscr{M}}_t-\mathscr{M}_t.
$$
This shows that the analysis of the {\tt EnKF} performance reduces to the convergence analysis of the
processes $\widehat{\mathscr{Z}}_t$  towards the process $\mathscr{Z}_t$ as the parameter $\epsilon$ (and thus $\overline{\epsilon}$) tends
to $0$. The stochastic analysis of the stochastic flows associated with the pair of diffusion processes $(\widehat{\mathscr{P}}_t,\widehat{\mathscr{Z}}_t)$ is developed in Section~\ref{2-d-section}.
Last, but not least observe that if $\overline{\epsilon}=0$, the process
$\widehat{\mathscr{Z}}_t$ reduces to the re-centered Kalman--Bucy filter
$$
\widecheck{\mathscr{Z}}_t=\widecheck{\mathscr{M}}_t-\mathscr{X}_t.
$$
In the above display, the Kalman-Bucy filters  $\widecheck{\mathscr{M}}_t$   are defined as $\mathscr{M}_t$ by replacing the solution $\mathscr{P}_t$ of the Riccati equation by the stochastic approximation $\widehat{\mathscr{P}}_t$.

\begin{cor}\label{cor-ref-int}
Assume that $\widehat{\mathscr{M}}_0=\mathscr{M}_0$ and $\widehat{\mathscr{P}}_0=\mathscr{P}_0$.
In this situation, for any $n\geq 1$ there exists some parameter $\epsilon_n\in \RR_+$ such that for any $\epsilon\in [0,\epsilon_n]$
$$
\vertiii{\widehat{\mathscr{M}}-\mathscr{M}}_{n}
\displaystyle\leq ~\frac{c_{1,n}}{\sqrt{N}}~\quad\mbox{and}\quad  \vertiii{\widecheck{\mathscr{M}}-\mathscr{M}}_{n}\leq \frac{c_{2,n}}{N+1}.
$$

\end{cor}

The above estimates are a direct consequence of Corollary~\ref{cor-2-EnKF} and Corollary~\ref{cor-1-EnKF}. From Theorem~\ref{theo--1intro}, we also have the uniform estimate,
$$
	\sqrt{N}~{\vertiii{ \widehat{\mathscr{P}}-\mathscr{P}}_n}\leq c_{1,n}
$$

Let $(\widehat{\mathscr{M}}_t,\widehat{\mathscr{P}}_t)$, and $(\widehat{\mathscr{M}}^{\,\prime}_t,\widehat{\mathscr{P}}^{\,\prime}_t)$ denote the {\tt EnKF} sample mean and sample variance starting from two possibly different initial conditions. In this situation we have
$$
\left.
\begin{array}{rcl}
\widehat{\mathscr{Z}}_t&=&(\widehat{\mathscr{M}}_t-\mathscr{X}_t)\\
\widehat{\mathscr{Z}}^{\,\prime}_t&=&(\widehat{\mathscr{M}}^{\,\prime}_t-\mathscr{X}_t)\end{array}
\right\}
~~\Longrightarrow ~~\widehat{\mathscr{Z}}_t-\widehat{\mathscr{Z}}^{\,\prime}_t=\widehat{\mathscr{M}}_t-\widehat{\mathscr{M}}_t^{\,\prime}.
$$
In what follows, $\epsilon=2/\sqrt{N}$, and $\widehat{\lambda}_{\epsilon} \,>0$ and $ \widehat{\lambda}_{\epsilon,\kappa} >0$  are the parameters defined in  (\ref{def-lambda-epsilon-iota}).

\begin{cor}\label{enkf-intro-1-proof}
For any $n\geq 1$, any time horizon $t\geq 0$, we have
$$
N/2> \zeta
\vee\left(\left(\imath_{n-2}+2\right)\,\overline{V}\right)
~~\Longrightarrow~~
 \vertiii{\widehat{\mathscr{P}}_t-\widehat{\mathscr{P}}_t^{\,\prime}}_n\leq c_n~ \exp{\left(
 -\widehat{\lambda}_{\epsilon} \,t\, \right)}.
$$

 \end{cor}
The above Riccati semigroup contraction estimates are direct consequence of Theorem~\ref{theo-4-intro}.
 When $(U,V)=(R,0)$ we find the estimate (\ref{enkf-intro-1}).

 The next corollary is a consequence of \eqref{first-contraction-formula} and gives the estimate \eqref{enkf-intro-2} when $(U,V)=(R,0)$.
 \begin{cor}\label{enkf-intro-2-proof}
For any $n\geq 1$ there exists some finite constants $c_n$ such that for
any  $t\geq 0$ we have
$$
N > 2\zeta_{2n}
~~\Longrightarrow~~
\vertiii{\widehat{\mathscr{M}}_t-\widehat{\mathscr{M}}_t^{\,\prime}}_n\leq c_{n}~\exp{\left[-\widehat{\lambda}_{\epsilon,2n}\,t/2\right]}.
$$
\end{cor}

\begin{figure}[!ht]
	\centering
	\resizebox*{0.32\textwidth}{0.2\textheight}{\includegraphics{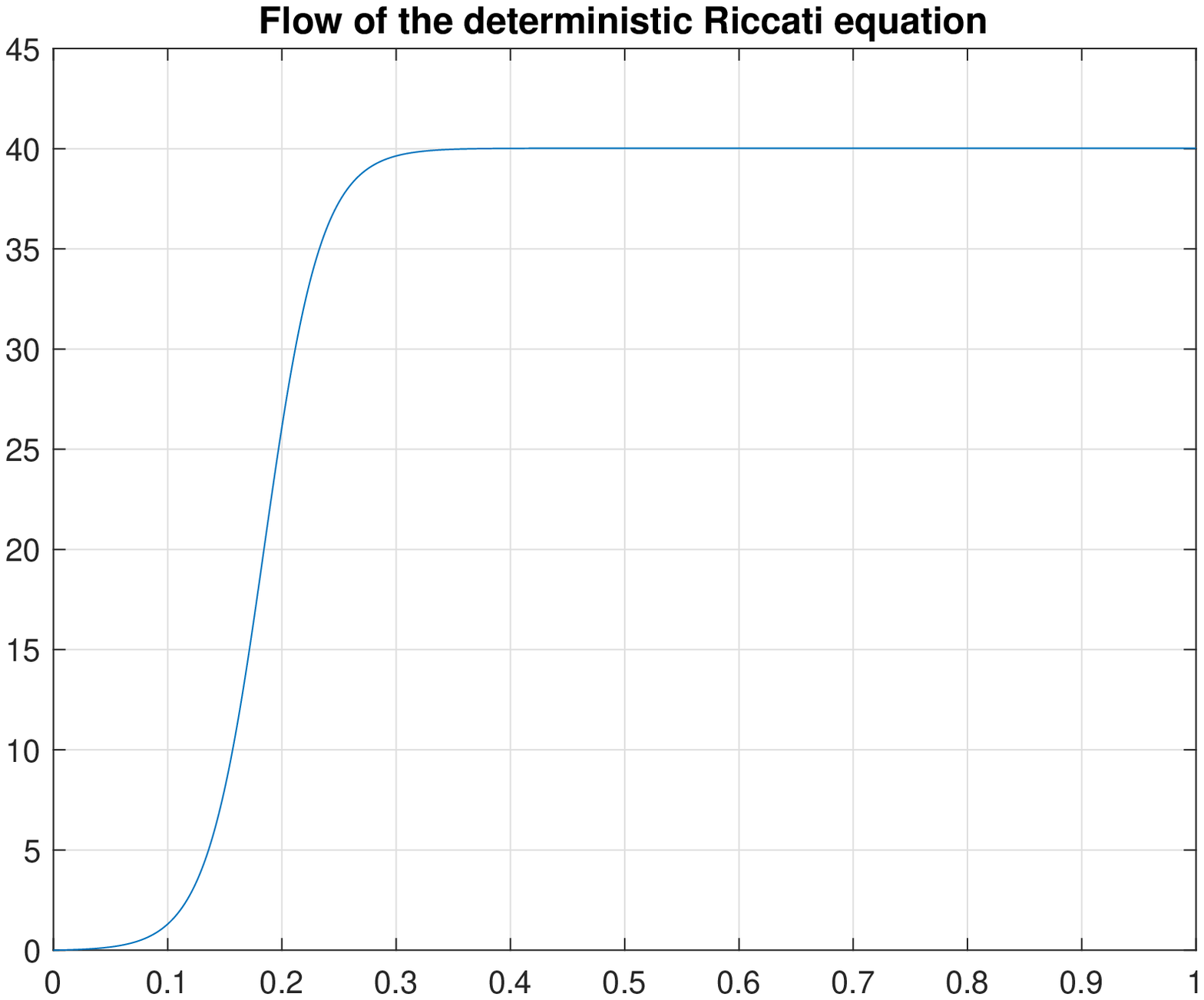}}
	\resizebox*{0.32\textwidth}{0.2\textheight}{\includegraphics{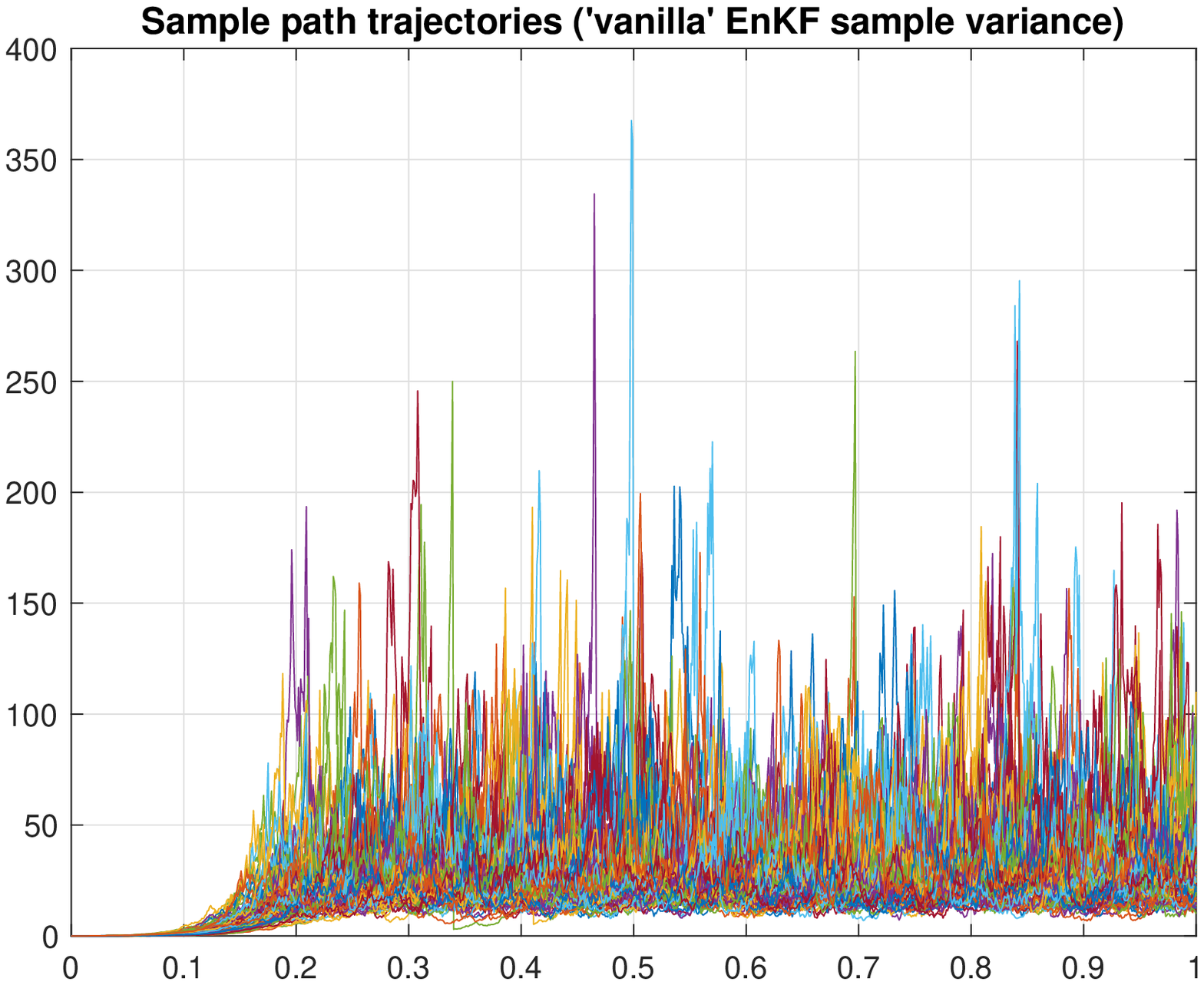}}
	\resizebox*{0.32\textwidth}{0.2\textheight}{\includegraphics{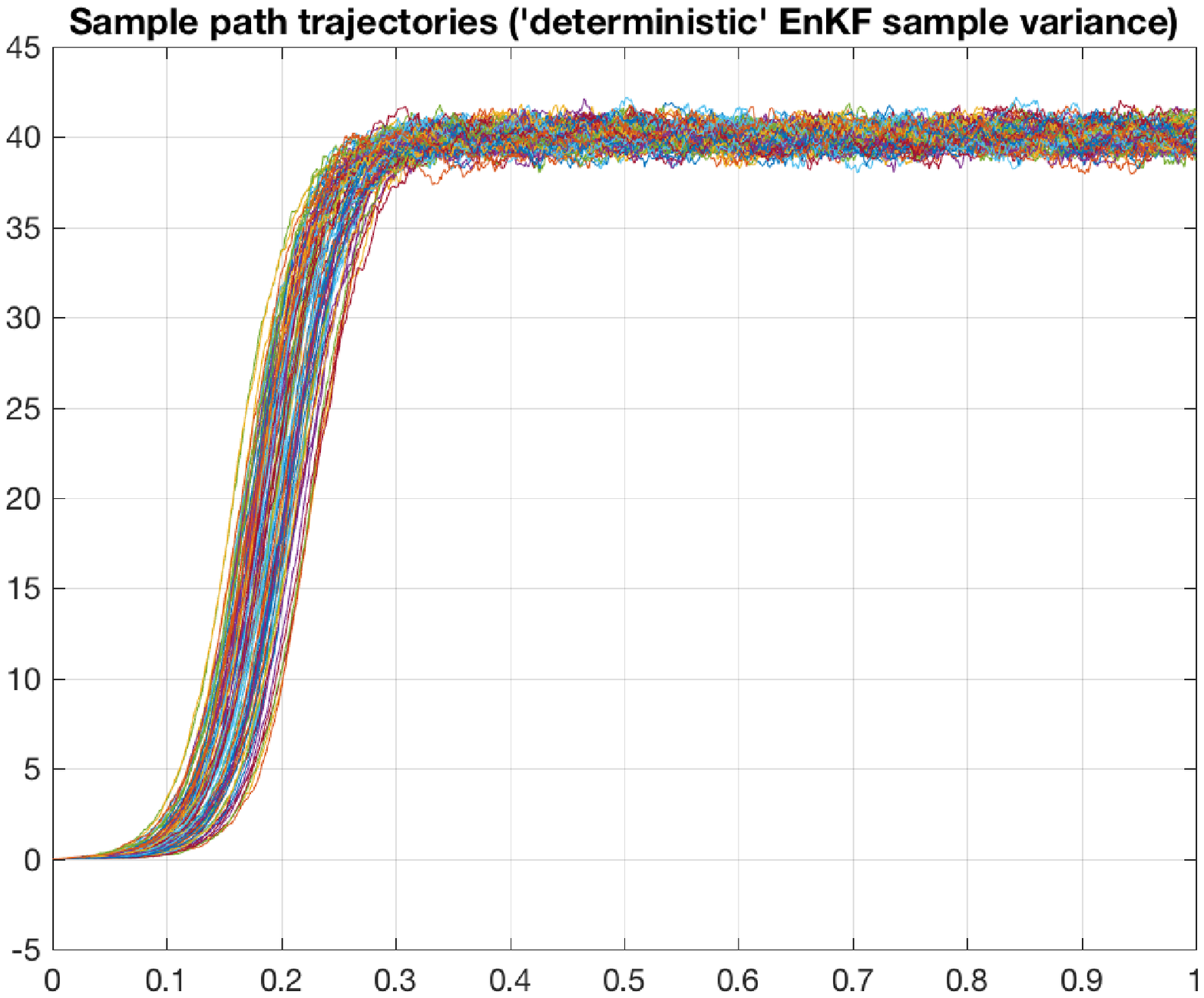}}
	\caption{Flow of the deterministic Riccati equation, $100$ sample paths of the `vanilla' {\tt EnKF} sample variance of case (1), and $100$ sample paths of the `deterministic' {\tt EnKF} sample variance of case (2).}
	\label{fig:detRiccatiflow}
\end{figure}

We consider an illustration of the fluctuation and stability properties of the sample variance in the different {\tt EnKF} variants. Consider again the model leading to Figure \ref{fig:invariantmeasureEnKFvsDEnKF}, and let $\widehat{\mathscr{P}}_0=0$. The deterministic Riccati flow ($\epsilon=0$) with the chosen model parameters is given in Figure \ref{fig:detRiccatiflow}, along with $100$ sample paths of the sample variances for both the `vanilla' {\tt EnKF} and the `deterministic' {\tt EnKF}.

\begin{figure}[!ht]
	\centering
	\resizebox*{0.495\textwidth}{.35\textheight}{\includegraphics{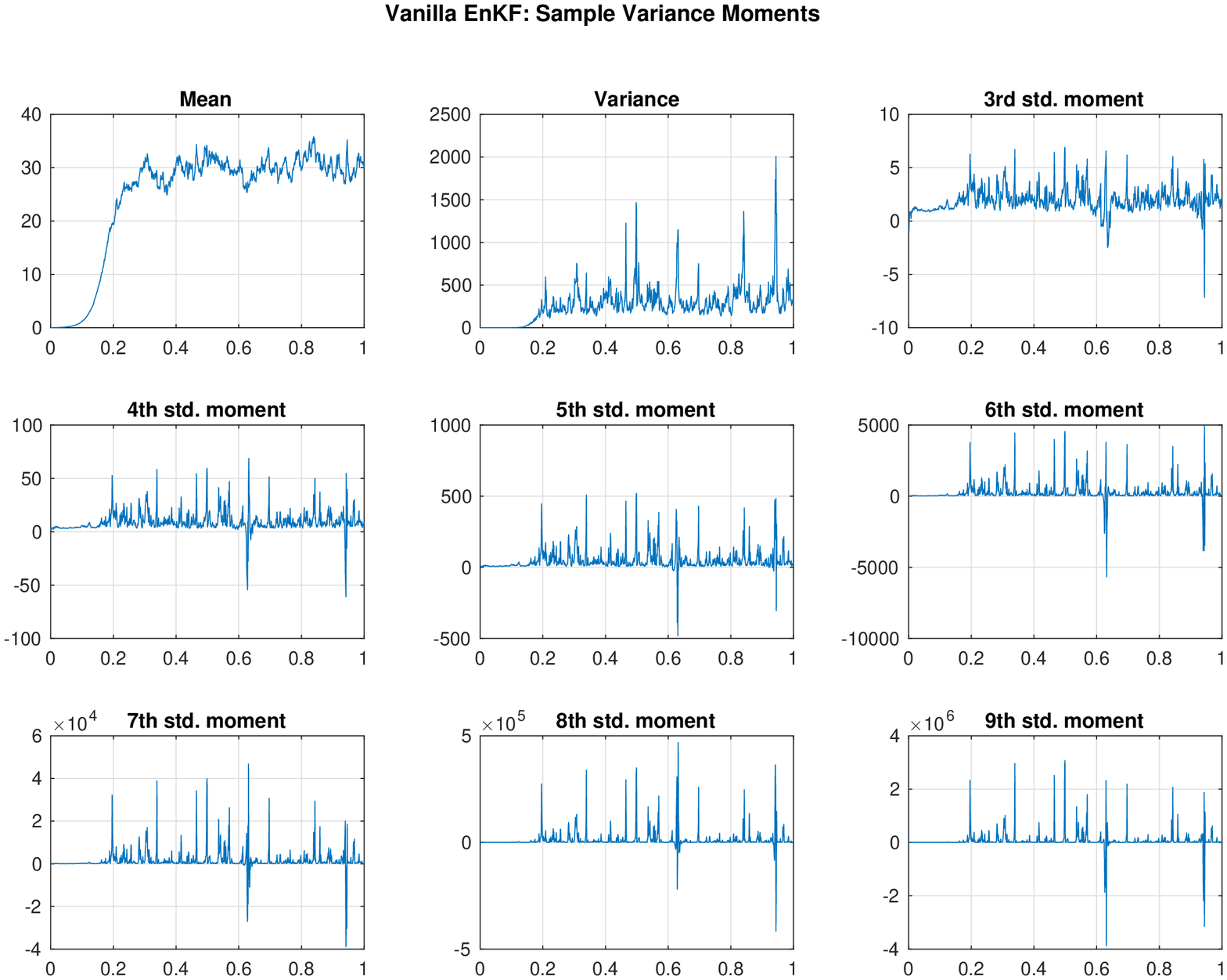}}
	\resizebox*{0.495\textwidth}{.35\textheight}{\includegraphics{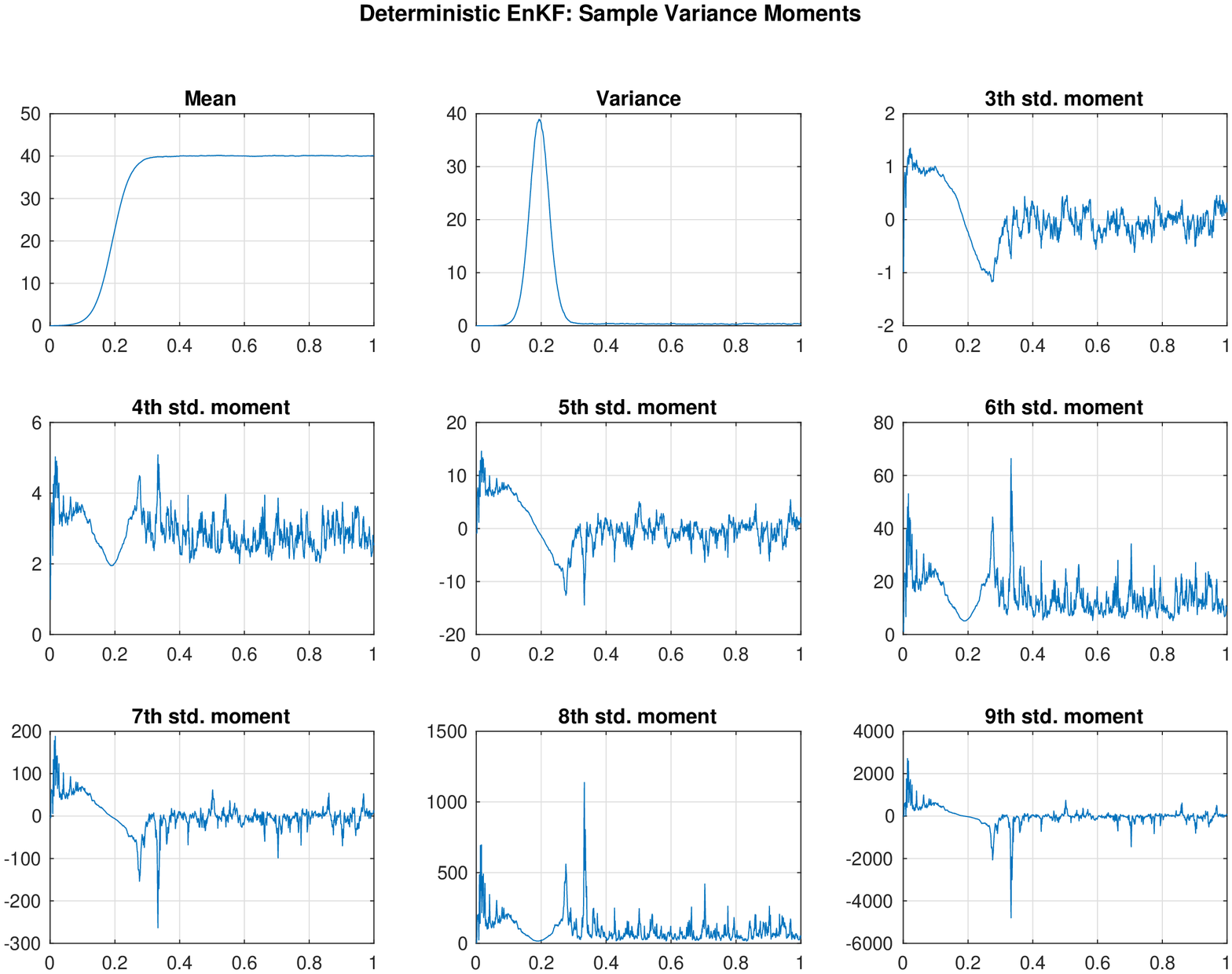}}\vspace{-.5cm}
	\caption{Flow of the sample variance moments for the `vanilla' {\tt EnKF} and the `deterministic' {\tt EnKF}.}
	\label{fig:momentFlow}
\end{figure}

Note in Figure \ref{fig:detRiccatiflow} the drastically reduced fluctuations in `deterministic' {\tt EnKF} sample variance sample paths. In Figure \ref{fig:momentFlow} we plot the flow of the first two central moments and the $3rd$ through the $9th$ standardised central moments for both the `vanilla' {\tt EnKF} sample variance, and the `deterministic' {\tt EnKF} sample variance distribution. Recall that $N=6$ in this case and we expect moments of the `vanilla' {\tt EnKF} sample variance in case (1) to exist up to $n=4$ with $n=5$ the boundary case; while all moments exist for the `deterministic' {\tt EnKF} of case (2).

We note in Figure \ref{fig:momentFlow} that the sample variance moments for the `vanilla' {\tt EnKF} in case (1) begin to destabilize around the $5th/6th$ moments as expected. Importantly, the mean of the sample variance for the `vanilla' {\tt EnKF} is very negatively biased in this case, while the mean of the `deterministic' {\tt EnKF} in case (2) is particularly accurate. We note also the very large variance in the sample variance for the 'vanilla' {\tt EnKF}.

Lastly, we note that $\widehat{\mathscr{P}}_t$ in (\ref{fv1-3-2}), (\ref{EnKF-1}) associated with the `vanilla' {\tt EnKF} of case (1), has non-globally Lipschitz coefficients. In particular, the drift is quadratic, while the diffusion has a polynomial growth of order $3/2$. It follows by \cite{hutzenthaler2011strong} that a basic Euler time-discretization may blow up, irregardless of the boundedness properties of the diffusion.

\appendix

\section{Proof of the Fluctuation Theorem}

This section is concerned with the proof of the moment and fluctuation estimates in Theorem \ref{theo--1intro}.

\subsection{Proof of the Moments Estimates (\ref{bias-n})}\label{defi-phi-n-eps-section}

\begin{lem}\label{lemma-1}
Assume that the flow $t\mapsto \varphi_t(x)$ satisfies a Riccati-type inequality of the form
$$
\partial_t\varphi_t(x)\leq 2A~\varphi_t(x)+R-S~\varphi_t(x)^2
$$
for any $t\geq 0$ and any $x\in \RR_+$.
In this situation, for any time horizon $t\geq 0$ and $x,y\geq 0$ we have
\begin{equation}\label{comparison-1}
\varphi_t(x)\leq \phi_{t}(y)+\exp{\left(2\int_0^t (A-S\phi_{s}(y))~ds\right)}~(x-y).
\end{equation}
where $\phi_{t}$ denotes the semigroup of the Riccati equation (\ref{nonlinear-KB-Riccati}).
\end{lem}

\proof
We have
\begin{eqnarray*}
\partial_t\left[\varphi_t(x)-\phi_{t}(y)\right]&\leq& 2A~\left[\varphi_t(x)-\phi_{t}(y)\right]-S~\left[\varphi_t(x)^2-\phi_{t}(y)^2\right]\\
&=&2\left(A-S~\phi_{s}(y)\right)~\left[\varphi_t(x)-\phi_{t}(y)\right]-S~\left[\varphi_t(x)-\phi_{t}(y)\right]^2\\
&\leq &2\left(A-S~\phi_{s}(y)\right)~\left[\varphi_t(x)-\phi_{t}(y)\right].
\end{eqnarray*}
This implies that
$$
\partial_t\left[\exp{\left(-2\int_0^t (A-S~\phi_{s}(y))~ds\right)}~\left[\varphi_t(x)-\phi_{t}(y)\right]\right]\leq 0\Longrightarrow (\ref{comparison-1}).
$$
This ends the proof of the lemma.
\qed

We are now in a position to prove the moments bias estimates stated in (\ref{bias-n}).

\noindent
{\bf Proof of (\ref{bias-n}):}
For any $n\geq 1$, we have
\begin{eqnarray*}
d\Phi_{t}^{\epsilon}(x)^n&=&n\left[\Phi_{t}^{\epsilon}(x)^{n-1}
(2A\Phi_{t}^{\epsilon}(x)+R-S\Phi_{t}^{\epsilon}(x)^2)~\right.\\
&&\hskip1cm\left.+\frac{n-1}{2}~\Phi_{t}^{\epsilon}(x)^{n-2}~\epsilon^2~ \Phi_{t}^{\epsilon}(x)~(U+V\Phi_{t}^{\epsilon}(x)^2)~\right]~dt\\
&&\hskip3cm+\epsilon~ n~\Phi_{t}^{\epsilon}(x)^{n-1}\sqrt{\Phi_{t}^{\epsilon}(x)~(U+V\Phi_{t}^{\epsilon}(x)^2)}~dW_t.
\end{eqnarray*}
This yields
\begin{eqnarray*}
n^{-1}~\partial_t\EE\left[\Phi_{t}^{\epsilon}(x)^n\right]&=&
2A~\EE\left[\Phi_{t}^{\epsilon}(x)^n\right]+\left(R+\frac{n-1}{2}~\epsilon^2~U\right)\EE\left[\Phi_{t}^{\epsilon}(x)^{n-1}\right]\\
&&\hskip5cm-\left(S-\frac{n-1}{2}~\epsilon^2~V\right)\EE\left[\Phi_{t}^{\epsilon}(x)^{n+1}\right]\\
&\leq &~2A~\EE\left[\Phi_{t}^{\epsilon}(x)^n\right]+\left(R+\frac{n-1}{2}~\epsilon^2~U\right)\EE\left[\Phi_{t}^{\epsilon}(x)^{n}\right]^{1-1/n}\\
&&\hskip5cm-\left(S-\frac{n-1}{2}~\epsilon^2~V\right)\EE\left[\Phi_{t}^{\epsilon}(x)^{n}\right]^{1+1/n}.
\end{eqnarray*}
Choosing $\epsilon$ small enough so that
$
(n-1)~\epsilon^2<2~S/V
$
we find that
\begin{eqnarray*}
\partial_t\EE\left[\Phi_{t}^{\epsilon}(x)^n\right]^{1/n}&=&n^{-1}~\EE\left[\Phi_{t}^{\epsilon}(x)^n\right]^{1/n-1}~\partial_t\EE\left[\Phi_{t}^{\epsilon}(x)^n\right]\\
&\leq &~2A~\EE\left[\Phi_{t}^{\epsilon}(x)^n\right]^{1/n}+\left(R+\frac{n-1}{2}~\epsilon^2~U\right)\\
&&\hskip2.5cm-\left(S-\frac{n-1}{2}~\epsilon^2~V\right)~
\left(\EE\left[\Phi_{t}^{\epsilon}(x)^{n}\right]^{1/n}\right)^2.
\end{eqnarray*}
The end of the proof of the upper bound in the l.h.s. estimate in (\ref{bias-n}) is now a consequence of  Lemma~\ref{lemma-1}.
Now we come to the proof of the r.h.s. estimate in (\ref{bias-n}).
By Jensen's inequality,
$$
\EE\left(\Phi^{-\epsilon}_t(x)^n\right)^{1/n}\geq \EE\left(\Phi^{-\epsilon}_t(x)\right)\geq\EE\left(\Phi_{t}^{\epsilon}(x)\right)^{-1}\geq \phi^{-1}_t(x).
$$
Arguing as above and using Lemma~\ref{inverse-ricc-lem}, we can also check that
$
\EE\left(\Phi^{-\epsilon}_t(x)^n\right)^{1/n}\leq 1/{\phi}^{(\epsilon,-n)}_t(x)$. This completes the proof of the r.h.s. estimate in (\ref{bias-n}).
Finally notice that
$$
1/{\phi}^{(\epsilon,-1)}_t(x)\geq \EE\left(\Phi^{-\epsilon}_t(x)\right)\geq\EE\left(\Phi_{t}^{\epsilon}(x)\right)^{-1}\Longrightarrow
{\phi}^{(\epsilon,-1)}_t(x)\leq \EE\left(\Phi_{t}^{\epsilon}(x)\right)\leq \EE\left(\Phi^{\epsilon}_t(x)^n\right)^{1/n}
$$
for any $n\geq 1$. This completes the proof.
  \qed

  \subsection{Proof of the Fluctuation Estimates (\ref{unif-est-intro})}\label{proof-theo-intro}

  The proof of the uniform estimates are based on the following perturbation lemma.
   \begin{lem}\label{lemma-2}
For any time horizon $t\geq 0$, any $x\in\RR_+$ and  $\epsilon\in \RR_+$ we have
\begin{eqnarray*}
\epsilon^{-1}\left[\Phi_{t}^{\epsilon}(x)-\phi_t(x)\right]&=&\int_0^t~\left(\partial \phi_{t-s}\right)(\Phi_{s}^{\epsilon}(x))~\left[\Phi^{\epsilon}_{s}(x)(U+V\Phi^{\epsilon}_{s}(x)^2)\right]^{1/2}~dW_s\\
&&\hskip1cm\displaystyle+\frac{\epsilon}{2}~\int_0^t~\left(\partial^2 \phi_{t-s}\right)(\Phi_{s}^{\epsilon}(x))~\left[\Phi^{\epsilon}_{s}(x)\left(U+V\,\Phi^{\epsilon}_{s}(x)^2\right)\right]~ds.
\end{eqnarray*}
\end{lem}
\proof
 By \cite[Proposition 2.2]{Bishop/DelMoral/Pathiraja:2017} we have
$$
\partial_s\phi_{t-s}(x)=-\Lambda\left(\phi_{t-s}(x)\right)=-\nabla\phi_{t-s}(x)~\Lambda(x).
$$
We fix $t$ and we use the interpolating path
$$
s\in [0,t]\mapsto \phi_{t-s}\left(\Phi_{s}^{\epsilon}(x)\right)\quad \mbox{\rm between}\quad
\Phi_{t}^{\epsilon}(Q) \quad \mbox{\rm and}\quad \phi_{t}(x).
$$
Applying Ito's formula on the interval $[0,t]$ we find that
\begin{eqnarray*}
d\phi_{t-s}\left(\Phi_{s}^{\epsilon}(x)\right)&=&-\left(\partial \phi_{t-s}\right)(\Phi_{s}^{\epsilon}(x))~\Lambda(\Phi_{s}^{\epsilon}(x))~ds+\partial \phi_{t-s}(\Phi_{s}^{\epsilon}(x))~d\Phi_{s}^{\epsilon}(x)\\
&&\hskip3cm+\frac{\epsilon^2}{2}~\left(\partial ^2\phi_{t-s}\right)(\Phi_{s}^{\epsilon}(x))~\left[\Phi^{\epsilon}_{t-s}(x)\left(U+V\,\Phi^{\epsilon}_{t-s}(x)^2\right)\right]~ds\\
&=&\frac{\epsilon^2}{2}~\left(\partial ^2\phi_{t-s}\right)(\Phi_{s}^{\epsilon}(x))~\left[\Phi^{\epsilon}_{s}(x)\left(U+V\,\Phi^{\epsilon}_{s}(x)^2\right)\right]~ds\\&&\hskip3cm+\epsilon~\partial \phi_{t-s}(\Phi_{s}^{\epsilon}(x))~\left[\Phi^{\epsilon}_s(Q)(U+V\Phi^{\epsilon}_{s}(Q)^2)\right]^{1/2}~dW_s.
\end{eqnarray*}
This ends the proof of the lemma.
\qed

We are now in a position to prove the uniform fluctuation estimates stated in Theorem~\ref{theo--1intro}.

\noindent
{\bf Proof of  (\ref{unif-est-intro}):}
For any $n\geq 2$, combining (\ref{estimate-first-derivative}) with Burkholder-Davis-Gundy and the generalized Minkowski inequalities we have
$$
\begin{array}{l}
\displaystyle\EE\left[\left\vert\int_0^t~\left(\partial \phi_{t-s}\right)(\Phi_{s}^{\epsilon}(x))~\left[\Phi^{\epsilon}_{s}(x)(U+V\Phi^{\epsilon}_{s}(x)^2)\right]^{1/2}~dW_s\right\vert^{n}\right]
\\
\\
\qquad \leq n^n~\displaystyle\EE\left[\left\vert\int_0^t~\left(\partial \phi_{t-s}\right)(\Phi_{s}^{\epsilon}(x))^2~\left[\Phi^{\epsilon}_{s}(x)(U+V\Phi^{\epsilon}_{s}(x)^2)\right]~ds\right\vert^{n/2}\right]\\
\\
\qquad\leq \displaystyle n^n~
\varpi^{2n}~\left[\int_0^t~e^{-2\lambda (t-s)}~
\left[U~
\EE\left(   \Phi^{\epsilon}_{s}(x)^{n/2}
       \right)^{2/n}+V~\EE\left(  \Phi^{\epsilon}_{s}(x)^{3n/2}
       \right)^{2/n}\right]~ds
\right]^{n/2}.
\end{array}
$$
Using the moment estimates (\ref{bias-n}) we find that
$$
\begin{array}{l}
\displaystyle\EE\left[\left\vert\int_0^t~\left(\partial \phi_{t-s}\right)(\Phi_{s}^{\epsilon}(x))~\left[\Phi^{\epsilon}_{s}(x)(U+V\Phi^{\epsilon}_{s}(x)^2)\right]^{1/2}~dW_s\right\vert^{n}\right]^{1/n}\\
\\
\qquad \leq \displaystyle \frac{n}{\sqrt{2\lambda }}~\varpi^2~\left[U~
 \phi^{(\epsilon,\frac{n}{2})}_{\star}(x)
     +V~ \phi^{(\epsilon,\frac{3n}{2})}_{\star}(x)^3\right]^{1/2}
~\leq ~\frac{n}{\sqrt{2\lambda }}~\varpi^2~\sigma_1\left(\phi^{(\epsilon,\frac{3n}{2})}_{\star}(x)\right).
\end{array}
$$
In the same vein, using (\ref{estimate-second-derivative}) we have
$$
\begin{array}{l}
\displaystyle\EE\left[\left\vert\int_0^t~\left(\partial ^2\phi_{t-s}\right)(\Phi_{s}^{\epsilon}(x))~\left[\Phi^{\epsilon}_{s}(x)\left(U+V\,\Phi^{\epsilon}_{s}(x)^2\right)\right]~ds
\right\vert^{n}\right]^{1/n}\\
\\
\qquad \leq \displaystyle \frac{\varpi^2}{\vert\varpi_-\vert}~\frac{2}{\lambda}~
\left[U~
 \phi^{(\epsilon,n)}_{\star}(x)
     +V~ \phi^{(\epsilon,3n)}_{\star}(x)^3\right]\leq \frac{\varpi^2}{\vert\varpi_-\vert}~\frac{2}{\lambda }~\sigma_1^2\left(\phi^{(\epsilon,3n)}_{\star}(x)\right),
\end{array}
$$
from which we conclude that
$$
\EE\left(\vert\VV^{\epsilon}_t(x)\vert^n\right)^{1/n}\leq \frac{\varpi^2}{\sqrt{2\lambda }}~\left[n~\sigma_1\left(\phi^{(\epsilon,\frac{3n}{2})}_{\star}(x)\right)+
 \frac{1}{\vert\varpi_-\vert}~\epsilon~\sqrt{\frac{2}{\lambda }}~\sigma_1^2\left(\phi^{(\epsilon,3n)}_{\star}(x)\right)
  \right].
  $$
This ends the proof of  the l.h.s. estimate in (\ref{unif-est-intro}). Now we come to the proof of  the second estimate in (\ref{unif-est-intro}).
Observe that
$$
\partial \sigma_1(x)=\frac{3}{2\sqrt{x}}\frac{U/3+Vx^2}{\sqrt{U+Vx^2}}\leq \frac{3}{2\sqrt{x}}~\sqrt{U+Vx^2}\le \frac{3}{2}~\left[\sqrt{U/x}+\sqrt{Vx}\right].
$$
We consider the functions
$$
\begin{array}{l}
\displaystyle0\leq \Sigma\left[x,y\right]:=\int_0^1~\left(\partial \sigma_1\right)(ux+(1-u)y)~du\leq \frac{3}{2}\left[2~\sqrt{\frac{U}{x\vee y}}+\sqrt{V~(x\vee y)} \right]\\
\\
\displaystyle\Longrightarrow\Sigma\left[x,y\right]^2\leq \left(\frac{3}{2}\right)^2~\left[4~\frac{U}{x\vee y}+V~(x\vee y)+4\sqrt{UV} \right]
\end{array}
$$
and
$$
 \frac{2}{\varpi_-}~\varpi^2~e^{-\lambda  t}\leq \Xi_{t}\left[x,y\right]:=\int_0^1~\left(\partial ^2\phi_{t}\right)\left(
ux+(1-u)y\right)~du\leq 0.
$$
With this notation, we have
\begin{eqnarray*}
\sigma_1\left[\Phi^{\epsilon}_{s}(x)\right]&=&\sigma_1\left[\phi_{s}(x)\right]+\left[\Phi^{\epsilon}_{s}(x)-\phi_{s}(x)\right]~\Sigma\left[\Phi^{\epsilon}_{s}(x),\phi_{s}(x)\right],\\
\left(\partial \phi_{t-s}\right)(\Phi_{s}^{\epsilon}(x))&=&\left(\partial \phi_{t-s}\right)(\phi_{s}(x))+\left[\Phi^{\epsilon}_{s}(x)-\phi_{s}(x)\right]~{\Xi}_{t-s}\left[\Phi^{\epsilon}_{s}(x),\phi_{s}(x)\right].
\end{eqnarray*}

By Lemma~\ref{lemma-2}
we have
$
\epsilon^{-1}\left[\VV_t^{\epsilon}(x)-\VV_t(x)\right]=\BB^{\epsilon}_t(x)+\CC^{\epsilon}_t(x)
$, with  the bias term
$$
\BB^{\epsilon}_t(x):=\frac{1}{2}~\int_0^t~\left(\partial ^2\phi_{t-s}\right)(\Phi_{s}^{\epsilon}(x))~\left[\Phi^{\epsilon}_{s}(x)\left(U+V\,\Phi^{\epsilon}_{s}(x)^2\right)\right]~ds.
$$
The centered remainder term $\CC^{\epsilon}_t(x)$ is given by
$$
\CC^{\epsilon}_t(x)=\epsilon~\CC^{1,\epsilon}_t(x)+\CC^{2,\epsilon}_t(x)+\CC^{3,\epsilon}_t(x),
$$
with
\begin{eqnarray*}
\CC^{1,\epsilon}_t(x)&:=&
\int_0^t~\VV_s^{\epsilon}(x)^2~\Sigma\left[\Phi^{\epsilon}_{s}(x),\phi_{s}(x)\right]~{\Xi}_{t-s}\left[\Phi^{\epsilon}_{s}(x),\phi_{s}(x)\right]~dW_s,\\
\CC^{2,\epsilon}_t(x)&:=&\int_0^t\left(\partial \phi_{t-s}\right)(\phi_{s}(x))~\VV_s^{\epsilon}(x)~\Sigma\left[\Phi^{\epsilon}_{s}(x),\phi_{s}(x)\right]~dW_s,\\
\CC^{3,\epsilon}_t(x)&:=&\int_0^t\sigma_1\left[\phi_{s}(x)\right]~\VV_s^{\epsilon}(x)~{\Xi}_{t-s}\left[\Phi^{\epsilon}_{s}(x),\phi_{s}(x)\right]~dW_s.
\end{eqnarray*}
Arguing as above, we have
$$
\begin{array}{l}
\displaystyle\EE\left(\vert\CC^{3,\epsilon}_t(x)\vert^{n}\right)^{1/n}
\displaystyle\leq \sqrt{\frac{2}{\lambda }}~n~v_{\epsilon,n}(x)~\frac{\varpi^{2}}{\vert\varpi_-\vert}~\sigma_1\left(\phi_{\star}(x)\right).
\end{array}
$$
Combining Cauchy-Schwartz inequality with the moment estimates (\ref{bias-n}) and (\ref{estimate-first-derivative}) we also check that
$$
\begin{array}{l}
\displaystyle\EE\left(\vert \CC^{2,\epsilon}_t(x)\vert^{n}\right)^{2/n}\\
\\
\displaystyle\leq \left(\frac{3n}{2}\right)^{2} \varpi^{4}~\EE\left(\vert \int_0^te^{-2\lambda  (t-s)} ~\VV_s^{\epsilon}(x)^2~\left[\left(4U~\phi^{-}_{\star}(x)+V
\phi_{\star}(x)+4\sqrt{UV} \right)+V~\Phi^{\epsilon}_{s}(x)\right]~ds\vert^{n/2}\right)^{2/n}\\
\\
\displaystyle\leq \left(\frac{3n}{2}\right)^{2}~\frac{2}{\lambda }~ \varpi^{4}~v_{\epsilon,2n}(x)^2~\left[U~\phi^{-}_{\star}(x)
+\sqrt{UV}+V~(\phi^{(\epsilon,n)}_{\star}(x)+\phi_{\star}(x))/4\right].
\end{array}
$$
The monotone properties (\ref{monotone-epsilon-n}) yield
$$
\displaystyle\EE\left(\vert \CC^{2,\epsilon}_t(x)\vert^{n}\right)^{1/n}
\displaystyle\leq \left(\frac{3}{2}\right)~\sqrt{\frac{2}{\lambda }}~n \varpi^{2}~v_{\epsilon,2n}(x)~\left[U~\phi^{-}_{\star}(x)+\sqrt{UV} +V~\phi^{(\epsilon,n)}_{\star}(x)/2\right]^{1/2}.
$$
In the same vein, we check that
$$
\begin{array}{l}
\displaystyle\EE\left(\vert \CC^{1,\epsilon}_t(x)\vert^{n}\right)^{2/n}\\
\\
\displaystyle\leq  \frac{(3n)^2}{\varpi_-^2}~\varpi^4~\EE\left(\vert \int_0^t~e^{-2\lambda  (t-s)}~\VV_s^{\epsilon}(x)^4\left[\left(4\,U\,\phi^{-}_{\star}(x)+V\,\phi_{\star}(x)+4\sqrt{UV}\right)+V\,\Phi^{\epsilon}_{s}(x) \right]
\,ds\vert^{n/2}\right)^{2/n}\\
\\
\displaystyle\leq  \frac{(3n)^2}{\varpi_-^2}~\varpi^4~\frac{2}{\lambda }~v_{\epsilon,4n}^4~\left[
\left(U~\phi^{-}_{\star}(x)+\sqrt{UV}\right)+V~\phi^{(\epsilon,n)}_{\star}(x)/2
\right].
\end{array}
$$
This yields the estimate
$$
\EE\left(\vert \CC^{1,\epsilon}_t(x)\vert^{n}\right)^{1/n}\leq  \frac{3n}{\vert\varpi_-\vert}~\varpi^2~\sqrt{\frac{2}{\lambda }}~v_{\epsilon,4n}^2~\left[
U~\phi^{-}_{\star}(x)+\sqrt{UV}+V~\phi^{(\epsilon,n)}_{\star}(x)/2
\right]^{1/2}.
$$
We conclude that
$$
\begin{array}{l}
\displaystyle
\EE\left(\left\vert \WW_t^{\epsilon}(x)\right\vert^n\right)^{1/n}
\displaystyle\leq \frac{\varpi^2}{\vert\varpi_-\vert}~\sqrt{\frac{2}{\lambda }}\left[~\frac{1}{\sqrt{2\lambda }}~\sigma_1^2\left(\phi^{(\epsilon,3n)}_{\star}(x)\right)+ ~n~v_{\epsilon,n}(x)~\sigma_1\left(\phi_{\star}(x)\right)\right]\\
\\
\displaystyle+  \frac{3n}{2}~\frac{\varpi^2}{\vert\varpi_-\vert}~\sqrt{\frac{2}{\lambda }}~~\left[
U~\phi^{-}_{\star}(x)+\sqrt{UV}+V~\phi^{(\epsilon,n)}_{\star}(x)/2
\right]^{1/2}~\left[2\epsilon~v_{\epsilon,4n}^2+v_{\epsilon,2n}(x)~\vert\varpi_-\vert\right].
\end{array}
$$
This ends the proof of the second estimate in (\ref{unif-est-intro}). Now we come to the proof of the uniform bias estimate
stated in the r.h.s. of (\ref{unif-est-intro}). Using (\ref{estimate-second-derivative}) we have
\begin{eqnarray*}
0\leq \overline{\Sigma}_{t}\left[x,y\right]&:=&\int_0^1~\partial \left(\sigma_1^2\right)(x)\left(
ux+(1-u)y\right)~du\leq  U+3V~(x+y)^2,\\
0\leq \overline{\Xi}_{t}\left[x,y\right]&:=&\int_0^1~\left(\partial ^2\phi_{t}\right)\left(
ux+(1-u)y\right)~du\leq 6~\varpi^2\vert\varpi_-\vert^{-2}~e^{-\lambda t}.
\end{eqnarray*}
On the other hand, we have
\begin{eqnarray*}
\sigma_1^2\left[\Phi^{\epsilon}_{s}(x)\right]&=&\sigma_1^2\left[\phi_{s}(x)\right]+\left[\Phi^{\epsilon}_{s}(x)-\phi_{s}(x)\right]~\overline{\Sigma}\left[\Phi^{\epsilon}_{s}(x),\phi_{s}(x)\right],\\
\left(\partial ^2\phi_{t-s}\right)(\Phi_{s}^{\epsilon}(x))&=&\left(\partial ^2\phi_{t-s}\right)(\phi_{s}(x))+\left[\Phi^{\epsilon}_{s}(x)-\phi_{s}(x)\right]~
\overline{\Xi}_{t-s}\left[\Phi^{\epsilon}_{s}(x),\phi_{s}(x)\right].
\end{eqnarray*}
This yields the decomposition
\begin{eqnarray*}
2\epsilon^{-1}~\left[\BB^{\epsilon}_t(x)-\WW_t(x)\right]&=&\int_0^t~\sigma_1^2\left[\phi_{s}(x)\right]~\VV^{\epsilon}_s(x)~\overline{\Xi}_{t-s}\left[\Phi^{\epsilon}_{s}(x),\phi_{s}(x)\right] ds\\
&&\hskip.5cm+\int_0^t~\left(\partial ^2\phi_{t-s}\right)(\phi_{s}(x))~\VV^{\epsilon}_s(x)~\overline{\Sigma}\left[\Phi^{\epsilon}_{s}(x),\phi_{s}(x)\right]~ds\\
&&\hskip1.5cm+\epsilon~\int_0^t~\overline{\Sigma}\left[\Phi^{\epsilon}_{s}(x),\phi_{s}(x)\right]~\VV^{\epsilon}_s(x)^2~\overline{\Xi}_{t-s}\left[\Phi^{\epsilon}_{s}(x),\phi_{s}(x)\right]~ds.
\end{eqnarray*}
Taking the expectations on both sides yields
\begin{eqnarray*}
\epsilon^{-1}~\left\vert\EE\left(\BB^{\epsilon}_t(x)\right)-\WW_t(x)\right\vert&\leq &\frac{3}{\lambda }~\varpi^2\vert\varpi_-\vert^{-2}
\sigma_1^2\left[\phi_{\star}(x)\right]~~v_{1,\epsilon}(x)~\\
&&\hskip.1cm+\frac{1}{\vert\varpi_-\vert}~\varpi^2~v_{\epsilon,2}(x)~\int_0^t~~e^{-\lambda  (t-s)}~\EE\left(
~\overline{\Sigma}\left[\Phi^{\epsilon}_{s}(x),\phi_{s}(x)\right]^2\right)^{1/2}~ds\\
&&\hskip.3cm+\epsilon~3~\varpi^2~\vert\varpi_-\vert^{-2}~v_{\epsilon,4}(x)^2~\int_0^t~\EE\left(\overline{\Sigma}\left[\Phi^{\epsilon}_{s}(x),\phi_{s}(x)\right]^2\right)^{1/2}~e^{-\lambda  (t-s)} ~ds.
\end{eqnarray*}
On the other hand, we have
$$
\sup_{t\geq 0}{\EE\left(\overline{\Sigma}\left[\Phi^{\epsilon}_{t}(x),\phi_{t}(x)\right]^2\right)^{1/2}}\leq
U+3V~\left(\phi_{\star}(x)+\phi^{(\epsilon,4)}_{\star}(x)\right)^2.
$$
This yields the estimate
$$
\begin{array}{l}
\epsilon^{-1}~\left\vert\EE\left(\BB^{\epsilon}_t(x)\right)-\WW_t(x)\right\vert\\
\\
\displaystyle\leq \frac{3}{\lambda }~\left({\varpi}/{\vert\varpi_-\vert}\right)^2\\
\\
\displaystyle\times\left(\sigma_1^2\left[\phi_{\star}(x)\right]~~v_{1,\epsilon}(x)~+~\left[U+3V~\left(\phi_{\star}(x)+\phi^{(\epsilon,4)}_{\star}(x)\right)^2\right]~\left[
\frac{\vert\varpi_-\vert}{3} ~v_{\epsilon,2}(x)~+\epsilon~v_{\epsilon,4}(x)^2\right]\right).
\end{array}
$$
\qed

\section{Proofs of Riccati Stability Theorems}\label{section-contraction}

This section is mainly concerned with the proof of the transition semigroup estimate stated in Theorem~\ref{theo-expo-sg-estimate} and the two contraction theorems, Theorem~\ref{theo-2-inside} and Theorem~\ref{theo-3-intro}. The proof of the Lipschitz-type stability estimate in Theorem \ref{theo-4-intro} is also given. We start with a brief review on exponential changes of probability measures.

\subsection{Some Changes of Measure}

\begin{defi}
Let $\La^b$ be the generator of a diffusion process $\Xa_t^{b}$ on some interval $I\subset \RR$ with a
drift function $b$ and diffusion function $\sigma$; that is
\begin{equation}\label{def-Lb}
\La^b=b~\partial +\frac{\sigma^2}{2}~\partial^2.
\end{equation}
\end{defi}
For any bounded measurable function $F$ on $C([0,t],I)$ and time horizon $t\in\RR_+$, and any smooth non-negative function $h$ on $I$ we have
$$
\EE\left(F\left(\Xa^b_{[0,t]}\right)~\frac{h(\Xa^b_t)}{h(\Xa^b_0)}~\exp{\left[-\int_0^t(h^{-1}\La^b(h))(\Xa^b_s)~ds\right]}\right)=\EE\left[F\left(\Xa^{b^h}_{[0,t]}\right)\right],
$$
where $\Xa_t^{b^h}$ is a diffusion with generator
$$
\La^{b^h}(f)=\La^b(f)+h^{-1}\Gamma_{\La^b}(h,f)=b^h~\partial f+\frac{\sigma^2}{2}~\partial ^2f\quad
\mbox{\rm with}\quad b^h:=b+\sigma^2 \partial (\log{h}).
$$
Inversely, we have
\begin{equation}\label{change-expo}
\EE\left(F\left(\Xa^a_{[0,t]}\right)\right)=
\EE\left[F\left(\Xa^{b}_{[0,t]}\right)~\frac{h(\Xa^{b}_t)}{h(\Xa^{b}_0)}~\exp{\left(-\int_0^t\left[h^{-1}\La^a(h)-\left[\sigma~\partial (\log{h})\right]^2\right]\left(\Xa^{b}_s\right)~ds\right)}\right],
\end{equation}
with the drift function
$$
\begin{array}{l}
\displaystyle b=a-\sigma^2 \partial (\log{h})\\
\Longrightarrow \begin{array}[t]{rcl}
h^{-1}\La^b(h)&=&\displaystyle h^{-1}\La^a(h)- \left[\sigma~\partial (\log{h})\right]^2=\displaystyle\sigma^2~\left[\frac{1}{2}~\frac{\partial ^2h}{h}-\left[\frac{\partial h}{h}\right]^2\right]~+a~\frac{\partial h}{h}.
\end{array}\end{array}
$$
We summarize the above discussing with the following lemma.
\begin{lem}
For any bounded measurable function $F$ on $C([0,t],I)$ and time horizon $t\in\RR_+$, any regular
potential function $\Va$, and
any smooth non-negative function $h$ on $I$, we have
\begin{equation}\label{change-expo-bis}
 \begin{array}{l}
\displaystyle\EE\left(F\left(\Xa^a_{[0,t]}\right)~\exp{\left[-\int_0^t~\Va(\Xa^a_s)~ds\right]}\right)=
\EE\left[F\left(\Xa^{b}_{[0,t]}\right)~\frac{h(\Xa^{b}_t)}{h(\Xa^{b}_0)}~\exp{\left(-\int_0^t\Ha\left(\Xa^{b}_s\right)~ds\right)}\right],
\end{array}
\end{equation}
with the potential function $\Ha$ on $I$ defined by
\begin{equation}\label{def-Ha}
\Ha=h^{-1}\La^b(h)+\Va\quad\mbox{and the drift functions}\quad
b=a-\sigma^2 ~\frac{\partial h}{h}.
\end{equation}
\end{lem}
The Feynman-Kac formula (\ref{change-expo-bis}) is valid for any pair of functions $(h,\Va)$ for which the expectation make sense.
For instance, let us assume that
the pair of functions $(h,\Va)$ is chosen so that
\begin{equation}
\Ha_{\star}:=\inf_{x\geq 0}{\Ha(x)}>0~~\Longrightarrow~~ \La^b(h)+\Va h\geq \Ha_{\star}\,h.
\label{condition-hV}
\end{equation}
In this situation the Feynman-Kac equation (\ref{change-expo-bis}) is well defined as
soon as $\EE\left[h(\Xa^{b}_t)/h(\Xa^{b}_0)\right]<\infty$. In addition, we have the estimates
$$
\exp{\left[-\int_0^t~\EE\left[\Va(\Xa^a_s)\right]~ds\right]}\leq
\displaystyle\EE\left(\exp{\left[-\int_0^t~\Va(\Xa^a_s)~ds\right]}\right)\leq
\EE\left[h(\Xa^{b}_t)/h(\Xa^{b}_0)\right]~\exp{\left[-\Ha_{\star}~t\right]}.
$$
We further assume that $I=\RR_+$ and the functions $(a,\sigma)$
 in  (\ref{change-expo-bis}) are chosen such that
$$
\sup_{x\geq 0}\partial^2a(x):= \partial^2a_{\star}<0\quad\mbox{and}\quad
\sigma\leq \sigma_{\epsilon}\quad\mbox{for some } \epsilon\in \RR_+.
$$
Also assume that the potential function $\Va$ is chosen so that
\begin{eqnarray}
\partial a+\delta^{-1}~\Va+\left[\beta~\theta+\gamma~\theta^{-1}\right]&\geq &\alpha,\label{hyp-a-prime-V}
\end{eqnarray}
for some $\delta>0$ and some parameters
$ (\alpha,\beta,\gamma)\in \RR^3$ such that
\begin{eqnarray*}
 a(0)\geq \left(\gamma+\frac{\epsilon^2}{2}~(1+\delta)~U\right)_+~&\mbox{and}&
 \frac{\vert \partial^2  a_{\star}\vert}{2}\geq \beta+
\frac{\epsilon^2}{2}~(1+\delta)~ V.
\end{eqnarray*}

\begin{lem}\label{lem-tech-aVh}

Under the assumption (\ref{hyp-a-prime-V}) the Feynman-Kac formula (\ref{change-expo-bis})  is
satisfied with the potential function $\Ha$ defined by
\begin{equation}\label{drift-ref-b}
\Ha=h^{-1}\La^b(h)+\Va \quad\mbox{with}\quad
h(x)=x^{\delta}\quad\mbox{and the drift function}\quad b=a-\sigma^2~\partial \log h.
\end{equation}
In addition, the minorisation property (\ref{condition-hV}) is satisfied with
\begin{equation}\label{minoration}
\Ha_{\star}\geq \delta~ \alpha+2~\delta~\sqrt{\left[a(0)-\gamma-\frac{\epsilon^2}{2}~(1+\delta)~U\right] \left[
\frac{\vert \partial^2  a_{\star}\vert}{2}-\beta-
\frac{\epsilon^2}{2}~(1+\delta)~ V\right]},
\end{equation}
as soon as $\epsilon$ is chosen sufficiently small so that the process $\Xa^b$ is well defined on $\RR_+$.
\end{lem}

\proof
When $h=\theta^{\delta}$ for some $\delta\geq 0$ we find that
$$
b=a-\sigma^2 h^{-1}\partial h\Longrightarrow
h^{-1}\La^b(h)=\delta~\theta^{-1}a-\frac{1}{2}~\delta(\delta+1)~\left(\theta^{-1}\sigma\right)^2.
$$
For instance, when the functions $(a,\sigma)=(\Lambda,\sigma_{\epsilon})$ are given by (\ref{mp-ref}),
we find that
$$
\Ha=\delta~\partial \Lambda+\Va+\delta\theta^{-1}~\left(R-\frac{\epsilon^2}{2}~(\delta+1)~U\right)+\delta\theta\left(S-\frac{\epsilon^2}{2}~(\delta+1)~V\right).
$$
Observe that
$
a(x)=\partial a(x)~x+a(0)-x^2\int_0^1~\partial ^2a(ux)~u~du$, and further assume that
$$
\overline{\tau}^2(x):=a(0)-x^2\int_0^1~\partial ^2a(ux)~u~du> 0.
$$
The above condition is clearly met for concave drift functions with $a(0)>0$. In this case,
\begin{eqnarray*}
\Ha= \overline{\tau}^2~\frac{\partial h}{h}+
\partial a~\frac{\theta\,\partial h}{h}+\Va-\sigma^2~
\frac{\partial h}{h}~
\left[\frac{\partial h}{h}-\frac{1}{2}~\frac{\partial ^2h}{\partial h}\right].
\end{eqnarray*}
Choosing the function $h=\theta^{\delta}$,
 we have the estimate
\begin{eqnarray*}
\Ha\geq
\delta~\frac{\overline{\tau}^2}{\theta}+\delta~\partial a+\Va-
\frac{\epsilon^2}{2}~\frac{\sigma_1^2}{\theta^2}~\delta~(1+\delta).
\end{eqnarray*}
Next, observe that
$$
 \partial^2  a_{\star}<0
\Longrightarrow
\frac{\overline{\tau}^2}{\theta}\geq \frac{a(0)}{\theta}+ \frac{\vert \partial^2  a_{\star}\vert}{2}~\theta
\quad\mbox{\rm
and}\quad
\frac{\sigma_1^2}{\theta^2}= \frac{U}{\theta}+V~\theta.
$$
This yields the estimate
$$
\Ha\geq \delta \left(\frac{a(0)}{\theta}+\frac{\vert \partial^2  a_{\star}\vert}{2}~\theta\right)+\delta~\partial a+\Va-\frac{\epsilon^2}{2}~\delta~(1+\delta)~ \left(\frac{U}{\theta}+V\theta\right).
$$
Rewriting the last inequality in a slightly different form, we get
$$
\Ha/\delta\geq \left[a(0)-\frac{\epsilon^2}{2}~(1+\delta)~U\right]\frac{1}{\theta}+ \left[\frac{\vert \partial^2  a_{\star}\vert}{2}-
\frac{\epsilon^2}{2}~(1+\delta)~ V\right]~\theta+ \partial a+\frac{1}{\delta}~\Va~.
$$
Using (\ref{hyp-a-prime-V}), one can check that
$$
\begin{array}{l}
\displaystyle\Ha/\delta\geq\alpha+
\left[a(0)-\gamma-\frac{\epsilon^2}{2}~(1+\delta)~U\right]\theta^{-1}+ \left[
\frac{\vert \partial^2  a_{\star}\vert}{2}-\beta-
\frac{\epsilon^2}{2}~(1+\delta)~ V\right]~\theta.
\end{array}
$$
This completes the proof.
\qed

We are now in a position to state and prove the main theorem of this section.

\begin{theo}\label{theo-ref-kappa}
For any $\kappa\geq 0$, we choose $\epsilon\in \RR_+$ such that
\begin{eqnarray*}
\epsilon^2~\left( 1\vee \imath_{\kappa}\right)~\overline{U}<1
&\mbox{and}&
\epsilon^2~(1+\jmath_{\kappa})~ \overline{V}<2.
\end{eqnarray*}
Then  the Feynman-Kac formula (\ref{change-expo-bis})  is met with the functions
$h=\theta^{ \imath_{\kappa}}$ and
\begin{eqnarray*}
 b&=&2A~\theta+R~\left(1-\epsilon^2~ \imath_{\kappa}~\overline{U}\right)-S~\left(1+\epsilon^2~ \imath_{\kappa}~\overline{V}\right)~\theta^2\qquad
\Va=-2\kappa~(A-S\theta),  \\
 \Ha/\kappa&=&2 \imath~A+R~\left(1+ \imath\right)
 \left(1-\frac{\epsilon^2}{2}~( \imath_{\kappa}+1)~\overline{U}\right)~\theta^{-1}+S~\left(1- \imath\right)~\left(1-\frac{\epsilon^2}{2}~(1+\jmath_{\kappa})~\overline{V}\right)~\theta.
\end{eqnarray*}
In addition, we have the minorisation property
\begin{equation*}
\Ha_{\star}\geq 2\kappa~\displaystyle\frac{A^2+RS~\ell_{\epsilon,\kappa} }{\sqrt{A^2+RS}}\geq \kappa~\widehat{\lambda}_{\epsilon,\kappa}
\end{equation*}
with the collection of non-negative parameters $\ell_{\epsilon,\kappa}$ defined by
$$
\ell_{\epsilon,\kappa}^2:=\left[1-
\frac{\epsilon^2}{2}~(1+ \imath_{\kappa})~\overline{U}\right] \left[
1
-
\frac{\epsilon^2}{2}~(1+\jmath_{\kappa})~ \overline{V}\right].
$$
\end{theo}
\proof
Applying the above lemma to
$$
\begin{array}{l}
a(x)=2Ax+R-Sx^2=2(A-Sx)x+R+Sx^2\\
\\
\Longrightarrow a(0)=R\qquad \partial a(x)=2(A-Sx)\quad\mbox{\rm and}
\quad \partial^2  a=-2S,
\end{array}$$
and
$
\Va(x)=-2\kappa~(A-Sx)
$,
we find that
\begin{eqnarray*}
\partial a(x)+\frac{1}{\delta}~\Va(x)&=&2(A-Sx)-\frac{2\kappa}{\delta}~(A-Sx)=2A\left[1-\frac{\kappa}{\delta}\right]-2S~\left[1-\frac{\kappa}{\delta}\right]x.
\end{eqnarray*}
This shows that condition (\ref{hyp-a-prime-V}) is met with
$$
\left(\alpha,\beta,\gamma\right)=\left(
2A\left[1-\frac{\kappa}{\delta}\right],2S~\left[1-\frac{\kappa}{\delta}\right],0\right)\qquad
a(0)=R>0\quad\mbox{\rm and}\quad \frac{\vert \partial^2  a_{\star}\vert}{2}-\beta=S~\left(\frac{2\kappa}{\delta}-1\right)>0,
$$
for any $0\leq \delta< 2\kappa$. We set $$ \imath={\delta}/{\kappa}-1\in [-1,1]\Longleftrightarrow
\delta= \imath_{\kappa}:=
\kappa( \imath+1)\quad\mbox{\rm and}\quad (1+\jmath_{\kappa}):=\frac{1+ \imath}{1- \imath}~(1+ \imath_{\kappa})~.
$$
By (\ref{minoration}) we conclude that
\begin{eqnarray}
\Ha_{\star}&\geq& \delta~2A\left[1-\frac{\kappa}{\delta}\right]+2~\delta~\sqrt{RS\left[1-
\frac{\epsilon^2}{2}~(1+\delta)~\overline{U}\right] \left[
~\left(\frac{2\kappa}{\delta}-1\right)
-
\frac{\epsilon^2}{2}~(1+\delta)~ \overline{V}\right]}\nonumber\\
&=&2\kappa~ \left\{A~ \imath+\sqrt{1- \imath ^2}~\sqrt{RS~\left[1-
\frac{\epsilon^2}{2}~(1+ \imath_{\kappa})~\overline{U}\right] \left[
1
-
\frac{\epsilon^2}{2}~(1+\jmath_{\kappa})~ \overline{V}\right]}\right\}.\label{proof-ref-exp}
\end{eqnarray}
On the other hand, for any given $\alpha\in\RR$ and $\beta\geq 0$, the maximal value $
f_{\alpha,\beta}^{\star}:=f_{\alpha,\beta}( \imath _{\alpha,\beta})$ of the function
$$
 \imath :\in [-1,1]\mapsto
f_{\alpha,\beta}( \imath ):=\alpha~ \imath +\sqrt{(1- \imath ^2)}~\beta
$$
is attained at
\begin{equation}\label{ref-iota-alpha-beta}
 \imath _{\alpha,\beta}=\frac{\alpha}{\sqrt{\alpha^2+\beta^2}}
\Longrightarrow
f_{\alpha,\beta}^{\star}=\sqrt{\alpha^2+\beta^2}.
\end{equation}
Choosing $ \imath:=\frac{A}{\sqrt{A^2+RS}}$ in (\ref{proof-ref-exp}) we find that
$$
\begin{array}{rcl}
\displaystyle (2\kappa)^{-1}\Ha_{\star}&\geq&\displaystyle\frac{A^2+ RS~\left(1-\left[1-\ell_{\epsilon,\kappa}\right]\right)}{\sqrt{A^2+RS}}\\
\\
&\geq&\displaystyle\sqrt{A^2+RS}~\left(1-\frac{RS}{A^2+RS}~\left(1-\sqrt{1-
\frac{\epsilon^2}{2}~\left[(1+ \imath_{\kappa})~\overline{U}+(1+\jmath_{\kappa})~ \overline{V}\right]}\right)\right).
\end{array}
$$
We conclude that
$$
\begin{array}{rcl}
\displaystyle\Ha_{\star}
&\geq& \displaystyle2\kappa~\sqrt{A^2+RS}~\left(1-\frac{\epsilon^2}{2}~\frac{RS}{A^2+RS}~\frac{
(1+ \imath_{\kappa})~\overline{U}+(1+\jmath_{\kappa})~ \overline{V}}{1+\sqrt{1-
\frac{\epsilon^2}{2}~\left[
(1+ \imath_{\kappa})~\overline{U}+(1+\jmath_{\kappa})~ \overline{V}\right]}}\right)\\
\\
&\geq&\displaystyle 2\kappa~\sqrt{A^2+RS}~\left(1-\frac{\epsilon^2}{2}~\zeta_{\kappa}\right),
\end{array}
$$
with the parameter
$$
\begin{array}{l}
\displaystyle(1+\jmath_{\kappa}):=\frac{1+ \imath}{1- \imath}~(1+ \imath_{\kappa})=
\frac{ \imath_1^2}{1- \imath^2}~(1+ \imath_{\kappa})= \imath_1^2(1+ \imath_{\kappa})(1+\jmath^2)\\
\\
\displaystyle~\Longrightarrow
\frac{(1+ \imath_{\kappa})~SU+(1+\jmath_{\kappa})~ RV}{A^2+RS}=(1+ \imath_{\kappa})~\left[
\frac{1}{1+\jmath^2}~~
\overline{U}+ \imath_1^2~ \overline{V}\right].
\end{array}
$$
Finally, to complete the proof, observe that
$$
b(x)=a(x)-\epsilon^2~ \imath_{\kappa}~(U+Vx^2)=2Ax+\left(R-\epsilon^2~ \imath_{\kappa}~U\right)-\left(S+\epsilon^2~ \imath_{\kappa}~V\right)~x^2.
$$
 \qed

\subsection{Proof of Theorem~\ref{theo-2-inside}}\label{theo-2-intro-proof}

Combining (\ref{tangent-epsilon}) with (\ref{derivative-sg}) we have
\begin{equation}\label{ref-contraction}
\sigma_{1}(x)~\left\vert
\partial P_t^{\epsilon}(f)(x)\left\vert\leq \exp{\left(-\lambda_{\epsilon}\,t\right)}~
\EE\left[\sigma_{1}(\Phi_{t}^{\epsilon}(x))~\right\vert \partial f(\Phi_{t}^{\epsilon}(x))\right\vert\right].
\end{equation}
On the other hand, using (\ref{ref-contraction})
for any $x,y\in \RR_+$  and any $f\in \mbox{\rm Lip}_{\sigma_1}(\RR_+)$ we have
$$
 P_t^{\epsilon}(f)(x)- P_t^{\epsilon}(f)(y)
= (x-y)~\int_0^1~\partial P_t^{\epsilon}(f)(u~x+(1-u)~y)~du.
$$
This implies that
$$
\begin{array}{l}
\displaystyle\displaystyle \left\vert P_t^{\epsilon}(f)(x)- P_t^{\epsilon}(f)(y)\right\vert\\
\\
\displaystyle\leq
\exp{\left(-\lambda_{\epsilon}\,t\right)}~ \left\vert (x-y)~\int_0^1~\frac{P^{\epsilon}_t(\sigma_1~|\partial f|~)(u~x+(1-u)~y)}{\sigma_1(u~x+(1-u)~y)}~du\right\vert\\
\displaystyle\leq
\exp{\left(-\lambda_{\epsilon}\,t\right)}~\left\vert (x-y)~\int_0^1~\frac{1}{\sigma_1(u~x+(1-u)~y)}~~du\right\vert=\exp{\left(-\lambda_{\epsilon}\,t\right)}~
d_{\sigma_1}(x,y).
\end{array}
$$
The last assertion comes from the fact that $\Vert \sigma_1~\partial f\Vert\leq 1$ and for any $x\geq y$
\begin{eqnarray*}
~(x-y)~
\int_0^1~\frac{1}{\sigma_1(u~x+(1-u)~y)}~du
=~\int_0^1~\partial_u\left[\int_0^{ux+(1-u)y}~\frac{1}{\sigma_1(z)}~dz\right]~du=d_{\sigma_1}(x,y).
\end{eqnarray*}
We may then conclude that
$$
\DD_{\sigma_1}\left(\delta_xP^{\epsilon}_t,\delta_yP^{\epsilon}_t\right)\leq \exp{\left(-\lambda_{\epsilon}\,t\right)}~
d_{\sigma_1}(x,y)\Longrightarrow (\ref{Wasserstein-estimate}).
$$

In terms of the ``carr\'e du champ'' operator, we have
$$
(\ref{ref-contraction})\Longrightarrow\Gamma_L(P^{\epsilon}_t(f),P^{\epsilon}_t(f))\leq  \exp{\left(-2\lambda_{\epsilon}\,t\right)}~P_t^{\epsilon}\left(\sqrt{\Gamma_L(f,f)}\right)^2\leq \exp{\left(-2\lambda_{\epsilon}\,t\right)}~P_t^{\epsilon}\left(\Gamma_L(f,f)\right).
$$
Next, consider the interpolating path
$$
\forall s\in [0,t]\qquad p_t(s):=P_{s}^{\epsilon}\left(P_{t-s}^{\epsilon}(f)^2\right)\quad\mbox{\rm between }\quad
P_{t}^{\epsilon}(f)^2\quad\mbox{\rm and}\quad P_{t}^{\epsilon}(f^2).
$$
It can be readily checked that
\begin{eqnarray*}
\partial_sp_t(s)&=&\left(\partial_sP_{s}^{\epsilon}\right)\left(P_{t-s}^{\epsilon}(f)^2\right)-P_{s}^{\epsilon}\left(\partial_s\left[P_{t-s}^{\epsilon}(f)^2\right]\right)\\
&=&P_{s}^{\epsilon}\left(L\left[P_{t-s}^{\epsilon}(f)^2\right]\right)-2P_{s}^{\epsilon}\left(P_{t-s}^{\epsilon}(f)~L\left[P_{t-s}^{\epsilon}(f)\right]\right)\\
&=&P_{s}^{\epsilon}\left(\Gamma_{L}\left(P_{t-s}^{\epsilon}(f),P_{t-s}^{\epsilon}(f)\right)\right)\leq \exp{\left(-2\lambda_{\epsilon}\,(t-s)\right)}~P_t^{\epsilon}\left(\Gamma_L(f,f)\right).
\end{eqnarray*}
This implies that
\begin{eqnarray*}
2\lambda_{\epsilon}~\left[P_{t}^{\epsilon}(f^2)-P_{t}^{\epsilon}(f)^2\right]&\leq& P_t^{\epsilon}\left(\Gamma_L(f,f)\right)~.
\end{eqnarray*}
Integrating with $\pi_{\epsilon}$ and letting $t\rightarrow\infty$, we find the Poincar\'e inequality
$$
\mbox{\rm Var}_{\pi_{\epsilon}}(f)=\pi_{\epsilon}(f^2)-\pi_{\epsilon}(f)^2
\stackrel{t\rightarrow\infty}{\longleftarrow\!\!\!-\!\!\!-\!\!\!-\!\!\!-\!\!\!-}\pi_{\epsilon}(f^2)-\pi_{\epsilon}\left(P_{t}^{\epsilon}(f)^2\right)\leq \frac{1}{2\lambda_{\epsilon}\,}~\pi_{\epsilon}\left(\Gamma_L(f,f)\right)~.
$$
Recalling that $\pi L=0$,  we  have that
\begin{eqnarray*}
\partial_t\mbox{\rm Var}_{\pi_{\epsilon}}\left(P_t^{\epsilon}(f)\right)&=&2~\pi_{\epsilon}\left[P_t^{\epsilon}(f)~\partial P_t^{\epsilon}(f)(x)\right]\\
&=&-\pi_{\epsilon}\left(\Gamma_L(P^{\epsilon}_t(f),P^{\epsilon}_t(f))\right)\leq -2\lambda_{\epsilon}\,~\mbox{\rm Var}_{\pi_{\epsilon}}\left(P^{\epsilon}_t(f)\right).
\end{eqnarray*}
This completes the proof of Theorem~\ref{theo-2-inside}.\qed

\subsection{Proof of Theorem~\ref{theo-3-intro}}\label{theo-3-intro-proof}
We apply Lemma~\ref{lem-tech-aVh} to the drift function
$
a
$ and to the potential function $\Va(x)=-2(A-Sx)$. To this end, observe that
$$
a(0)=R+\frac{\epsilon^2}{2}~U>0\quad\mbox{\rm and}\quad
\partial ^2a_{\star}=-2\left(S-\frac{3\epsilon^2}{2}~V\right)<0.
$$
In addition, we have
\begin{eqnarray*}
\partial a(x)+\frac{1}{\delta}~\Va(x)&=&2A\left(1-\frac{1}{\delta}\right)-2\left(S\left(1-\frac{1}{\delta}\right)-\frac{3\epsilon^2}{2}~V\right)~x.
\end{eqnarray*}
This shows that (\ref{hyp-a-prime-V}) is met with
$$
\begin{array}{l}
\displaystyle\left(\alpha,\beta,\gamma\right)=\left(2A\left[1-\frac{1}{\delta}\right],2\left[S\left(1-\frac{1}{\delta}\right)-\frac{3\epsilon^2}{2}~V\right],0\right)\\
\\
\displaystyle\Longrightarrow
 \frac{\vert \partial ^2a_{\star}\vert}{2}-\beta=S\left(
 \frac{2}{\delta}-1\right)+\frac{3\epsilon^2}{2}~V>0\quad\mbox{\rm for any } \delta\in [0,2[.
\end{array}$$
We conclude that  the minorisation property (\ref{minoration}) is satisfied with
\begin{eqnarray*}
\Ha_{\star}/2
&\geq& A\left[\delta-1\right]+\sqrt{\delta(2-\delta)}~\sqrt{\left[R-\frac{\epsilon^2}{2}~\delta~U\right]S}.
\end{eqnarray*}
Rewriting the last inequality in a slightly different way with $ \imath :=\delta-1\in [-1,1]$,  we have
$$
\Ha_{\star}/2\geq A~ \imath +\sqrt{1- \imath ^2}~\sqrt{RS~\left[1-\frac{\epsilon^2}{2}~(1+ \imath) ~\overline{U}\right]}.
$$
Arguing as in (\ref{ref-iota-alpha-beta}) we choose
$$
 \imath =\frac{A}{\sqrt{A^2+RS}}\quad\mbox{\rm and}\quad 2R>\epsilon^2~U~ \imath_1,
$$
we find that
$$
\begin{array}{l}
\displaystyle\Ha_{\star}
\displaystyle \geq 2~\sqrt{A^2+RS}~\left(1-\frac{\epsilon^2}{2}~\zeta\right)
\end{array}
\quad\mbox{\rm with}\quad \zeta= \imath_1~\frac{SU}{A^2+RS}.
$$
On the other hand, the drift function (\ref{drift-ref-b}) is given by
\begin{eqnarray*}
b(x)&=&a(x)-\epsilon^2~ \imath_1~(U+Vx^2)\\
&=&2Ax+\left(R+\epsilon^2~\left(\frac{1}{2}- \imath_1\right)~U\right)-\left(S+\epsilon^2\left[ \imath_1-\frac{3}{2}~\right]~V\right)~x^2=2Ax+\widehat{R}_{\epsilon}-\widehat{S}_{\epsilon}~x^2.
\end{eqnarray*}
The potential function in (\ref{drift-ref-b}) is given by
$$
\Ha(x)=2 \imath  A+\frac{1+ \imath }{x}~\left[R-\frac{\epsilon^2}{2}~(1+ \imath )~U\right]+(1- \imath )~x\left[
S+\frac{\epsilon^2}{2}~(1+ \imath )~V\right]=\widehat{\Ha}_{\epsilon}(x).
$$
This ends the proof of the theorem.\qed

\subsection{Proof of Theorem~\ref{theo-expo-sg-estimate}}\label{theo-expo-sg-estimate-proof}

Combining (\ref{deterministic-expo-A-PS}) with Jensen's inequality, for any $\kappa>0$ we have
\begin{eqnarray*}
\EE\left(\Ea_{t}(x)^{2\kappa}\right)^{1/\kappa}&\geq&  \tau_t(x)\geq   \varpi(x)^2~e^{-\lambda \,t}.
\end{eqnarray*}
This proves the l.h.s. estimate in (\ref{unif-laplace-estimates-bis}). Next, set
$$
n_{\kappa}:=m_{\kappa}-( \imath_{\kappa}-1)_+=
\imath_{\kappa}-\left(\imath_{\kappa}-1\right)_-\leq
m_{\kappa}:=2\imath_{\kappa}-1.
$$
We assume  that $\epsilon$ is chosen so that
$${R}^{(\epsilon,-n_{\kappa})}\wedge{S}^{(\epsilon,-n_{\kappa})}\wedge
{R}^{(\epsilon,-m_{\kappa})}\wedge{S}^{(\epsilon,-m_{\kappa})}> 0,
$$
so that
the Riccati semigroups $\phi_t^{(\epsilon,-n_\kappa)}(x)$ and $\phi_t^{(\epsilon,-m_\kappa)}(x)$ are well defined.
We also let $\Phi^{(\epsilon,-m_\kappa)}_t(x)$ be the stochastic Riccati flow associated with the parameters
$\left({R}^{(\epsilon,-m_{\kappa})},{S}^{(\epsilon,-m_{\kappa})}\right)$.
By the exponential change of probability measure discussed in Theorem~\ref{theo-ref-kappa}, we find that
$$
\displaystyle
\EE\left(\Ea_{t}^{\epsilon}(x)^{2\kappa}\right)\\
\\
\leq \EE\left[
\left(
{{\Phi}^{(\epsilon,-m_\kappa)}_{t}(x)}/{x}
\right)^{  \imath_{\kappa}}\right]
\exp{\left[-\kappa\widehat{\lambda}_{\epsilon,\kappa} ~t\right]}.
$$
When $ \imath_{\kappa}\in \left[0,1\right]\cap[0,2\kappa]$, applying Jensen's inequality we check that
$$
 \EE\left[
\left(
{{\Phi}^{(\epsilon,-m_\kappa)}_{t}(x)}/{x}
\right)^{ \imath_{\kappa}}\right]\leq  \EE\left[
\left(
{{\Phi}^{(\epsilon,-m_\kappa)}_{t}(x)}/{x}
\right)\right]^{  \imath_{\kappa}}\leq\left[
{{\phi}^{(\epsilon,-n_\kappa)}_{\star}(x)}/{x}
\right]^{  \imath_{\kappa}}.
$$
On the other hand when $ \imath_{\kappa}\in \left[1,\infty\right]\cap[0,2\kappa]$, using (\ref{bias-n}) we
check that
$$
\begin{array}{l}
\left({R}^{(\epsilon,-m_\kappa)},{S}^{(\epsilon,-m_\kappa)}\right)+\frac{\epsilon^2}{2}~( \imath_{\kappa}-1)~\left(U,-V\right)
=({R}^{(\epsilon,-\imath_\kappa)},{S}^{(\epsilon,-\imath_\kappa)})\\
\\
\Longrightarrow \EE\left[
\left(
{{\Phi}^{(\epsilon,-m_\kappa)}_{t}(x)}/{x}
\right)^{  \imath_{\kappa}}\right]\leq
\left[
{{\phi}^{(\epsilon,-n_\kappa)}_{\star}(x)}/{x}
\right]^{  \imath_{\kappa}}.
\end{array}
$$
This implies that
\begin{eqnarray*}
\EE\left(\Ea_{t}^{\epsilon}(x)^{2\kappa}\right)&\leq&
\left[
\frac{{\varpi}^{(\epsilon,-n_\kappa)}_+}{x}\vee 1
\right]^{  \imath_{\kappa}}~\exp{\left[-\kappa~\widehat{\lambda}_{\epsilon,\kappa} ~t\right]}\leq
\left[
\frac{\varpi_+}{x}\vee 1
\right]^{  \imath_{\kappa}}~\exp{\left[-\kappa~\widehat{\lambda}_{\epsilon,\kappa} ~t\right]}.
\end{eqnarray*}
For small values of $t$ we also have
\begin{equation}\label{pre-estimate-st}
\EE\left(\Ea_{t}^{\epsilon}(x)^{2\kappa}\right)^{1/\kappa}\leq \exp{(2 A t)}.
\end{equation}
More generally, for any $0\leq s\leq t$ we have
$$
\displaystyle
\EE\left(\Ea_{s,t}^{\epsilon}(x)^{2\kappa}\right)\\
\\
\leq  \exp{\left[-\kappa~\widehat{\lambda}_{\epsilon,\kappa} ~(t-s)\right]}~\EE\left[
\left(
{{\phi}^{(\epsilon,-n_\kappa)}_{t-s}(\Phi^{\epsilon}_s(x))}/{\Phi^{\epsilon}_s(x)}
\right)^{  \imath_{\kappa}}\right].
$$
When $ \imath_{\kappa}\in \left[0,1\right]\cap[0,2\kappa]$, we check that
$$
\displaystyle
\EE\left[
\left(
{{\phi}^{(\epsilon,-n_\kappa)}_{t-s}(\Phi^{\epsilon}_s(x))}/{\Phi^{\epsilon}_s(x)}
\right)^{  \imath_{\kappa}}\right]
\\
\\
\leq   \left(\EE\left[
{{\phi}^{(\epsilon,-n_\kappa)}_{\star}(\Phi^{\epsilon}_s(x))}/{\Phi^{\epsilon}_s(x)}
\right]\right)^{  \imath_{\kappa}}.
$$
On the other hand, using the uniform estimate (\ref{unif-phi-star-2}) we have
$$
\frac{{\phi}^{(\epsilon,-n_\kappa)}_{\star}(\Phi^{\epsilon}_s(x))}{\Phi^{\epsilon}_s(x)}\leq  1+ \frac{\varpi^{(\epsilon,-n_\kappa)}_+}{\Phi^{\epsilon}_s(x)}\leq 1+ \frac{\varpi_+}{\Phi^{\epsilon}_s(x)}.
$$
Combining the above estimate with  (\ref{bias-n})  we conclude that
$$
\begin{array}{l}
\displaystyle
\EE\left[
\left(
{{\phi}^{(\epsilon,-n_\kappa)}_{t-s}(\Phi^{\epsilon}_s(x))}/{\Phi^{\epsilon}_s(x)}
\right)^{  \imath_{\kappa}}\right]
\displaystyle\leq
  \left[ 1+\varpi_+/\,{\phi}^{(\epsilon,-1)}_{s}(0)
\right]^{  \imath_{\kappa}}.
\end{array}
$$
Arguing as above, when $ \imath_{\kappa}\in \left[1,\infty\right]\cap[0,2\kappa]$ we check that
$$
  \EE\left[
\left(
{{\phi}^{(\epsilon,-n_\kappa)}_{t-s}(\Phi^{\epsilon}_s(x))}/{\Phi^{\epsilon}_s(x)}
\right)^{  \imath_{\kappa}}\right]
\le
\left[
 1
+\varpi_+/\,{\phi}^{(\epsilon,- \imath_{\kappa})}_{s}(0)\right]^{ \imath_{\kappa}}.
$$
This implies that
\begin{equation}\label{unif-laplace-estimates}
\begin{array}{l}
\displaystyle
\sup_{x\geq 0}{\EE\left(\Ea_{s,t}^{\epsilon}(x)^{2\kappa}\right)^{1/\kappa}}
\displaystyle\leq~\rho_{\epsilon,\kappa,s}~
\exp{\left[-\widehat{\lambda}_{\epsilon,\kappa} ~(t-s)\right]},
\end{array}
\end{equation}
with the constant
$$
\rho_{\epsilon,\kappa,s}:=\left[
1
+\varpi_+/~{\phi}^{\,(\epsilon,-( \imath_{\kappa}\vee 1))}_{s}(0)\right]
^{ \imath_1}.
$$
Combining H\"older inequality with  (\ref{pre-estimate-st}) and (\ref{unif-laplace-estimates}),
we  have
\begin{eqnarray*}
\displaystyle
\EE\left(\Ea_{t}^{\epsilon}(x)^{2\kappa}\right)
&\leq &
\EE\left(\Ea_{\upsilon}^{\epsilon}(x)^{2\kappa (1+1/u)}\right)^{u/(1+u)}\times \EE\left(
\Ea_{\upsilon,t}^{\epsilon}(x)^{2\kappa (1+u)}
\right)^{1/(1+u)}\\
&\leq& e^{2\kappa A \upsilon}~
\left[
1
+\varpi_+/~{\phi}^{\,(\epsilon,-( \imath _{\kappa (1+u)}\vee 1))}_{\upsilon}(0)\right]
^{\kappa(1+ \imath)}
\exp{\left[-\kappa~\widehat{\lambda}_{\epsilon,\kappa (1+u)}~(t-\upsilon)\right]},
\end{eqnarray*}
for any $u>0$.
When $0\leq t\leq \upsilon$ we also have
$$
\EE\left(\Ea_{t}^{\epsilon}(x)^{2\kappa}\right)^{1/\kappa}\leq e^{\upsilon \left[A+\widehat{\lambda}_{\epsilon,\kappa (1+u)}\right]}~\exp{\left[-\widehat{\lambda}_{\epsilon,\kappa (1+u)}~t\right]}.
$$
We end the proof of (\ref{unif-laplace-estimates-bis}) by letting $u\rightarrow0$ and choosing $\upsilon=1$. The proof of the theorem is now completed.
\qed

\subsection{Proof of Theorem~\ref{theo-4-intro}}\label{theo-4-intro-proof}

Combining the Feynman-Kac formula (\ref{tangent-overline}) with the estimate (\ref{tangent-overline-estimate})
 for any $n\geq 1$ and taking any $x_1>x_2$ we have the Taylor integral formula
 \begin{eqnarray*}
 \Phi^{\epsilon}_t(x_1)-\Phi^{\epsilon}_t(x_2)&=&(x_1-x_2)~\int_0^1~\Ta^{\epsilon}_t(ux_1+(1-u)x_2)~du,
 \end{eqnarray*}
as well as the estimates
 $$
\begin{array}{l}
\vertiii{\Phi^{\epsilon}_t(x_1)-\Phi^{\epsilon}_t(x_2)}_n\\
\\
 \displaystyle\leq ~\exp{\left[-\widehat{\lambda}_{\epsilon} \,t\right]}~~\int_0^1~\frac{(x_1-x_2)
}{\widehat{\sigma}(ux_1+(1-u)x_2)}~ \vertiii{\widehat{\Phi}^{\epsilon}_t(ux_1+(1-u)x_2)}_{ \imath_n}^{ \imath_1}~du.
\end{array}
$$
On the other hand, we have
$$
\vertiii{\widehat{\Phi}^{\epsilon}_t(ux_1+(1-u)x_2)}_{ \imath_n}\leq \widehat{\phi}^{(\epsilon, \imath_n)}_t(ux_1+(1-u)x_2)\leq \widehat{\phi}^{(\epsilon, \imath_n)}_t(x_1).
$$
Using (\ref{unif-phi-star}) for any time $t\geq \upsilon>0$ we  check that
$$
\begin{array}{l}
\vertiii{\widehat{\Phi}^{\epsilon}_t(ux_1+(1-u)x_2)}_{ \imath_n}\leq ~ c_{1,\upsilon}~(2\widehat{\varpi}^{(\epsilon, \imath_n)}_{+}-\widehat{\varpi}^{(\epsilon, \imath_n)}_{-})
\\
\\
 \qquad\displaystyle\Longrightarrow \vertiii{\Phi^{\epsilon}_t(x_1)-\Phi^{\epsilon}_t(x_2)}_n\leq ~c_{2,\upsilon}~\left[2\widehat{\varpi}^{(\epsilon, \imath_n)}_{+}-\widehat{\varpi}^{(\epsilon, \imath_n)}_{-}\right]^{ \imath_1}~d_{\widehat{\sigma}}(x_1,x_2)~\exp{\left[-\widehat{\lambda}_{\epsilon} \,t\right]}~.
\end{array}$$
More generally we have
$$
\begin{array}{l}
\vertiii{\Phi^{\epsilon}_t(x_1)-\Phi^{\epsilon}_t(x_2)}_n~\leq~d_{\widehat{\sigma}}(x_1,x_2)~\left(x_1\vee x_2\vee  \widehat{\varpi}^{(\epsilon, \imath_n)}_+\right)^{ \imath_1}~\exp{\left[-\widehat{\lambda}_{\epsilon} \,t\right]}~.
\end{array}
$$
This ends the proof of the theorem. \qed

\section{Other Proofs}

In this Appendix we give the proof of Proposition \ref{FK-intro-prop}, the proof of Lemma \ref{lem-3}, the proof of Theorem \ref{theo-OU-proc}, and the proof of Corollary \ref{cor-2-EnKF} in order.

\subsection{Proof of Proposition~\ref{FK-intro-prop}}\label{FK-intro-prop-proof}

We have
$$
\EE\left[f(\Phi_{t}^{\epsilon}(x))~\Ta_t^{\epsilon}(x)\right]=
\EE\left[f(\Phi_{t}^{\epsilon}(x))~\exp{\left[\int_0^t\partial \Lambda\left(\Phi^{\epsilon}_s(x)\right)~ds\right]}~M^{\epsilon}_t(x)\right],
$$
with the exponential martingale $M^{\epsilon}_t(x)$ defined by
$$
\sigma_{\epsilon}(x)~M^{\epsilon}_t(x):=\sigma_{\epsilon}(\Phi_{t}^{\epsilon}(x))~\exp{\left[-\int_0^t~\left(\sigma^{-1}_{\epsilon}L\sigma_{\epsilon}\right)(\Phi^{\epsilon}_s(x))~ds\right]}.
$$
This implies that
\begin{eqnarray*}
\EE\left[f(\Phi_{t}^{\epsilon}(x))~\Ta_t^{\epsilon}(x)\right]&=&\EE\left[f(\overline{\Phi}_{t}^{\epsilon}(x))~\exp{\left[\int_0^t~\partial \Lambda\left(\overline{\Phi}^{\epsilon}_s(x)\right)~ds\right]}\right],
\end{eqnarray*}
where $\overline{\Phi}_{t}^{\epsilon}(x)$ stands for the stochastic flow associated with a diffusion with generator
$$
\La(f)=L(f)+\sigma_{\epsilon}^{-1}\Gamma_L\left(\sigma_{\epsilon},f\right)=L(f)+\sigma_{\epsilon}~\partial \sigma_{\epsilon}~\partial f.
$$
In other words $\overline{\Phi}_{t}^{\epsilon}(x)$ is the stochastic Riccati flow associated
with the diffusion function $\sigma_{\epsilon}$ and the drift function $a$ given by
$$
\sigma_{\epsilon}(x)~\partial \sigma_{\epsilon}(x)=\frac{\epsilon^2}{2}~ (U+3Vx^2)\Longrightarrow a(x)=2Ax+\left(R+\frac{\epsilon^2}{2}~U\right)-\left(S-\frac{3\epsilon^2}{2}~V\right)~x^2.
$$
This ends the proof of the proposition.\qed

\subsection{Proof of Lemma~\ref{lem-3}}\label{proof-lem-3}

For any $x>0$ we have
\begin{eqnarray*}
\partial \sigma_{\epsilon}(x)&=&\frac{\epsilon}{2}~\frac{U+3Vx^2}{\sqrt{x(U+Vx^2)}}\Longrightarrow \frac{\partial \sigma_{\epsilon}(x)}{\sigma_{\epsilon}(x)}=\frac{1}{2}~\frac{U+3Vx^2}{x(U+Vx^2)}.
\end{eqnarray*}
In addition, we have
$$
\begin{array}{l}
\displaystyle\partial ^2\sigma_{\epsilon}(x)=\frac{\epsilon}{2}~\frac{1}{\sqrt{x(U+Vx^2)}}
\left[6Vx-\frac{(U+3Vx^2)^2}{2x(U+Vx^2)}\right]\\
\\
\displaystyle\Longrightarrow \partial ^2\sigma_{\epsilon}(x)\sigma_{\epsilon}(x)=\frac{\epsilon^2}{2}~
\left[6Vx-\frac{(U+3Vx^2)^2}{2x(U+Vx^2)}\right].
\end{array}
$$
This implies that
$$
\Ha^{\epsilon}(x)=\frac{U+3Vx^2}{U+Vx^2}~\left(A+\frac{1}{2}~\frac{R}{x}-\frac{S}{2}~x\right)-2~(A-Sx)+\frac{\epsilon^2}{4}~
\left[6Vx-\frac{1}{x}~\frac{(U+3Vx^2)^2}{U+Vx^2}\right].
$$
For instance when $V= 0$ we have
$$
\Ha^{\epsilon}(x)=~-A+\frac{1}{2} \left[\left(R-\frac{\epsilon^2}{2}~U\right)~\frac{1}{x}+3S~x\right]\geq \lambda_{\epsilon}\,:=-A+
\sqrt{3
 \left(R-\frac{\epsilon^2}{2}~U\right)S},
$$
as soon as
$
2R\geq \epsilon^2~U
$.
When $U=0$ we have
$$
\Ha^{\epsilon}(x)=A+\frac{3R}{2}~\frac{1}{x}+x~\frac{1}{2}~\left(
S-\frac{3\epsilon^2}{2}~V~\right)\geq \lambda_{\epsilon}\,=A+\sqrt{3R\left(
S-\frac{3\epsilon^2}{2}~V~\right)},
$$
as soon as $2S\geq 3\epsilon^2V$.
We further assume that $V\wedge U>0$.
In this case, we have
\begin{eqnarray*}
\Ha^{\epsilon}(x)&=&\left[\frac{U+3Vx^2}{U+Vx^2}-2\right]~A\\
&&+\frac{1}{x}~\frac{U+3Vx^2}{U+Vx^2}~\left[\frac{R}{2}-\frac{\epsilon^2}{4}~U-\frac{3\epsilon^2}{4}~Vx^2\right]~+Sx~
\left(2-\frac{1}{2}~\frac{U+3Vx^2}{U+Vx^2}+~\frac{3}{2}~\epsilon^2~V~\right).
\end{eqnarray*}
This yields the decomposition
\begin{eqnarray*}
\Ha^{\epsilon}(x)&=&\left[ \iota (x)-2\right]~A+\frac{1}{x}~ \iota (x)~\left[\frac{R}{2}-\frac{\epsilon^2}{4}~U\right]
+x~
\left[S\left(2-\frac{1}{2}~ \iota (x)\right)+~\frac{3}{2}~\epsilon^2~V\left(S-\frac{1}{2}~ \iota (x)\right)\right],
\end{eqnarray*}
with the increasing function
$$
x\in [0,1]~\mapsto~ \iota (x)=\frac{U+3Vx^2}{U+Vx^2}=1+\frac{2Vx^2}{U+Vx^2}\in [1,3]\Longrightarrow
 \iota (x)-2=\frac{Vx^2-U}{U+Vx^2}\in [-1,1].
$$
Observe that
\begin{eqnarray*}
\Ha^{\epsilon}(x)&\geq &\frac{(\;x\,\sqrt{V/U}\,)^2-1}{1+(\,x\,\sqrt{V/U}\,)^2}~A+\frac{1}{x\,\sqrt{V/U}}~r_{\epsilon}
+s_{\epsilon}~\left(\sqrt{V/U}~x\right)~:=\Ha_{-}^{\epsilon}\left(x\sqrt{V/U}\right),
\end{eqnarray*}
with
$$
r_{\epsilon}:=\frac{1}{2}~\sqrt{V/U}~\left[R-\frac{\epsilon^2}{2}~U\right]\geq 0\quad\mbox{\rm and}\quad
s_{\epsilon}:=\frac{1}{2\sqrt{V/U}}~\left[S+3~\epsilon^2~V\left(S-\frac{3}{2}\right)\right]\geq 0.
$$
When $A=0$, for any $x\in \RR_+$ we have
$$
\Ha_{-}^{\epsilon}\left(x\right)=\frac{1}{x}~r_{\epsilon}
+s_{\epsilon}~x\geq 2~\sqrt{r_{\epsilon}s_{\epsilon}}\Longrightarrow\lambda_{\epsilon}\,\geq \sqrt{\left[R-\frac{\epsilon^2}{2}~U\right]\left[S+3~\epsilon^2~V\left(S-\frac{3}{2}\right)\right]},
$$
as soon as $2R\geq \epsilon^2~U$ and $S\geq \frac{9}{2}~\frac{\epsilon^2}{1+3\epsilon^2V}~V$.
We further assume that $A\not=0$. Observe that
\begin{eqnarray*}
\Ha_{-}^{\epsilon}\left(1/x\right)&= &\frac{x^2-1}{1+x^2}~(-A)+\frac{1}{x}~s_{\epsilon}
+r_{\epsilon}~x.
\end{eqnarray*}
This shows that there is no loss of generality in  assuming that $A>0$, up to changing  $(r_{\epsilon},s_{\epsilon})$ by $(s_{\epsilon},r_{\epsilon})$.
We further assume that $A>0$. In this case, we have
\begin{eqnarray*}
\overline{\Ha}_{\,-}^{\,\epsilon}\left(x\right)&:=&\Ha_{-}^{\epsilon}\left(x\right)/A\\
&=& \frac{x^2-1}{1+x^2}+\frac{1}{x}~\overline{r}_{\epsilon}
+\overline{s}_{\epsilon}~x\geq  \frac{x^2-1}{1+x^2}+\left(\frac{1}{x}~
+~x\right)~\left(\overline{r}_{\epsilon}\wedge \overline{s}_{\epsilon}\right)\quad\mbox{\rm with}\quad (\overline{r}_{\epsilon},\overline{s}_{\epsilon})=(r_{\epsilon},s_{\epsilon})/A.
\end{eqnarray*}
This yields the estimate
$$
\begin{array}{l}
\displaystyle\overline{\Ha}_{\,-}^{\,\epsilon}\left(x\right)
=-1+ \frac{2x^2}{1+x^2}+\frac{1}{x}~\overline{r}_{\epsilon}+\overline{s}_{\epsilon}~x\geq -1+2\sqrt{\overline{r}_{\epsilon}\overline{s}_{\epsilon}}\\
\\
\displaystyle\Longrightarrow
\lambda_{\epsilon}\,\geq -A+\sqrt{\left[R-\frac{\epsilon^2}{2}~U\right]\left[S+3~\epsilon^2~V\left(S-\frac{3}{2}\right)\right]}.
\end{array}
$$
This completes the proof of the lemma.
\qed

\subsection{Proof of Theorem~\ref{theo-OU-proc}}\label{theo-OU-proc-proof}

The first estimate  is a direct consequence of the Ornstein-Uhlenbeck formula (\ref{OU-formula}) combined with
the exponential semigroup estimates stated in Theorem~\ref{theo-expo-sg-estimate}. It is  also readily checked  that
$$
\sup_{\overline{\epsilon}\in [0,1]}{\vertiii{\Psi^{(\epsilon,\overline{\epsilon})}_t(x,z_1)-\Psi^{(\epsilon,\overline{\epsilon})}_t(x,z_2)}_n}\leq c_{1,n}~
\exp{\left[-\widehat{\lambda}_{\epsilon, {n/2}} ~t/2\right]}~\vert z_1-z_2\vert.
$$
We fix the parameters $(z,\epsilon,\overline{\epsilon})$  and we set
\begin{eqnarray*}
\Delta_t^{\psi}(x_1,x_2)&:=&\Psi^{(\epsilon,\overline{\epsilon})}_t(x_1,z)-\Psi^{(\epsilon,\overline{\epsilon})}_t(x_2,z)\qquad
\Delta_t^{\phi}(x_1,x_2):=\Phi^{\epsilon}_t(x_1)-\Phi^{\epsilon}_t(x_2),\\
\Delta_t(x_1,x_2)&:=&\frac{S\left[\Phi^{\epsilon}_t(x_1)+\Phi^{\epsilon}_t(x_2)\right]+\overline{\epsilon}^2
\left[U+V\left(\Phi^{\epsilon}_t(x_1)^2+\Phi^{\epsilon}_t(x_1)\Phi^{\epsilon}_t(x_2)+\Phi^{\epsilon}_t(x_2)^2\right)\right]}{\varsigma_{\overline{\epsilon}}\left(\Phi^{\epsilon}_t(x_1)\right)+\varsigma_{\overline{\epsilon}}\left(\Phi^{\epsilon}_t(x_2)\right)}\\
&\leq &\frac{R^{-1/2}}{2}~\left[S\left[\Phi^{\epsilon}_t(x_1)+\Phi^{\epsilon}_t(x_2)\right]+\overline{\epsilon}^2
\left[U+V\left[\Phi^{\epsilon}_t(x_1)+\Phi^{\epsilon}_t(x_2)\right]^2\right]\right].
\end{eqnarray*}
Using the norm estimates (\ref{bias-n}) we check that
\begin{eqnarray*}
\vertiii{\Delta(x_1,x_2)}_n&\leq&
\frac{R^{-1/2}}{2}~\left[S\left[\phi^{(\epsilon,n)}_{\star}(x_1)+\phi^{(\epsilon,n)}_{\star}(x_2)\right]+\overline{\epsilon}^2
\left[U+V\left[\phi^{(\epsilon,2n)}_{\star}(x_1)+\phi^{(\epsilon,2n)}_{\star}(x_2)\right]^2\right]\right]\\
&\leq &c_{2,n}~(1+x_1+x_2)^2.
\end{eqnarray*}
On the other hand combining the Laplace estimates (\ref{unif-laplace-estimates-bis}) with the decomposition (\ref{OU-formula})
and Burkholder-Davis-Gundy inequality we check that
$$
\vertiii{\Psi^{(\epsilon,\overline{\epsilon})}_t(x,z)}_n\leq c_{3,n}
\exp{\left[-2\widehat{\lambda}_{\epsilon, {n/2}} ~t\right]}~\vert z\vert+n~\EE\left[\left\vert\int_0^t~
\Ea_{s,t}^{\epsilon}(x)^2~\varsigma^2_{\overline{\epsilon}}(\Phi^{\epsilon}_s(x))~ds\right\vert^{n/2}\right]^{1/n}.
$$
Using the generalized Minkowski inequality this implies that
$$
\begin{array}{l}
\displaystyle\vertiii{\Psi^{(\epsilon,\overline{\epsilon})}_t(x,z)}_n\\
\\
\displaystyle\leq c_{3,n}
\exp{\left[-2\widehat{\lambda}_{\epsilon, {n/2}} ~t\right]}~\vert z\vert+n~\left[\int_0^t~
\EE\left[\Ea_{s,t}^{\epsilon}(x)^{2n}\right]^{1/n}
\EE\left[\varsigma^{2n}_{\overline{\epsilon}}(\Phi^{\epsilon}_s(x))\right]^{1/n}
~ds\right]^{1/2}\\
\\
\displaystyle\leq  c_{3,n}
\exp{\left[-2\widehat{\lambda}_{\epsilon, {n/2}} ~t\right]}~\vert z\vert+c_{4,n}~~\left(\int_0^t~
\exp{\left[-\widehat{\lambda}_{\epsilon, n} ~(t-s)\right]}~
\EE\left[\varsigma^{2n}_{\overline{\epsilon}}(\Phi^{\epsilon}_s(x))\right]^{1/n}
~ds\right)^{1/2}.
\end{array}
$$
Next, observe that
\begin{eqnarray*}
\EE\left[\varsigma^{2n}_{\overline{\epsilon}}(\Phi^{\epsilon}_s(x))\right]^{1/n}&\leq& \left(R+S~
\phi^{(\epsilon,2n)}_{\star}(x)^{2}\right)+\overline{\epsilon}^2~ \EE\left[\left(
\Phi^{\epsilon}_s(x)(U+V\Phi^{\epsilon}_s(x)^2)\right)^n\right]^{1/n}\\
&\leq &\left(R+S~
\phi^{(\epsilon,2n)}_{\star}(x)^{2}\right)+\overline{\epsilon}^2~\left(U\phi^{(\epsilon,n)}_{\star}(x)+V\phi^{(\epsilon,3n)}_{\star}(x)^3\right).
\end{eqnarray*}
We may then conclude that
$$
\begin{array}{l}
\displaystyle\vertiii{\Psi^{(\epsilon,\overline{\epsilon})}_t(x,z)}_n
\leq
c_{5,n}
\exp{\left[-2\widehat{\lambda}_{\epsilon,n/2}~t\right]}~\vert z\vert\\
\\
\displaystyle\hskip2cm+c_{6,n}~\left(\left(R+S~
\phi^{(\epsilon,2n)}_{\star}(x)^{2}\right)+\overline{\epsilon}^2~\left(U\phi^{(\epsilon,n)}_{\star}(x)+V\phi^{(\epsilon,3n)}_{\star}(x)^3\right)\right)^{1/2}.
\end{array}
$$
This ends the proof of (\ref{first-unif-Psi-formula}).
The proof of the second estimate in (\ref{first-contraction-formula}) is based on the formula
$$
\begin{array}{l}
\displaystyle
d\Delta_t^{\psi}(x_1,x_2)\\
\\
=\displaystyle\frac{1}{2}~\partial \Lambda\left(\Phi^{\epsilon}_t(x_1)\right)~\Delta_t^{\psi}(x_1,x_2)~dt-S~
\Delta_t^{\phi}(x_1,x_2)~\Psi^{(\epsilon,\overline{\epsilon})}_t(x_2,z)~dt
+\displaystyle\Delta_t^{\phi}(x_1,x_2)~\Delta_t(x_1,x_2)~d\Wa^{\prime}_t.
\end{array}
$$
This implies that
$$
\begin{array}{l}
\displaystyle
\Delta_t^{\psi}(x_1,x_2)
\\
=\displaystyle-S~\int_0^t~\Ea_{s,t}(x_1)~
\Delta_s^{\phi}(x_1,x_2)~\Psi^{(\epsilon,\overline{\epsilon})}_s(x_2,z)~ds
+\displaystyle \int_0^t~\Ea_{s,t}(x_1)~\Delta_s^{\phi}(x_1,x_2)~\Delta_s(x_1,x_2)~d\Wa^{\prime}_s,
\end{array}
$$
from which we check that
$$
\begin{array}{l}
\displaystyle
\vertiii{\Delta_t^{\psi}(x_1,x_2)}_n
\leq \displaystyle c_{1,n}~(1+\vert z\vert+x_2^{3/2})~\int_0^t~\EE\left(\Ea_{s,t}(x_1)^{3n}\right)^{1/(3n)}
\EE\left(\Delta_s^{\phi}(x_1,x_2)^{3n}\right)^{1/(3n)}~ds\\
\\
\hskip2cm+\displaystyle c_{2,n}~(1+x_1+x_2)^2~\left[\int_0^t~
\EE\left[\Ea_{s,t}(x_1)^{3n}\right]^{2/(3n)}
~
\EE\left[\Delta_s^{\phi}(x_1,x_2)^{3n}\right]^{2/(3n)}~ds\right]^{1/2}.
\end{array}
$$
Using (\ref{unif-laplace-estimates-bis}) and
(\ref{Wasserstein-estimate-overline-bis}), we obtain that
$$
\begin{array}{l}
\displaystyle
\vertiii{\Delta_t^{\psi}(x_1,x_2)}_n\\
\\
\leq \displaystyle c_{1,n}~(1+\vert z\vert+x_2^{3/2})~\frac{1+x_1+x_2}{x_1\wedge x_2}~\vert x_1-x_2\vert~
~\int_0^t~ \exp{\left(
 -\widehat{\lambda}_{\epsilon} \,s\right)}
\exp{\left[-\widehat{\lambda}_{\epsilon,3n/2}~(t-s)/2\right]}
~~ds\\
\\
\hskip1cm+\displaystyle c_{2,n}~\frac{(1+x_1+x_2)^3}{x_1\wedge x_2}~\vert x_1-x_2\vert~\left[\int_0^t~
 \exp{\left(
 -2\widehat{\lambda}_{\epsilon} \,s\right)}~
\exp{\left[-\widehat{\lambda}_{\epsilon,3n/2}~(t-s)\right]}~ds\right]^{1/2}.
\end{array}
$$
As a result, we conclude that
$$
\displaystyle
\vertiii{\Delta_t^{\psi}(x_1,x_2)}_n
\leq \displaystyle c_{3,n}~(1+\vert z\vert)~\frac{(1+x_1+x_2)^3}{x_1\wedge x_2}~\vert x_1-x_2\vert~t~
\exp{\left[-\widehat{\lambda}_{\epsilon,3n/2}~t/2\right]},
$$
completing the proof of the theorem.
\qed

\subsection{Proof of Corollary~\ref{cor-2-EnKF}}\label{cor-2-EnKF-proof}

The proof is based on the decomposition
\begin{eqnarray*}
\Psi^{(\epsilon,0)}_t(x,z)-\Psi^{(0,0)}_t(x,z)&=&
\left[\Ea_{t}^{\epsilon}(x)-\Ea_{t}(x)\right]~z+\int_0^t~
\left(\Ea_{s,t}^{\epsilon}(x)-\Ea_{s,t}(x)\right)~\varsigma(\Phi^{\epsilon}_s(x))~d\Wa^{\prime}_s\\
&&
\hskip1cm\displaystyle+\epsilon~S\int_0^t~\Ea_{s,t}^{\epsilon}(x)~\VV^{\epsilon}_s(x)~\frac{
\Phi^{\epsilon}_s(x)+\phi_s(x)}{\varsigma(\Phi^{\epsilon}_s(x))+\varsigma(\phi_s(x))}~d\Wa^{\prime}_s.
\end{eqnarray*}
On the other hand we have
$$
\begin{array}{l}
\displaystyle
\left\vert \Ea_{s,t}^{\epsilon}(x)-\Ea_{s,t}(x)\right\vert
\displaystyle\leq \epsilon~S~
\left[\Ea_{s,t}^{\epsilon}(x)+\Ea_{s,t}(x)\right]~\left[\int_s^t\vert \VV^{\epsilon}_u(x)\vert~du\right].
\end{array}
$$
Arguing as above and using the uniform fluctuation estimates (\ref{unif-est-intro}), we get  that
$$
\begin{array}{l}
\displaystyle
\EE\left[\left\vert \Ea_{s,t}^{\epsilon}(x)-\Ea_{s,t}^{0}(x)\right\vert^n\right]^{1/n}\\
\\
\displaystyle\leq
 \epsilon~S~v^{\epsilon}_{2n}(x)~(t-s)
\left(\EE\left(\Ea_{s,t}^{\epsilon}(x)^{2n}\right)^{1/(2n)}+\Ea_{s,t}^{0}(x)~\right)\\
\\
\displaystyle\leq c_{1,n}~
 \epsilon~v^{\epsilon}_{2n}(x)~(t-s)
\left(
\exp{\left[-\widehat{\lambda}_{\epsilon, n} ~(t-s)/2\right]}+ \exp{\left[-\lambda (t-s)/2\right]}~\right).
\end{array}
$$
The last assertion is a direct consequence of  (\ref{deterministic-expo-A-PS}) and the uniform Laplace estimate (\ref{unif-laplace-estimates-bis}).
Combining Burkholder-Davis-Gundy inequality and the generalized Minkowski inequality, we have
$$
\begin{array}{l}
\displaystyle
\EE\left[\left\vert \int_0^t~
\left(\Ea_{s,t}^{\epsilon}(x)-\Ea_{s,t}^{0}(x)\right)~\varsigma(\Phi^{\epsilon}_s(x))~d\Wa^{\prime}_s\right\vert^n\right]^{2/n}\\
\\
\displaystyle\leq n^2~\int_0^t~
\EE\left(\left(\Ea_{s,t}^{\epsilon}(x)-\Ea_{s,t}^{0}(x)\right)^{2n}\right)^{1/n}~\left[R+S~
\EE\left(\Phi^{\epsilon}_s(x)^{2n}\right)^{1/n}\right]
~ds.
\end{array}
$$
Using the uniform moment estimates (\ref{bias-n}), we find that
$$
\begin{array}{l}
\displaystyle
\EE\left[\left\vert \int_0^t~
\left(\Ea_{s,t}^{\epsilon}(x)-\Ea_{s,t}^{0}(x)\right)~\varsigma(\Phi^{\epsilon}_s(x))~d\Wa^{\prime}_s\right\vert^n\right]^{2/n}\\
\\
\displaystyle\leq c_{2,n}~\epsilon^2~v^{\epsilon}_{4n}(x)^2\left[1+~
\phi^{(\epsilon,2n)}_{\star}(x)^2\right]~\int_0^t~
 ~s^2
\left(
\exp{\left[-\widehat{\lambda}_{\epsilon,2n} ~s\right]}+\exp{\left[-\lambda  s\right]}~\right)
~ds\\
\\
\displaystyle\leq~c_{3,n}~\epsilon^2~v^{\epsilon}_{4n}(x)^2\left[1+
\phi^{(\epsilon,2n)}_{\star}(x)^2\right].
\end{array}
$$
In much the same way one can check that
$$
\begin{array}{l}
\displaystyle
\EE\left[\left\vert \int_0^t\Ea_{s,t}^{\epsilon}(x)~\VV^{\epsilon}_s(x)~\frac{
\Phi^{\epsilon}_s(x)+\phi_s(x)}{\varsigma(\Phi^{\epsilon}_s(x))+\varsigma(\phi_s(x))}~d\Wa^{\prime}_s\right\vert^n\right]^{2/n}\\
\\
\displaystyle\leq c_{4,n}~\int_0^t\EE\left[
\Ea_{s,t}^{\epsilon}(x)^{n}~\VV^{\epsilon}_s(x)^{n}~
\left[\Phi^{\epsilon}_s(x)^2+\phi_s(x)^2\right]^{n/2}\right]^{2/n}~ds\\
\\
\displaystyle\leq c_{5,n}~v^{\epsilon}_{3n}(x)^2~\left[\phi_{\star}(x)^2+\phi_{\star}^{(\epsilon,3n)}(x)^2
\right]~\int_0^t~
 ~s^2
\left(
\exp{\left[-\widehat{\lambda}_{\epsilon, 3n} ~s\right]}+\exp{\left[-\lambda  s\right]}~\right)
~ds.
\end{array}
$$
This implies that
$$
\begin{array}{l}
\displaystyle
\EE\left[\left\vert \int_0^t\Ea_{s,t}^{\epsilon}(x)~\VV^{\epsilon}_s(x)~\frac{
\Phi^{\epsilon}_s(x)+\phi_s(x)}{\varsigma(\Phi^{\epsilon}_s(x))+\varsigma(\phi_s(x))}~d\Wa^{\prime}_s\right\vert^n\right]^{2/n}
\displaystyle\leq c_{6,n}~v^{\epsilon}_{3n}(x)^2~\left[\phi_{\star}(x)^2+\phi_{\star}^{(\epsilon,3n)}(x)^2
\right].
\end{array}
$$
Recalling that $\sup_{x\geq 0}{(xe^{-\alpha x})}=1/(\alpha e)$  for any $\alpha$, we  find the uniform estimate
$$
\begin{array}{l}
\displaystyle
\epsilon^{-1}{\vertiii{\Psi^{(\epsilon,0)}(x,z)-\Psi^{(0,0)}(x,z)}_{n}}
\displaystyle\leq
c_{6,n}~v^{\epsilon}_{4n}(x)~(z\vee 1)~\left[1\vee
\phi^{(\epsilon,3n)}_{\star}(x)\right]~.
\end{array}
$$
This completes the proof of the corollary.
\qed

\end{document}